\documentclass[twoside,11pt]{article}
\usepackage[top=2.5cm,bottom=2.5cm,left=3cm,right=3cm]{geometry}
\usepackage[dvipsnames,svgnames,x11names]{xcolor}
\usepackage{amsmath,amssymb,mathtools}
\usepackage{graphicx}
\usepackage{float}
\usepackage{booktabs,tabularx,multirow,colortbl}
\usepackage{caption}
\usepackage{hyperref}
\usepackage{fancyhdr}
\usepackage{multicol}
\usepackage{enumitem}

\hypersetup{
    colorlinks=true,
    linkcolor=blue,
    citecolor=blue,
    urlcolor=blue
}

\usepackage{setspace}
\usepackage{xurl}

\begin{document}

\thispagestyle{empty}

\begin{center}
    \vspace*{1cm}
    \vspace{0.2cm}
    \vspace{0.8cm}
    
    \vspace{1cm}
    \vspace{2.5cm}
    
    \noindent
        
    
    \vspace{1.2cm}
    
    \vspace{0.2cm}
    {\large\bfseries \url{}} \\
    \vspace{0.5cm}
\end{center}

\vspace{0.8cm}

\begin{multicols}{2}
    \footnotesize
    
    \columnbreak
    
\end{multicols}

\vspace{0.4cm}
\begin{center}
    \vspace{0.2cm}
\end{center}

\thispagestyle{empty}
\setcounter{page}{1}

\begin{center}
    {\Large\bfseries Canonical metric on a mildly singular Kähler varieties with an intermediate log Kodaira dimension} \\
    \vspace{0.5cm}
    {\large \textsc{Hassan Jolany}} \\
    \vspace{1cm}
    
    \begin{minipage}{0.85\textwidth}
        \small
        \noindent \textbf{Abstract.} Existence of canonical metric on a canonical model of projective singular variety was a long standing conjecture and the major part of this conjecture is about varieties which do not have definite first Chern class(most of the varieties do not have definite first Chern class). There is a program which is known as Song-Tian program for finding canonical metric on canonical model of a projective variety by using Minimal Model Program. In this paper, we apply Song-Tian program for mildly singular pair (X, D) via Log Minimal Model Program where D is a simple normal crossing divisor on X with conic singularities. We show that there is a unique $C^\infty$-fiberwise conical Kähler-Einstein metric on (X, D) with vanishing Lelong number which is twisted by logarithmic Weil-Petersson metric and an additional term of Fujino-Mori [72] as soon as we have fiberwise KE-stability or Kawamata's condition of Theorems 2.28, or 2.30(in $C^0$ case). In final we highlight that how the complete answer of this question can be reduced to CMA equation corresponding to fiberwise Calabi-Yau foliation(due to H.Tsuji) which is still open.
        
        \vspace{0.5cm}
        \noindent \textcolor{blue}{53C44, 32Q15, ; 53C55, 58J35, 14E30,}
    \end{minipage}
\end{center}

\vspace{1cm}

\section{Introduction}

In this paper we study conical Kähler metrics. This analysis uses the language and the basic theory of Kähler currents which extensively studied by Jean-Pierre Demailly, Claire Voisin, Nessim Sibony, Tien-Cuong Dinh, Vincent Guedj, Gang Tian, and Jian Song, see [1, 2, 3, 4, 5] and references therein. Conical Kähler metrics introduced by Tian on the quasi-projective varieties X $\setminus$ D where X is a smooth projective variety and D $\subset$ X is a simple normal crossing divisor and also by Thurston in connection with geometrization conjecture. Conical Kähler metrics have recently become a key ingredient in the solution of the Tian-Yau-Donaldson conjecture about the existence of Kähler-Einstein metrics of positive first Chern class [6, 7, 8, 9]. Conical Kähler-Einstein metrics are unique, i.e. canonically attached to pair (X, D). These type of metrics can be used to probe the paired variety (X, D) using differential-geometric tools, for instance Tian's proof on logarithmic version of Miyaoka-Yau Chern number inequality.


One of the applications of conical Kähler Ricci flow on holomorphic fiber spaces is in Analytical Minimal Model Program (AMMP) which introduced by Song and Tian [10],[11] and a different approach by Tsuji [13] also by using dynamical systems of Bergman kernels and finding singular canonical measure which its inverse is Analytic Zariski Decomposition. Robert Berman in [14], also gave a new statistical-mechanical approach to the study of canonical metrics and measures on a complex algebraic variety $X$ with positive Kodaira dimension. He has shown that the canonical random point processes converges in probability towards a canonical deterministic measure on $X$, coinciding with the canonical measure of Song-Tian.

Firstly, we give an overview on Analytical Minimal Model Program and explain how we can use of conical Kähler Ricci flow on holomorphic fiber spaces in AMMP to get a log canonical measure and hence a canonical metric in pair $(X, D)$.

We assume in this paper $(X,g)$ is a Kähler manifold. We say the Kähler metric g is Kähler Einstein metric on X, if it satisfies in $Ric(g)=\lambda g\in c_{1}(X)$ where $\lambda$ is constant and it can be normalized as $\lambda\in\{-1,0,1\}$.

Now, given a metric $g$, we can define a matrix valued 2-form $\Omega$ by writing its expansion in local coordinates, as follows

\begin{equation*}
\Omega_{i}^{j} = \sum_{i, p = 1}^{n} g^{j\bar{p}} R_{i\bar{p}k\bar{l} }dz^{k} \wedge d\bar{z}^{l}
\end{equation*}

In fact, this expansion for $\Omega$ gives a well-defined $(1,1)$ form to be called the curvature form of the metric $g$.

Nex, we consider the following expansion

\begin{equation*}
\det \left(Id + \frac{t\sqrt{-1}}{2\pi} \Omega\right) = 1 + t\phi_1(g) + t^2\phi_2(g) + \dots
\end{equation*}

each of the forms $\phi_i(g)$ is a $(i,i)$-form and is called $i$-th Chern form of the metric $g$. The cohomology class represented by each $\phi_i(g)$ is independent on the metric $g$. Since we need the analysis of first Chern class, we restrict our attention to

\begin{equation*}
\phi_1(g) = \frac{\sqrt{-1}}{2\pi} \sum_{i=1}^n \Omega_i^i = \sum_{i,p=1}^n g^{i\bar{p}} R_{i\bar{p}k\bar{l}} dz^k \wedge d\bar{z}^l
\end{equation*}

but the right hand side is just $Ric(g)$. So we have $c_1(X) = [Ric(g)]$


The Ricci curvature of the form $\omega = \sqrt{-1} g_{i\bar{j}}dz_i\wedge d\bar{z}_j$ , is satisfies in the following identity

\begin{equation*}
Ric(\omega) = -\sqrt{-1}\partial\bar{\partial}\log\det g_{i\bar{j}} = -\sqrt{-1}\partial\bar{\partial}\log\omega^n
\end{equation*}

It is clear that $\omega^n$ can be seen as a section of $K_{X} \otimes \bar{K}_{X}$, in other words, $\omega^n$ is a Hermitian metric on the holomorphic line bundle $K_{X}^{-1}$. Note that by Chern-Weil theory, because Ricci form of any Kähler metric on $X$ defines the same de Rham cohomology class, we can define the first Chern class by Ricci form $c_{1}(X) = c_{1}(K_{X}^{-1}) = [Ric(\omega)]$.

The Kähler condition requires that $\omega$ is a closed positive (1, 1)-form. In other words, the following hold

\begin{equation*}
\frac{\partial g_{i\bar{k}}}{\partial z^j} = \frac{\partial g_{j\bar{k}}}{\partial z^i} \quad \text{and} \quad \frac{\partial g_{k\bar{l}}}{\partial z^{\bar{j}}} = \frac{\partial g_{k\bar{j}}}{\partial z^{\bar{l}}}
\end{equation*}

Now, because we are in deal with singularities, so we use of $(1,1)$-current instead of $(1,1)$-forms which is singular version of forms. A current is a differential form with distribution coefficients. Let, give a definition of current here. We recall a singular metric $h_{\text{sing}}$ on a Line bundle L which locally can be written as $h_{\text{sing}} = e^{\phi}h$ where h is a smooth metric, and $\phi$ is an integrable function. Then one can define locally the closed current $T_{L,h_{\text{sing}}}$ by the following formula

\begin{equation*}
T_{L,h_{\text{sing}}} = \omega_{L,h} + \frac{1}{2i\pi}\partial\bar{\partial}\log\phi
\end{equation*}

The current Geometry is more complicated than symplectic geometry. For instance in general one can not perform the wedge product of currents due to this fact that one can not multiply the coefficients which are distributions and in general the product of two distributions is not well defined. However, in some special cases one can define the product of two currents. Here we mention the following important theorem about wedge product of two currents

\vspace{0.3cm}
\noindent \textbf{Theorem 1.1} \textit{Let $\Theta$ be a positive $(p,p)$-current and $T$ be a positive $(1,1)$-current. Assume, for simplicity, that one of these currents has smooth coefficients. Then the wedge product $\Theta \wedge T$ is a positive $(p + 1,p + 1)$-current}
\vspace{0.3cm}

In the theorem above it is important that one of the currents is of type (1, 1). Note that for currents of higher bidegree this theorem is not true. Dinew showed that the wedge


product of two smooth positive (2, 2)-currents in $\mathbb{C}^4$ may fail to be positive. One then simply defines the space of currents to be the dual of space of smooth forms, defined as forms on the regular part $X_{\text{reg}}$ which, near $X_{\text{sing}}$, locally extend as smooth forms on an open set of $\mathbb{C}^N$ in which $X$ is locally embedded.

The topology on space of currents are so important. In fact the space of currents with weak topology is a Montel space, i.e., barrelled, locally convex, all bounded subsets are precompact which here barrelled topological vector space is Hausdorff topological vector space for which every barrelled set in the space is a neighbourhood for the zero vector.

Also because we use of push-forward and Pull back of a current and they can cont be defined in sense of forms, we need to introduce them. If $f : X \to Y$ be a holomorphic map between two compact Kähler manifolds then one can push-forward a current $\omega$ on X by duality setting

\begin{equation*}
\langle f_* \omega , \eta \rangle := \langle \omega , f^* \eta \rangle
\end{equation*}

In general, given a current $T$ on $Y$, it is not possible to define its pull-back by a holomorphic map. But it is possible to define pull-back of positive closed currents of bidegree (1, 1). We can write such currents as $T = \theta + dd^c\varphi$ where $\theta \in T$ is a smooth form, and thus one defines the pull-back of current $T$ as follows

\begin{equation*}
f^* T := f^* \theta + d d^c \varphi \circ f
\end{equation*}

Let $X$ and $Y$ be compact Kähler manifolds and let $f: X \to Y$ be the blow up of $Y$ with smooth connected center $Z$ and $\omega \in H^{1,1}(X, \mathbb{R})$. Demailly showed that

\begin{equation*}
\omega = f^* f_* \omega + \lambda E
\end{equation*}

where $E$ is the exceptional divisor and $\lambda \geq -\nu(\omega, Z)$ where $\nu(\omega, Z) = \inf_{x \in Z} \nu(\omega, x)$ and $\nu(\omega, x)$ is the Lelong number.

As an application, the pushforward $f_*\Omega$ of a smooth nondegenerate volume form $\Omega$ on X with respect to the holomorphic map $\pi: X \to X_{\text{can}}$ is defined as follows: From definition of pushforward of a current by duality, for any continuous function $\psi$ on $X_{\text{can}}$, we have

\begin{equation*}
\int_{X_{\text{can}}} \psi f_* \Omega = \int_{X} (f^* \psi) \Omega = \int_{y \in X_{\text{can}}} \int_{\pi^{-1} (y)} (f^* \psi) \Omega
\end{equation*}


and hence on regular part of $X_{can}$ we have

\begin{equation*}
\pi_* \Omega = \int_{\pi^{-1}(y)} \Omega
\end{equation*}

Note that if $\omega$ is Kähler then,

\begin{equation*}
d\mathrm{Vol}_{\omega_y}(X_y) = df_*(\omega^n) = f_*(d\omega^n) = 0
\end{equation*}

So, $\mathrm{Vol}(X_y) = C$ for some constant $C > 0$ for every $y \in X_{can}$ where $\pi^{-1}(y) = X_y$. See [1][2]. Moreover direct image of volume form $f_*\omega_X^n = \sigma\omega_{X_{can}}^m$ where $\sigma \in L^{1+\epsilon}$ for some positive constant $\epsilon$, see[32],[11].

\textbf{Theorem 1.2} \textit{If $T$ is a positive $(1,1)$-current then it was proved in [2] that locally one can find a plurisubharmonic function $u$ such that}

\begin{equation*}
\sqrt{-1}\partial\bar{\partial}u = T
\end{equation*}

\textit{Note that, if $X$ be compact then there is no global plurisubharmonic function $u$.}

Here we mention an important Riemann extension theorem for plurisubharmonic functions:

\textbf{Riemann extension theorem [2]} : Let $D$ is a proper subvariety of $X$ and $\dot{\psi}$ is a plurisubharmonic function on $X \setminus D$. Assume that for each point $x \in D$, there exists a neighborhood $U$ such that $\psi$ is bounded from above on $U \setminus (U \cap D)$. Then $\dot{\psi}$ extends uniquely to a plurisubharmonic function over $X$. Here $\psi|_{X \setminus D} = \dot{\psi}$

The central problem in Kähler geometry is to finding existence and uniqueness of Kähler Einstein metrics and more generally finding canonical metrics on Kähler manifolds. The canonical way to finding Kähler Einstein metrics is deforming metrics by Hamilton's Ricci flow[76] and in Kähler setting this flow is named as Kähler Ricci flow by S.T.Yau in his celebrating work on the proof of Calabi conjecture [15]. The Kähler Ricci flow is defined by

\begin{equation*}
\frac{\partial}{\partial t} g(t) = -Ric(g)
\end{equation*}

The nice thing on running Kähler Ricci flow on Kähler manifold is that if the initial metric is Kähler, as long as you have smooth solution, then the evolving metric will be Kähler and in fact we can say the Kähler condition is preserved by Kähler Ricci flow. The important fact on Ricci flow is that Ricci flow does surgeries by itself and the


Kähler Ricci flow with surgeries should be a continious path in the Gromov-Hausdorff moduli space even can be stronger after rescaling. Our general philosophy on surgery theory is that PDE surgery by Kähler Ricci flow is exactly the same as Geometric surgery and also Algebraic surgery by flips and flops, albeit algebraic surgery may not exists.

By rescaling metric in the same time, we can define the normalized Kähler Ricci flow as follows

\begin{equation*}
\frac{\partial}{\partial t} g_{i\bar{j}}(t) = -Ric(g_{i\bar{j}}) + \lambda g_{i\bar{j}}
\end{equation*}

So, if this flow converges then we can get Kähler Einstein metric. We have classical results for existence of Kähler-Einstein metrics when the first Chern class is negative, zero or positive. Cao [16] by using parabolic estimates of Aubin [17] and Yau showed that if $c_1(X) < 0$ or canonical line bundle is positive then the normalized Kähler Ricci flow corresponding to $\lambda = -1$ converges exponentially fast to Kähler Einstein metric for any initial Kähler metric. Cao also showed that if canonical line bundle be torsion or $c_1(X) = 0$, then the normalized Kähler Ricci flow corresponding to $\lambda = 0$ converges to the Ricci flat Kähler metric in the Kähler class of the initial Kähler metric which gives an alternative proof of Calabi conjecture. Perelman [18][19] showed that when the canonical line bundle is negative or $c_1(X) > 0$, (such type of manifolds are called Fano manifolds) the normalized Kähler Ricci flow corresponding to $\lambda = 1$ converges exponentially fast to Kähler Einstein metric if the Fano manifold $X$ admits a Kähler Einstein metric and later Tian-Zhu [20] extended this result to Kähler Ricci soliton. In Fano case, there is a non-trivial obstruction to existence of Kähler Einstein metric such as Futaki invariant [22], Tian's $\alpha$-invariant [21]. In 2012, Chen, Donaldson, and Sun [6][7][8] and independently Tian [9] proved that in this case existence is equivalent to an algebro-geometric criterion called Tian's stability or called K-stability. So in these three cases we always required that the first Chern class either be negative, positive or zero. But the problem is that most Kähler manifolds do not have definite or vanishing first Chern class. So the question is does exists any canonical metrics on such manifolds and how does the Kähler Ricci flow behave on such manifolds.J.Song and G.Tian started a program to give answer to this question which is known as Song-Tian program now. For answering to these questions let first give a naive description on Mori's Minimal Model Program [23, 24]. Let $X_0$ be a projective variety with canonical line bundle $K \to X_0$ of Kodaira dimension

\begin{equation*}
\kappa(X_0) = \limsup \frac{\log \dim H^0(X_0, K^{\otimes \ell})}{\log \ell}
\end{equation*}


This can be shown to coincide with the maximal complex dimension of the image of $X_0$ under pluri-canonical maps to complex projective space, so that $\kappa(X_0) \in \{-\infty, 0, 1, \dots, m\}$. Also since in general we work on Singular Kähler variety we need to notion of numerical Kodaira dimension instead of Kodaira dimension.

\begin{equation*}
\kappa_{num}(X) = \sup_{k \geq 1} \left[ \limsup_{m \to \infty} \frac{\log \dim_{\mathbb{C}} H^0(X, mK_X + kL)}{\log m} \right]
\end{equation*}

where $L$ is an ample line bundle on $X$.Note that the definition of $\kappa_{num}(X)$ is independent of the choice of the ample line bundle $L$ on $X$. Siu formulated that the abundance conjecture is equivalent as

\begin{equation*}
\kappa_{kod}(X) = \kappa_{num}(X)
\end{equation*}

Numerical dimension is good thing.

It is worth to mention that if $f: X \to Y$ be an algebraic fibre space and $\kappa(X) \geq 0$, $\kappa(Y) = \dim Y$, (for example Iitaka fibration), then $\kappa(X) = \kappa(Y) + \kappa(F)$, where $F$ is a general fibre of $F$.

A Calabi-Yau $n$-fold or Calabi-Yau manifold of (complex) dimension $n$ is defined as a compact $n$-dimensional Kähler manifold $M$ such that the canonical bundle of $M$ is torsion, i.e. there exists a natural number $m$ such that $K_M^{\otimes m}$ is trivial. In the case when we have the pair $(X,D)$, then log Calabi-Yau variety can be defined when $m(K_X + D)$ is trivial for some $m \in \mathbb{N}$ as $\mathbb{Q}$ Cartier. As for the adjunction formula, it definitely works as long as $K_X + D$ is $Q$-Cartier and it works up to torsion if it is $Q$-Cartier. If it is not $\mathbb{Q}$-Cartier, it is not clear what the adjunction formula should mean, but even then one can have a sort of adjunction formula involving $Ext$'s but this is almost Grothendieck Duality then. The log Minimal Model and abundance conjectures would imply that every variety of log Kodaira dimension $\kappa(X) = 0$ is birational to a log Calabi-Yau variety with terminal singularities.

From the definition of Gang Tian and Shing-Tung Yau [25]. In the generalization of the definition of Calabi-Yau variety to non-compact manifolds, the difference $(\Omega \wedge \bar{\Omega} - \omega^n / n!)$ must vanish asymptotically. Here, $\omega$ is the Kähler form associated with the Kähler metric, $g$.

For singular Calabi-Yau variety we can take $\kappa_{num}(X)=0$ as definition of singular Calabi-Yau variety when $K_X$ is pseudo-effective.


\subsection*{1.1 \quad Song-Tian program}

A compact complex $m$-manifold $X_0$ is said to be of general type if $\kappa(X_0) = m$. If $K_{X_0} \geq 0$ then $X_0$ is a Minimal Model by definition. Now, let the canonical line bundle $K_{X_0}$ is not be semi-positive, then $X_0$ can be replaced by sequence of varieties $X_1, \dots, X_m$ with finitely many birational transformations, i.e., $X_i$ isomorphic to $X_0$ outside a codimension 1 subvariety such that

\begin{equation*}
K_{X_m} \geq 0
\end{equation*}

and we denote by $X_{min} = X_m$ the minimal model of $X_0$. Hopefully using minimal model program $X_0$ of non-negative Kodaira dimension can be deformed to its minimal model $X_{min}$ by finitely many birational transformations and we can therefore classify projective varieties by classifying their minimal models with semi-positive canonical bundle. To deal with our goal, we need to explain Abundance conjecture. Roughly speaking, Abundance conjecture tells us that if a minimal model exists, then the canonical line bundle $X_{min}$ induces a unique holomorphic map

\begin{equation*}
\pi : X_{min} \to X_{can}
\end{equation*}

where $X_{can}$ is the unique canonical model of $X_{min}$. The canonical model completely determined by variety X as follows,

\begin{equation*}
X_{can} = \operatorname{Proj} \left(\bigoplus_{m \geq 0} H^0(X, K_{X_{min}}^m)\right)
\end{equation*}

So combining MMP and Abundance conjecture, we directly get

\begin{equation*}
X \to X_{min} \to X_{can}
\end{equation*}

which $X$ here might not be birationally equivalent to $X_{can}$ and dimension of $X_{can}$ is smaller than of dimension $X$. Note also that, Minimal model is not necessary unique but $X_{can}$ is unique. Note that if the canonical Ring $R(X, K_X)$ is finitely generated then the pluricanonical system induces an algebraic fibre space $\pi': X \to X_{can}$. For example if $K_X$ be semi-ample then $R(X, K_X)$ is finitely generated and we have an algebraic fibre space $\pi': X \to X_{can}$.

Now we explain how it is related to Kähler Einstein metric and normalized Kähler Ricci flow. By Song-Tian program it turns out that the normalized Kähler Ricci flow doing exactly same thing to replace X by its minimal model by using finitely many geometric

\newpage

surgeries and then deform minimal model to canonical model such that the limiting of canonical model is coupled with generalized Kähler Einstein metric twisted with Weil-Petersson metric gives canonical metric for $X$. Song-Tian program on MMP is that if $X_0$ be a projective variety with a smooth Kähler metric $g_0$, we apply the Kähler Ricci flow with initial data $(X_0, g_0)$, then there exists $0 < T_1 < T_2 \dots < T_{m+1} \leq \infty$ for some $m \in \mathbb{N}$, such that

\begin{equation*}
(X_0, g_0) \xrightarrow{t \to T_1} (X_1, g_1) \xrightarrow{t \to T_2} \dots \xrightarrow{t \to T_m} (X_m, g_m)
\end{equation*}

after finitely many surgeries in Gromov-Hausdorff topology and either $dim_{\mathbb{C}}X_m < n$, or $X_m = X_{min}$ if not collapsing. In the case $dim_{\mathbb{C}}X_m < n$, then $X_m$ admits Fano fibration which is a morphism of varieties whose general fiber is a Fano variety of positive dimension in other words has ample anti-canonical bundle. In the case $X_m = X_{min}$ is a minimal model and after appropriate normalization, we have long-time existence for the solutions of normalized Kähler Ricci flow to canonical metric

\begin{equation*}
(X_{min}, g_m) \xrightarrow{t \to \infty} (X_{can}, g_{can})
\end{equation*}

Now we are ready to give Song-Tian algorithm via MMP along Kähler Ricci flow.

$\mathbb{I}$-Start with a pair $(X,g_0)$ where $X$ is a projective variety with with log terminal singularities. Apply the Kähler Ricci flow

\begin{equation*}
\frac{\partial}{\partial t} g(t) = -Ric(g)
\end{equation*}

starting with $g_0$ and let $T > 0$ be the first singular time

$\mathbb{III}$- If the first singular time is finite, $T < \infty$, then

\begin{equation*}
(X, g(t)) \xrightarrow{G-H} (Y, d_Y)
\end{equation*}

(G-H here means Gromov-Hausdorff topology). Here $(Y, d_Y)$ is a compact metric space homeomorphic to the projective variety Y. Now there is three possibilities which we mention here.

$\mathbb{A}$ -In the case $dimY = dimX$. Then $(X, g(t))$ is deformed continuously to a compact metric space $(X^+, g_{X^+}(t))$ in Gromov-Hausdorff topology as $t$ passes through $T$, where is a normal projective variety satisfying the following diagram

\begin{figure}[H]
  \centering
  \includegraphics[width=0.35\textwidth,height=0.82\textheight,keepaspectratio]{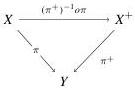}
\end{figure}

and $X^{+}$ is $\pi^{+}$-ample. $\pi^{+}: X^{+} \to Y$ is a general flip of $X$ and we have

\begin{equation*}
X^{+} = Proj \left( \bigoplus_{m} \pi_{*} (\mathcal{O}_{X}(K_{X}^{\otimes m})) \right)
\end{equation*}

So essentially the Kähler Ricci flow should deform $X$ to a new variety $Y$ in Gromov-Hausdorff topology and if $Y$ has mild singularities then the Kähler Ricci flow smoothout the metric and we continue the flow on $Y$. In fact Kähler Ricci flow is a non-linear parabolic heat equation and so it must smoothout singularities and in this case $Y$ and $X^{+}$ are the same thing. If $Y$ has very bad singularities then we can still apply the Ricci flow on $Y$ and the flow immediately resolve the singularities to a new variety $X^{+}$ and the resolution of singularities is in Gromov-Hausdorf topology and we have continuous resolution of singularities and $Y$ continuously jump to $X^{+}$ and we can continue the Kähler Ricci flow. Gang Tian jointly with Gabriele La Nave started a program to studying finite time singularity of the Kähler-Ricci flow, for those varieties that flips can be achieved through variations of symplectic quotients via V-soliton metrics.

$\mathbb{B}$ -In the case $0 < dimY < dimX$.

$(X, D)$ admits log Fano fibration over $Y$. Then $(X \setminus D, g(t))$ collapses and converges to a metric space $(Y, g_{Y})$ in Gromov-Hausdorff topology as $t \to T$ and we replace $(X, g_{0})$ by $(Y, g_{Y})$ and continue the flow.

$\mathbb{C}$ -In the case $dimY = 0$.

Then $(X, D)$ is the log Fano and the flow becomes extinct at $t \to T$ and if we rescale it then the flow converges to general version of conical Kähler Ricci soliton.

So the previous situations were for finite time singularity $T < \infty$. If we have long time existence $T = \infty$, then we also have the following three situations.

$\mathbb{A}$ -In the case $\kappa(X, D) = dimX$.


Then $\frac{1}{t} g(t)$ converges to unique conical Kähler Einstein metric $g_{KE}$ on $X_{can}^{D}$ in Gromov-Hausdorff topology as $t \to \infty$.

In the case when $X$ has log-terminal singularities, then Song and Tian gave an affirmative answer to this part and later Odaka [60] showed that If $(X, L)$ admits Kähler Ricci flat metric then the singularities of $X$ must be mild.

$\mathbb{B}$ -In the case $0 < \kappa (X,D) < dimX$.

Then $X$ admits a Calabi-Yau fibration over its canonical model $X_{can}$. Moreover, $\frac{1}{t} g(t)$ collapses and converges to the unique generalized Kähler Einstein metric $g_{can}$ satisfying to the generalized kähler-Einstein equation on $X_{can}$ in Gromov-Hausdorff topology as $t \to \infty$ as soon as we have fiberwise KE-stability or Kawamata's condition of Theorems 2.28, or 2.30 (for $C^0$-case).

\begin{equation*}
Ric(\omega_{can}) = -\omega_{can} + \omega_{WP}^{D} + \sum_{P} (b(1 - t_{P}^{D})) [ \pi^{*} (P) ] + [ B^{D} ]
\end{equation*}

where $B^{D}$ is $\mathbb{Q}$-divisor on $X$ such that $\pi_{*}\mathcal{O}_{X}([iB_{+}^{D}]) = \mathcal{O}_{B} (\forall i > 0)$. Here $s_{P}^{D} := b(1 - t_{P}^{D})$ where $t_{P}^{D}$ is the log-canonical threshold of $\pi^{*}P$ with respect to $(X, D - B^{D}/b)$ over the generic point $\eta_{P}$ of $P$. i.e.,

\begin{equation*}
t_{P}^{D} := \max \left\{t \in \mathbb{R} \mid (X, D - B^{D} / b + t \pi^{*} (P)) \text{ is sub log canonical over } \eta_{P} \right\}
\end{equation*}

and $\omega_{can}$ has zero Lelong number. We conjecture that $\omega_{can}$ is a good metric in the sense of Mumford

Recently Song and Tian [26] solved the weak convergence in the case when $X$ has singularity and canonical bundle is semi ample. They showed that the flow converges and collapses with uniformly bounded scalar curvature. Note that if we accept Abundance conjecture then it is always true.

$\mathbb{C}$ -In the case $\kappa(X,D)=0$.

Then the conical Kähler Ricci flow converges to the unique Ricci flat Kähler metric in the initial Kähler class in Gromov-Hausdorff topology as $t \to \infty$. Song-Tian solved the weak convergence for this case.

So if we have long time existence $T = \infty$ then the main goal of this paper is to extend the case $\mathbb{B}$ when $0 < \kappa(X, D) < dimX$ for pair $(X, (1 - \beta)D)$ where $D$ is divisor on $X$, $0 < \beta < 1$ and we can write


\begin{equation*}
Ric(\omega_{can}) = -\omega_{can} + \omega_{WP}^D + \sum_P (b(1 - t_P^D))[\pi^*(P)] + [B^D]
\end{equation*}

where $B^D$ is $\mathbb{Q}$-divisor on $X$ such that $\pi_* \mathcal{O}_X([iB_+^D]) = \mathcal{O}_B$ ($\forall i > 0$). Here $s_P^D := b(1 - t_P^D)$ where $t_P^D$ is the log-canonical threshold of $\pi^* P$ with respect to $(X, D - B^D/b)$ over the generic point $\eta_P$ of $P$. i.e.,

\begin{equation*}
t_P^D := \max\{t \in \mathbb{R} \mid (X, D - B^D/b + t\pi^*(P)) \text{ is sub log canonical over } \eta_P\}
\end{equation*}

and $\omega_{can}$ has zero Lelong number and $\pi$ induces a function $f : X_0 \to \mathcal{M}_{CY}^D$ where $X_0$ is a dense-open subset over

\begin{equation*}
X_{can}^D = \text{Proj}R(X; K_X + D) = \text{Proj} \bigoplus_{m \ge 0} H^0(X, mK_X + \lfloor mD \rfloor)
\end{equation*}

which the fibers of $\pi$ are smooth. Note that the moduli space of log Calabi Yau pair $\mathcal{M}_{LCY}^D$ is essentially smooth (if exists!) and admits a canonical metric, that is, the logarithmic Weil-Peterson metric.

It is worth to mention that in Minimal Model Program

\begin{equation*}
X = X_1 \dashrightarrow \cdots \dashrightarrow X_k \to Y
\end{equation*}

\noindent 1) If $\kappa(X) = -\infty$ then $X_k \to Y$ is Mori fiber space with fiber being Fano variety of Picard number 1.

\noindent 2) If $0 \le \kappa(X) < n$ then $X_k \to Y$ is a Calabi-Yau fiber space, i.e., fibers are Calabi-Yau varieties

\noindent 3) If $\kappa(X) = n$ then $Y = X^{can}$ is a canonically polarized variety.

\subsection*{1.2 \quad Lelong number:}

Let $W \subset \mathbb{C}^n$ be a domain, and $\Theta$ a positive current of degree $(q, q)$ on $W$. For a point $p \in W$ one defines

\begin{equation*}
\mathfrak{v}(\Theta, p, r) = \frac{1}{r^{2(n-q)}} \int_{|z-p|<r} \Theta(z) \wedge (dd^c|z|^2)^{n-q}
\end{equation*}

The Lelong number of $\Theta$ at $p$ is defined as

\begin{equation*}
\mathfrak{v}(\Theta, p) = \lim_{r \to 0} \mathfrak{v}(\Theta, p, r)
\end{equation*}

\newpage

Let $\Theta$ be the curvature of singular hermitian metric $h = e^{-u}$, one has

\begin{equation*}
\mathfrak{v}(\Theta , p) = \sup \left\{\lambda \geq 0: u \leq \lambda \log \left(\left| z - p \right| ^{2}\right) + O (1) \right\}
\end{equation*}

\textbf{Siu's Decomposition[61]}: Let $X$ be a complex manifold and $T$ a closed positive current of bidimension $(p, p)$. Then, there is a unique decomposition of $T$ as a (possibly finite) weakly convergent series

\begin{equation*}
T = \sum_{j \geq 1} \lambda_{j} [ A _{j} ] + R, \quad \lambda_{j} \geq 0,
\end{equation*}

where $[A_{j}]$ is the current of integration over an irreducible p-dimensional analytic set $A_{j} \subset X$ and where R is a closed positive current with the property that $\dim E_{c}(R) < p$ for every $c > 0$. Here $E_{c}(R)$ is the set of point $x \in X$ such that the Lelong number $\nu(R, x)$ of R at x is greater than or equal to c. Later, when we study AZD and Zariski Decomposition, we can see the importance of this decomposition.

Lelong number give a lot of information to us. See [13][40].

\noindent \textbf{A1)} If the line bundle $L \to X$ be a nef and big then there exists a hermitian metric on $L$ with vanishing Lelong number.

\noindent \textbf{A2)} Let $X$ be a projective manifold and $(L, h)$ a positive, singular hermitian line bundle, whose Lelong numbers vanish everywhere. Then $L$ is nef.

\noindent \textbf{A3)} Let $X$ be a smooth projective variety and let $D$ be a nef and big $\mathbb{R}$-divisor on $X$. Then the first Chern class $c_{1}(D)$ can be represented by a closed positive (1,1)-current $T$ with $\mathfrak{v}(T) = 0$

\noindent \textbf{A4)} Let $X$ be a smooth projective variety and let $D$ be a $\mathbb{R}$-divisor on $X$ such that $c_{1}(D)$ can be represented by a closed positive (1, 1) current $T$ with $\mathfrak{v}(T) = 0$. Then $D$ is nef.

\noindent \textbf{A5)} Suppose that there exists a modification $f: Y \to X$ such that there exists a Zariski decomposition $f^{*}L = P + N$ of $f^{*}L$ on $Y$. Then there exists a closed positive (1, 1) current $S$ such that $c_{1}(P) = [S]$ and $\mathfrak{v}(S) = 0$

By using previous A1--A5, we have the following remark.

\noindent \textbf{Theorem}: Canonical metric (Kähler Einstein metric, constant scalar curvature, twisted Kähler Einstein metric along $X \to X_{can}$) has zero Lelong number.


\section*{2 \quad Generalized K\"ahler-Einstein metric along Long Canonical Model}

In this section by assuming finite generation of log canonical ring, we try to construct a unique logarithmic canonical measure and hence we try to find a unique twisted Kähler Einstein metric on pair $(X, D)$. When the log canonical ring is finitely generated, the canonical model can be defined by Proj of the log canonical ring as follows

\begin{equation*}
X_{can}^{D} = \text{Proj} R(X, K_{X} + D) = \text{Proj} \bigoplus_{m \geq 0} H^{0}(X, mK_{X} + \lfloor mD \rfloor)
\end{equation*}

here $\lfloor x \rfloor = \max \{m \in \mathbb{Z} \mid m \leq x\}$. If the log canonical ring is not finitely generated, then $\text{Proj} R(X, K_{X} + D)$ is not a variety, and so it cannot be birational to $(X, D)$; in particular, $(X, D)$ admits no canonical model. If the log canonical ring be finitely generated, then the log canonical model $X_{can}^{D}$ is unique. Moreover we show that our log canonical measure is birational invariant.

Recently, Birkar, Cascini, and Hacon,[27] and independently Siu [28] showed the following theorems

\noindent \textbf{Theorem 2.1} \textit{Let $X$ be a smooth projective variety. Then the canonical ring $R(X; K_{X})$ is a finitely generated $\mathbb{C}$-algebra and hence canonical model exists and is unique. If $D$ be a smooth divisor on $X$ then log canonical ring is finitely generated also.}

A pair $(X, D)$ is said Kawamata log terminal if for any proper birational morphism $\phi : Y \to X$ we may write

\begin{equation*}
K_{Y} + D_{Y} = \phi^{*}(K_{X} + D)
\end{equation*}

where $D_{Y}$ is a $\mathbb{Q}$-divisor whose coefficients are strictly less than one.

\noindent \textbf{Theorem 2.2} \textit{Let $(X, D)$ be a projective Kawamata log terminal pair. Then the log canonical ring $R(X; K_{X} + D)$ is finitely generated and hence log canonical model exists and is unique}

Because we are in deal with singular Hermitian metrics, for our next results, we need to consider a kind of Zariski decomposition which is called by Analytic Zariski Decomposition (AZD). In fact for our purpose, we are interested to study $H^{0}(X, \mathcal{O}_{X}(mL) \otimes \mathcal{I}(h^{m}))$. Now, If a line bundle $L$ admits an Analytic Zariski Decomposition $h$, then to study $H^{0}(X, \mathcal{O}_{X}(mL))$ it is sufficient to study $H^{0}(X, \mathcal{O}_{X}(mL) \otimes \mathcal{I}(h^{m}))$. In fact this definition of Analytic Zariski Decomposition of Tsuji [13] is inspired by Zariski. Zariski for

considering the finite generation of canonical rings introduced the notion of Zariski decomposition, which says that if $X$ be a smooth projective variety and $D$ be a pseudoeffective divisor on $X$, then $D = P + N$ is said to be Zariski decomposition of $D$, if $P$ be nef and $N$ be effective and we have isomorphism
\begin{equation*}
H^{0}(X, \mathcal{O}_{X}([mP])) \to H^{0}(X, \mathcal{O}_{X}([mD])).
\end{equation*}
So if the divisor $D$ be Zariski decomposition then $R(X, \mathcal{O}_X([P])) \simeq R(X, \mathcal{O}_X([D]))$ holds true.

\textbf{Definition 2.3} Let $L$ be a line bundle on a complex manifold $X$. A singular Hermitian metric $h$ is given by
\begin{equation*}
h = e^{-\varphi} \cdot h_{0}
\end{equation*}
where $h_0$ is a $C^\infty$-Hermitian metric on $L$ and $\varphi \in L_{loc}^1(X)$ is an arbitrary function on $X$. The curvature current $\Theta_h$ of the singular Hermitian line bundle $(L, h)$ is defined by
\begin{equation*}
\Theta_{h} := \Theta_{h_{0}} + \sqrt{-1} \partial \bar{\partial} \varphi
\end{equation*}
where $\partial\bar{\partial}$ is taken in the sense of a current as introduced before.

To give an another equivalent definition for singular hermitian metric. Let $L \to X$ be a holomorphic line bundle over a complex manifold $X$ and fix an open cover $X = \cup U_{\alpha}$ for which there exist local holomorphic frames $e_{\alpha}: U_{\alpha} \to L$. The transition functions $g_{\alpha \beta} = e_{\beta} / e_{\alpha} \in \mathcal{O}_X^*(U_\alpha \cap U_\beta)$ determine the Cech 1-cocycle $\{(U_\alpha, g_{\alpha \beta})\}$. If $h$ is a singular Hermitian metric on $L$ then $h(e_{\alpha}, e_{\alpha}) = e^{-2\varphi_{\alpha}}$, where the functions $\varphi_{\alpha} \in L_{loc}^{1}(U_{\alpha})$ are called the local weights of the metric $h$. We have $\varphi_{\alpha} = \varphi_{\beta} + \log |g_{\alpha \beta}|$ on $U_{\alpha} \cap U_{\beta}$ and the curvature of $h$,
\begin{equation*}
c_{1}(L, h)|_{U_{\alpha}} = d d^{c} \varphi_{\alpha}
\end{equation*}
is a well defined closed $(1,1)$ current on $X$.

One of the important example of singular hermitian metric is singular hermitian metric with algebraic singularities. Let $m$ be a positive integer and $\{S_{i}\}$ a finite number of global holomorphic sections of $mL$. Let $\varphi$ be a $C^{\infty}$-function on $X$. Then
\begin{equation*}
h := e^{-\varphi} \cdot \frac{1}{(\sum_{i} |S_{i}|^{2})^{1/m}}
\end{equation*}

defines a singular hermitian metric on $L$. We call such a metric $h$ a singular hermitian metric on $L$ with algebraic singularities. To giving another examples: Let $S \in H^0(X; L)$ then $h := h_0 e^{-\varphi} = \frac{h_0}{|S|_{h_0}^2}$ by taking $\varphi := \log |S|_{h_0}^2$ give a singular hermitian metric. Let $D$ be a $\mathbb{Q}$-effective divisor on smooth projective variety $X$ and $S \in H^0(X, \mathcal{O}_X(mD))$ be a global section then,
\begin{equation*}
h := \frac{h_0}{\left(h_0^m(S, S)\right)^{1/m}}
\end{equation*}
defines an singular hermitian metric on $mD$ where $h_0^m$ means smooth hermitian metric on $mD$. See [2].

For defining Analytic Zariski Decomposition we need to introduce the notion of multiplier ideal sheaf.

\textbf{Definition 2.4} The $L^2$-sheaf $\mathcal{L}^2(L, h)$ of the singular Hermitian line bundle $(L, h)$ is defined by
\begin{equation*}
\mathcal{L}^2(L, h) := \{\sigma \in \Gamma(U, \mathcal{O}_X(L)) \mid h(\sigma, \sigma) \in L_{loc}^1(U)\}
\end{equation*}
where $U$ runs over open subsets of $X$. There exists an ideal sheaf $\mathcal{I}(h)$ such that
\begin{equation*}
\mathcal{L}^2(L, h) = \mathcal{O}_X(L) \otimes \mathcal{I}(h)
\end{equation*}
More explicitly, by taking $h = e^{-\varphi} \cdot h_0$ the multiplier ideal sheaf can be defined by
\begin{equation*}
\mathcal{I}(h) = \mathcal{L}^2(\mathcal{O}_X, e^{-\varphi}).
\end{equation*}
If $h$ is a smooth metric with semi-positive curvature, then $\mathcal{I}(h) = \mathcal{O}_X$, but the converse is not true. In general, let $(L, h_L)$ be a singular Hermitian $\mathbb{Q}$-line bundle on a smooth projective variety $X$, then $(L, h_L)$ is said to be Kawamata log terminal, if its curvature current is semipositive and $\mathcal{I}(h) = \mathcal{O}_X$. For instance, for Weil-Petersson metric we have $\mathcal{I}(h_{WP}) = \mathcal{O}_X$, so $(L_{WP}, h_{WP})$ is Kawamata log terminal Hodge $\mathbb{Q}$-line bundle.

\textbf{Remark:} It is worth to mention that the corresponding Monge-Ampere equation of the Song-Tian metric
\begin{equation*}
Ric(\omega_{can}) = -\omega_{can} + \omega_{WP}^D + \sum_P (b(1 - t_P^D))[\pi^*(P)] + [B^D]
\end{equation*}
where $B^D$ is $\mathbb{Q}$-divisor on $X$ such that $\pi_* \mathcal{O}_X([iB_+^D]) = \mathcal{O}_{X_{can}^D}$ ($\forall i > 0$). Here $s_P^D := b(1 - t_P^D)$ where $t_P^D$ is the log-canonical threshold of $\pi^* P$ with respect to $(X, D - B^D/b)$ over the generic point $\eta_P$ of $P$. i.e.,

\begin{equation*}
t_{P}^{D} := \max \left\{t \in \mathbb{R} \mid \left(X, D - B^{D} / b + t \pi^{*}(P)\right) \text{ is sub log canonical over } \eta_{P} \right\}
\end{equation*}
$\omega_{can}$ can not have algebraic singularities.

Note that $\mathcal{I}(h)$ is not coherent sheaf of $\mathcal{O}_X$-ideals but if $\Theta_h$ be positive then Nadel [29, 30] showed that $\mathcal{I}(h)$ is coherent sheaf of $\mathcal{O}_X$-ideals

Now we are ready to define the notion of Analytic Zariski Decomposition.

\textbf{Definition 2.5} Let $X$ be a compact complex manifold and let $L$ be a line bundle on $X$. A singular Hermitian metric $h$ on $L$ is said to be an analytic Zariski decomposition, if the following hold.
\begin{enumerate}
\item the curvature $\Theta_h$ is a closed positive current.
\item for every $m \geq 0$, the natural inclusion
\begin{equation*}
H^{0}(X, \mathcal{O}_{X}(mL) \otimes \mathcal{I}(h^{m})) \to H^{0}(X, \mathcal{O}_{X}(mL))
\end{equation*}
is an isomorphism, where $\mathcal{I}(h^{m})$ denotes the multiplier ideal sheaf of $h^{m}$.
\end{enumerate}

Now, we mention one of the important theorems of Tian-Yau and Tsuji about AZD and existence of Kähler Einstein metric on pair $(X,D)$ where $D$ is simple normal crossing divisor and $K_{X}+D$ is ample.

Let $\sigma_{D}$ be a nontrivial global holomorphic section of divisor $D$. We denote the singular hermitian metric on $D$ by
\begin{equation*}
h_{D} = \frac{1}{\|\sigma_{D}\|^{2}} := \frac{h_{0}}{h_{0}(\sigma_{D}, \sigma_{D})}
\end{equation*}
where $h_0$ is an arbitrary $C^\infty$-hermitian metric on $D$. We can define the Bergman kernel of $H^0(X, \mathcal{O}_X(D) \otimes \mathcal{I}(h_D))$ with respect to the $L^2$-inner product:
\begin{equation*}
(\sigma, \sigma^{\prime}) := (\sqrt{-1})^{n^{2}} \int_{X} h_{D} \cdot \sigma \wedge \bar{\sigma}^{\prime}
\end{equation*}
as
\begin{equation*}
K(K_{X} + D, h_{D}) = \sum_{i = 0}^{N} |\sigma_{i}|^{2}
\end{equation*}
where $\{\sigma_0,\sigma_1,\dots ,\sigma_N\}$ is a complete orthonormal basis of $H^0 (X,\mathcal{O}_X(D)\otimes \mathcal{I}(h_D))$.

\textbf{Theorem 2.6} Let $X$ be a smooth projective variety and let $D$ be a divisor with simple normal crossing on $X$ such that $K_{X} + D$ is ample then there exists a Kähler Einstein current $\omega_{E}$ with $\operatorname{Ric}(\omega_{E}) = -\omega_{E} + [D]$. The metric
\begin{equation*}
h := (\omega_{E}^{n})^{-1}
\end{equation*}
is a singular hermitian metric on $K_{X} + D$ with strictly positive curvature on $X$. Let $D = \sum_{i} D_{i}$ be the irreducible decomposition of $D$ and let $S_{i}$ be a nontrivial global section of $\mathcal{O}_X(D_i)$ with divisor $D_{i}$. Then $h$ is analytic zariski decomposition on $K_{X} + D$ and $h$ has logarithmic singularities along $D$ which means there exists a $C^\infty$-hermitian metric $h_0$ on $K_{X} + D$ such that
\begin{equation*}
h = h_{0} \prod_{i} \left| \log \|S_{i}\|^{2} \right|^{2}
\end{equation*}
where $\|S_{i}\|$ denotes the hermitian norm of $S_{i}$ with respect to a $C^{\infty}$-hermitian metric on $\mathcal{O}_{X}(D_{i})$ respectively. Hence this show that the singular hermitian metric $h$ blows up along simple normal crossing divisor $D$.

As an application, let $\omega$ be a complete Kähler-Einstein form on $M$ with $\dim X = n$ such that $Ric(\omega) = -\omega + [D]$, then
\begin{equation*}
h = \omega^{n} \cdot h_{D}
\end{equation*}
is an Analytic Zariski Decomposition on $K_{X} + D$.

\textbf{Definition 2.7} We say a line bundle on a projective complex manifold $X$ is pseudo-effective if and only if it admits a singular Hermitian metric $h_{\mathrm{sing}}$ whose associated Chern current $T_{L,h_{\mathrm{sing}}}$ is positive.

J.-P. Demailly in [2] showed that a line bundle $L$ on a projective manifold $X$ admits an Analytic Zariski Decomposition, if and only if $L$ is pseudo-effective.

For the Kawamata log terminal pair $(X,D)$ such that $X$ is smooth projective variety and $K_{X}+D$ is pseudoeffective then there exists a singular hermitian metric $h_{can}$ on $K_{X}+D$ such that $h_{can}$ is Analytic Zariski decomposition on $K_{X}+D$.

Tsuji in [13] showed that under some certain condition, AZD exists. He proved that if $L$ be a big line bundle, i.e.,
\begin{equation*}
\operatorname{Vol}(L) = \limsup_{m \to \infty} \frac{\dim H^{0}(X, L^{\otimes m})}{m^{\dim X}} > 0
\end{equation*}
then $L$ has an analytic Zariski decomposition.

To state our results we need also to introduce log-Iitaka fibration.

\textbf{Definition 2.8} We say $(X, D)$ is a Kawamata log-terminal pair or shortened as klt pair, if $X$ is a normal projective variety over $\mathbb{C}$ of dimension $n$, and $D$ is an arbitrary $\mathbb{Q}$-divisor such that $K_X + D$ is $\mathbb{Q}$-Cartier, and for some (or equivalently any) log-resolution $\pi : X' \to X$, we have:
\begin{equation*}
K_{X^{\prime}} = \pi^{*}(K_{X} + D) + \sum a_{i} E_{i}
\end{equation*}
where $E_{i}$ are either exceptional divisors or components of the strict transform of $D$, and the coefficients $a_{i}$ satisfy the inequality $a_{i} > -1$.

Log Iitaka fibration is the most naive geometric realization of the positivity of the Log canonical ring.

\textbf{Definition 2.9} Let $X$ be a smooth projective variety with $\operatorname{Kod}(X) \geq 0$. Then for a sufficiently large $m > 0$, the complete linear system $|m!K_X|$ gives a rational fibration with connected fibers:
\begin{equation*}
f: X \dashrightarrow Y.
\end{equation*}
We call $f: X \dashrightarrow Y$ the Iitaka fibration of $X$. Iitaka fibration is unique in the sense of birational equivalence. We may assume that $f$ is a morphism and $Y$ is smooth. For Iitaka fibration $f$ we have
\begin{enumerate}
\item For a general fiber $F$, $\operatorname{Kod}(F) = 0$ holds.
\item $\dim Y = \operatorname{Kod}(X)$.
\end{enumerate}

Similarly, the log-Iitaka fibration of the pair $(X, D)$ gives the asymptotic analysis of the following rational mappings.
\begin{equation*}
\phi_{k}: X \dashrightarrow Y_{k} \subset \mathbb{P}\left(H^{0}(X, (K_{X} + D)^{\otimes k})^{*}\right)
\end{equation*}

\textbf{Theorem 2.10} Let $(X, D)$ be a normal projective variety such that log-Kodaira dimension
\begin{equation*}
\kappa(X, D) = \limsup_{m \in \mathbb{N}} \frac{\log h^{0}(X, K_{X} + D)}{\log m} > 0
\end{equation*}
Then for all sufficiently large in semigroup $k \in N(K_X + D) := \{m \geq 0 : h^0(X, m(K_X + D)) \neq 0\}$ there exists a commutative diagram

\begin{equation*}
\begin{array}{ccc}
X_{\infty, D} & \xrightarrow{\quad u_{\infty} \quad} & X \\
\varphi_{\infty, D} \bigg\downarrow & & \bigg\downarrow \varphi_{k, D} \\
\text{Iitaka}(X, K_X + D) & \xrightarrow{\quad v_{k, D} \quad} & \text{Im}(\varphi_{k, D})
\end{array}
\end{equation*}

where the horizontal maps are birational. One has $\dim(\text{Iitaka}(X, K_X + D)) = \kappa(X, D)$. Moreover if we set $(K_X + D)_{\infty} = u_{\infty}^*(K_X + D)$ and $F$ is the very general fiber of $\varphi_{\infty,D}$, we have $\kappa(F, (K_X + D)_{\infty}|_F) = 0$.

Note that if $K_{X} + D$ be a semi-ample line bundle over a projective manifold $X$. Then there is an algebraic fibre space $\Phi_{\infty}: X \to Y$ such that for any sufficiently large integer $m$ with $m(K_{X} + D)$ being globally generated, $Y_{m} = Y$ and $\Phi_{\infty} = \Phi_{m}$, where $Y$ is a normal projective variety and $\Phi_{m}: X \to P(H^{0}(X, (K_{X} + D)^{\otimes m}))$ and here $Y_{m} := \Phi_{m}(X)$. Here $\Phi_{\infty}$ is called Iitaka fibration. In particular if $K_{X} + D$ be semi-ample then the algebraic fibre space $f: X \to X_{can}$ is called Iitaka fibration. It is worth to mention that if $f: (X, D) \to Y$ be Iitaka fibration then the logarithmic Weil-Petersson metric $\omega_{WP}$ indicate the isomorphism

\begin{equation*}
R(X, K_X + D)^{(m)} = R(Y, L_{X/Y, D})^{(m)}
\end{equation*}

for some positive integer $m$, where for graded ring $R = \bigoplus_{i=0}^{\infty} R_i$, we set

\begin{equation*}
R^{(m)} = \bigoplus_{i=0}^{\infty} R_{mi}
\end{equation*}

and here

\begin{equation*}
L_{X/Y, D} := \frac{1}{m!} \left(f_* \mathcal{O}_X(m!(K_{X/Y} + D))\right)^{**}
\end{equation*}

here $m$ chosen such that $m!(K_X + D)$ become Cartier. We set $K_{X/Y} = K_X \otimes \pi^* K_Y^{-1}$ and call it the relative canonical bundle of $\pi : X \to Y$. Note that the direc image $f_*\mathcal{O}_X(m!(K_{X/Y} + D))$ can be realized as the pull back of an ample vector bundle on the moduli space of Calabi-Yau fibers via moduli map.

Eyssidieux, Guedj, and Zeriahi in [32] showed that if $(X,D)$ be a projective variety of general type (which in this case log canonical ring is finitely generated), then $X_{can}$ will have only canonical singularities and there exists a Kähler Einstein metric on $X_{can}$ with a continuous potential. Furthermore it is continuous on $X$ and smooth on a Zariski open dense set of $X$. Note that only condition of finite generation of canonical

ring doesn't give Kähler Einstein metric. For example, take $\mathbb{CP}^1\times \Sigma_g$ where $\Sigma_{g}$ is a curve with positive genus, then it has no Kähler Einstein metric, but its canonical ring is finitely generated. Our main goal of this section is to prove the following theorem

\textbf{Theorem 2.11} \textit{Let $X$ be an $n$-dimensional projective variety with a simple normal crossing divisor $D$ and with log Kodaira dimension $0 < \kappa < n$. If the log canonical ring is finitely generated, then there exists a unique canonical measure $\Omega_{can}^{D}$ on the pair $(X, D)$ satisfying:}

\textit{(1)- We have $0 < \frac{\Psi_X^D}{\Omega_{can}^D} < \infty$ where}
\begin{equation*}
\Psi_X^D = \sum_{m=0}^M \sum_{j=0}^{d_m} |S_{m,j}|^{\frac{2}{m}}
\end{equation*}
\textit{where $\{S_{m,j}\}_{j=0}^{d_m}$ spans $H^0(X, \mathcal{O}_X(mD) \otimes \mathcal{I}(mh_D))$}

\textit{(2)- $\Omega_{can}^{D}$ is continuous on $X$ and smooth on Zariski open set of $X$}

\textit{(3)- Let $\Phi : (X, D) \dashrightarrow X_{can}^{D}$ be the pluricanonical map. Then there exists a unique closed positive $(1,1)$-current $\omega_{can}^{D}$ with bounded local potential on $X_{can}^{D}$ such that $\Phi^{*}(\omega_{can}^{D}) = \sqrt{-1}\partial\bar{\partial}\log\Omega_{can}^{D}$ outside of the base locus of the pluricanonical system. Furthermore, on $X_{can}^{D}$,}
\begin{equation*}
(\omega_{can}^D)^{\kappa} = \Phi_* \Omega_{can}
\end{equation*}
\textit{and let we have fiberwise KE-stability or Kawamata's condition of Theorems 2.27, or 2.28, then}
\begin{equation*}
\text{Ric}(\omega_{can}) = -\omega_{can} + \omega_{WP}^D + \sum_P (b(1 - t_P^D))[\pi^*(P)] + [B^D]
\end{equation*}
\textit{where $B^D$ is $\mathbb{Q}$-divisor on $X$ such that $\pi_*\mathcal{O}_X([iB_+^D]) = \mathcal{O}_{X_{can}^D}$ ($\forall i > 0$). Here $s_P^D := b(1 - t_P^D)$ where $t_P^D$ is the log-canonical threshold of $\pi^*P$ with respect to $(X, D - B^D/b)$ over the generic point $\eta_P$ of $P$. i.e.,}
\begin{equation*}
t_P^D := \max \left\{ t \in \mathbb{R} \mid \left(X, D - B^D/b + t\pi^*(P)\right) \text{ is sub log canonical over } \eta_P \right\}
\end{equation*}

\textit{(4)- $\Omega_{can}^{D}$ is invariant under birational transformations.}

For proving parts 2),3), and 4) of Theorem 2.11., we need to prove some lemmas. Part 1) of Theorem 2.11 by the same method of Song-Tian is trivial and we get the desired result. So we prove parts 2),3), and 4)

\textbf{Definition 2.12} \textit{Suppose that $X$ is an $n$-dimensional projective variety with simple normal crossing divisor $D$ of log Kodaira dimension $0 < \kappa < n$. Let the log canonical ring $R(X, K_X + D)$ is finitely generated and $\Phi : (X, D) \dashrightarrow X^D_{can}$ be the pluricanonical map. There exists a non-singular model $(X^{\dagger}, D^{\dagger} = \iota^*D)$ (where $\iota : X^{\dagger} \to X$ is the normalization) of $(X, D)$ and the following diagram holds:}

\begin{figure}[H]
  \centering
  \includegraphics[width=0.4\textwidth,height=0.82\textheight,keepaspectratio]{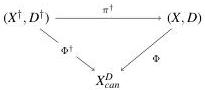}
\end{figure}

\textit{where $\pi^{\dagger}$ is birational and the generic fiber of $\Phi^{\dagger}$ has Kodaira dimension zero.}

\textit{Note that $\kappa(X, D) := \kappa(X^{\dagger}, D^{\dagger})$ and of if $\kappa(X, D) = \dim X - 1$ then the general fibers are elliptic curves, and of if $\kappa(X, D) = \dim X - 2$ then the general fibers are surfaces of log Kodaira dimension zero.}

\begin{enumerate}
\item \textit{Let $\Omega^D$ be a measure associated to a orthonormal base of the}
\begin{equation*}
H^0(X, \mathcal{O}_X(mD) \otimes \mathcal{I}(mh_D))
\end{equation*}
\textit{then the pushforward measure $\Phi_*\Omega^D$ on $X^D_{can}$ is defined by}
\begin{equation*}
\Phi_*\Omega^D = (\Phi^{\dagger})_* \left( (\pi^{\dagger})^* \Omega^D \right)
\end{equation*}

\item \textit{Let $\Phi = \Phi_m$ be the pluricanonical map associated to a basis $\{\sigma_{j_m}\}_{j_m=0}^{d_m}$ of the linear system $| m(K_X + D) |$, associated to the $H^0(X, \mathcal{O}_X(mD) \otimes \mathcal{I}(mh_D))$, for $m$ sufficiently large. Also let $\Omega^D_m = \left( \frac{1}{|S|_{h_D}^{2(1-\beta)}} \sum_{j_m=0}^{d_m} \sigma_{j_m} \otimes \bar{\sigma}_{j_m} \right)^{\frac{1}{m}}$ and $\omega_{FS}$ be the Fubini-Study metric of $\mathbb{CP}^{d_m}$ restricted on $X_{can}$ associated to $\Phi_m$. Then the generalized log Weil-Petersson metric is defined by}
\begin{equation*}
\tilde{\omega}_{WP}^D = \frac{1}{m} \omega_{FS} - \sqrt{-1} \partial \bar{\partial} \log \Phi_* \Omega_m^D
\end{equation*}
\end{enumerate}

In particular, $\tilde{\omega}_{WP}^{D}$ coincides with $\omega_{WP}^{D}$ which introduced before on a Zariski open set of $X_{can}^{D}$

\textbf{Lemma 2.13} \textit{$\Phi_{*}\Omega^{D}$ is independent of the choice of the diagram in Definition 1.4}

\textbf{Proof} Let $\rho$ be a test function on $X_{can}^{D}$. Then
\begin{equation*}
\int_{X_{can}} \rho \Phi_* \Omega^D = \int_{X^{\dagger}} \left((\Phi^{\dagger})^* \rho\right) (\pi^{\dagger})^* \Omega^D = \int_X (\Phi^* \rho) \Omega^D
\end{equation*}
\noindent\hfill $\square$
Since the log generic fiber of $\Phi^{\dagger}$ has log Kodaira dimension zero, we can derive the following lemma. Let us recall Hartogs' extension theorem which is as same as Serre's S2 property in general setting.

\textbf{Lemma 2.14} \textit{Let $f$ be a holomorphic function on a set $G \setminus K$, where $G$ is an open subset of $\mathbb{C}^n$ ($n \geq 2$) and $K$ is a compact subset of $G$. If the relative complement $G \setminus K$ is connected, then $f$ can be extended to a unique holomorphic function on $G$.}

\textbf{Lemma 2.15} \textit{Let $\{\sigma_{j_p}^{(1)}\}_{j_p = 0}^{d_p}$ and $\{\sigma_{j_q}^{(2)}\}_{j_q = 0}^{d_q}$ be bases of the linear systems $|p(K_X + D)|$ and $|q(K_X + D)|$, for $p$ and $q$ sufficiently large. Let $\Omega_p^D = \left(\frac{1}{|S|_{h_D}^{2p(1 - \beta)}}\sum_{j_p = 0}^{d_p}\sigma_{j_p}\otimes \bar{\sigma}_{j_p}\right)^{\frac{1}{p}}$ and $\Omega_q^D = \left(\frac{1}{|S|_{h_D}^{2q(1 - \beta)}}\sum_{j_q = 0}^{d_q}\sigma_{j_q}\otimes \bar{\sigma}_{j_q}\right)^{\frac{1}{q}}$. Then $(\pi^\dagger)^*\left(\frac{\Omega_p^D}{\Omega_q^D}\right)$ is constant on any generic fiber and so}
\begin{equation*}
(\pi^{\dagger})^* \left(\frac{\Omega_{p}^{D}}{\Omega_{q}^{D}}\right) = (\Phi^{\dagger})^* \left(\frac{\Phi_* \Omega_{p}^{D}}{\Phi_* \Omega_{q}^{D}}\right)
\end{equation*}

\textbf{Proof} Take $F = \frac{\Omega_p^D}{\Omega_q^D}$ then $F$ is smooth. We consider the following diagram
\begin{equation*}
\begin{array}{ccc}
(X^{\dagger}, D^{\dagger}) & \xrightarrow{\quad f \quad} & (X, D) \\
\Phi^{\dagger} \bigg\downarrow & & \bigg\downarrow \Phi \\
(Y^{\dagger}, D_Y^{\dagger}) & \xrightarrow{\quad g \quad} & X_{can}^D
\end{array}
\end{equation*}

where $\Phi^{\dagger}$ is an log Iitaka fibration. A very general log fiber of $\Phi^{\dagger}$ is nonsingular of log Kodaira dimension zero. Let $F_{s_0} = ((\Phi^{\dagger})^{-1}(s_0), D_{s_0}^{\dagger})$ be a very general log fiber. Take $B$ as an open neighborhood such that for every $s \in B$ log fiber is nonsingular. Let $\eta$ be a nowhere-vanishing holomorphic $\kappa$-form on $B$. Then
\begin{equation*}
\frac{(\pi^{\dagger})^{*}\sigma_{p,j}}{\eta^{p}}\Big|_{s_{0}}\in H^{0}(X_{s_{0}}^{\dagger},p!(K_{X_{s_{0}}^{\dagger}}+D_{s_{0}}^{\dagger}))
\end{equation*}
and
\begin{equation*}
\frac{(\pi^{\dagger})^{*}\sigma_{q,j}}{\eta^{q}}\Big|_{s_{0}}\in H^{0}(X_{s_{0}}^{\dagger},q!(K_{X_{s_{0}}^{\dagger}}+D_{s_{0}}^{\dagger}))
\end{equation*}

But $\dim \mathbb{P}H^{0}(X_{s_{0}}^{\dagger},k(K_{X_{s_{0}}^{\dagger}}+D_{s_{0}}^{\dagger}))=0$ for any $k\geq 1$ and hence $\frac{\Omega_{p}^{D}}{\Omega_{q}^{D}}$ must be constant on each log fiber $F_{s_{0}}$. Hence $f^{*}F$ is constant on very general log fiber of $\Phi^{\dagger}$. Therefore, by applying Hartogs extension theorem, $f^{*}F$ is smooth on $X^{\dagger}$ and so $f^{*}F$ is the pullback of a function on a Zariski open set of $Y^{\dagger}$. By commutative diagram, on a Zariski open set of $X$, $F$ is the pullback of a function on $X_{can}^{D}$ and so $F$ has to be constant on a very general fiber of $\Phi$. \hfill $\square$

Now we show that our definition of $\tilde{\omega}_{WP}^{D}$ only depends on $(X,D)$.

\textbf{Theorem 2.16} \textit{The definition of generalized logarithmic Weil-Petersson metric $\tilde{\omega}_{WP}^{D}$ only depends on pair $(X,D)$.}

\textbf{Proof} Because $\Phi_{*}\Omega^{D}$ is independent of the choice of the diagram in

\begin{figure}[H]
  \centering
  \includegraphics[width=0.4\textwidth,height=0.82\textheight,keepaspectratio]{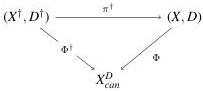}
\end{figure}

so the generalized logarithmic Weil-Petersson metric $\tilde{\omega}_{WP}^{D}$ does not depend on the choice of the diagram. Let $\{\sigma_{j_{p}}^{(1)}\}_{j_{p}=0}^{d_{p}}$ and $\{\sigma_{j_{q}}^{(1)}\}_{j_{q}=0}^{d_{q}}$ are bases of the linear systems $|p(K_{X}+D)|$ and $|q(K_{X}+D)|$, for $p$ and $q$ sufficiently large. Let $\Omega_{p}^{D}=\left(\frac{1}{|S|_{h_D}^{2p(1 - \beta)}}\sum_{j_p = 0}^{d_p}\sigma_{j_p}\otimes \bar{\sigma}_{j_p}\right)^{\frac{1}{p}}$ and $\Omega_q^D = \left(\frac{1}{|S|_{h_D}^{2q(1 - \beta)}}\sum_{j_q = 0}^{d_q}\sigma_{j_q}\otimes \bar{\sigma}_{j_q}\right)^{\frac{1}{q}}$. Assume that $\omega_{FS}^{(1),D}$ and $\omega_{FS}^{(2),D}$ are the log Fubini-Study metrics of $(\mathbb{CP}^{d_p},D)$ and $(\mathbb{CP}^{d_q},D)$ restricted on $X_{can}^D$ associated to $\Phi_p$ and $\Phi_q$. Then by avoiding the base locus of $R(X^{\dagger},K_{X^{\dagger}} + D^{\dagger})$, there exists a Zariski open set $V$ of $X^{\dagger}$, such that
\begin{equation*}
\frac {1}{p} (\Phi^ {\dagger}) ^ {*} \omega_ {F S} ^ {(1), D} = \sqrt {- 1} \partial \bar {\partial} \log (\pi^ {\dagger}) ^ {*} \Omega_ {p} ^ {D}
\end{equation*}

and

\begin{equation*}
\frac {1}{q} (\Phi^ {\dagger}) ^ {*} \omega_ {F S} ^ {(2), D} = \sqrt {- 1} \partial \bar {\partial} \log (\pi^ {\dagger}) ^ {*} \Omega_ {q} ^ {D}
\end{equation*}

and hence there exists a Zariski open set $U$ of $X_{can}$ such that,

\begin{equation*}
\frac {1}{p} \omega_ {F S} ^ {(1), D} - \frac {1}{q} \omega_ {F S} ^ {(2), D} = \sqrt {- 1} \partial \bar {\partial} \log \left(\frac {\Phi_ {*} \Omega_ {p} ^ {D}}{\Phi_ {*} \Omega_ {q} ^ {D}}\right)
\end{equation*}

But $\frac{1}{p}\omega_{FS}^{(1),D}$ and $\frac{1}{q}\omega_{FS}^{(2),D}$ are in the same cohomology class and $\log \left(\frac{\Phi_*\Omega_p^D}{\Phi_*\Omega_q^D}\right)$ is in $L^{\infty}(X_{can}^{D})$, hence previous equality holds true everywhere on $X_{can}^{D}$ and we have $\tilde{\omega}_{WP}^{(1)} = \tilde{\omega}_{WP}^{(2)}$ and proof is complete. \hfill $\Box$

\vspace{1em}
Now, we are ready to prove Theorem 2.11. First we prove part (3).

\vspace{1em}
\noindent\textbf{Proof} Since the log canonical Ring $R(X, K_X + D)$ is finitely generated, hence the following diagram exists

\begin{figure}[H]
\centering
\includegraphics[width=0.45\textwidth,height=0.82\textheight,keepaspectratio]{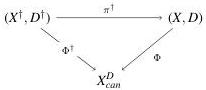}
\end{figure}

\noindent where $X_{can}^{D}$ is the log canonical model of the pair $(X, D)$, and $(X^{\dagger}, D^{\dagger})$ is the log resolution of the stable base locus of the pluricanonical systems such that $(\pi^{\dagger})^{*}(m(K_{X} + D)) = L + E$ for sufficiently large $m$, where $L$ is semi-ample, so it is globally generated and $E$ is a simple normal crossing divisor. Here, $(X^{\dagger}, D^{\dagger})$ is a log Iitaka fibration over the log canonical model $X_{can}^D$ such that the log generic fiber has log Kodaira dimension zero.Let $\{\sigma_j\}_{j=0}^{d_m}$ be a basis of $H^0(X,m(K_X+D)\otimes\mathcal{I}(h_D^m))$ and $\{\zeta_j\}_{j=0}^{d_m}$ be a basis of $H^0(X,(L+E)\otimes\mathcal{I}(h_E))$. Take

\begin{equation*}
\Omega=(\pi^\dagger)^\star(\sum_{n=0}^m\sum_{j_n=0}^{d_n}|\sigma_{j_n}|^{\frac{2}{n}})
\end{equation*}

\noindent be a degenerate smooth volume form on $(X^{\dagger},D^{\dagger})$ and $\omega=\frac{1}{M}\sqrt{-1}\partial\bar{\partial}\log\left(\sum_{j=0}^{d_{m}}|\zeta_{j}|^{2}\right)$. Then the following Monge-Ampère equation has a unique continuous solution $\varphi$ on $X_{can}^{D}$,

\begin{equation*}
(e^{-t}\omega_{WP}^{D}+(1-e^{-t})\omega_{0}+Ric(h_{N})+\sqrt{-1}\partial\bar{\partial}\varphi)^{\kappa}=e^{\varphi}\frac{|S_{E}|_{h_{E}}^{2}}{|S_{D}|_{h_{D}}^{2}}(\Phi^{\dagger})_{*}\Omega
\end{equation*}

\noindent where $\omega_{WP}^{D}$ is the logarithmic Weil-Petersson metric and $N$ is an effective divisor correspond to the canonical bundle formula of Fujino-Mori [72]. So by taking $\sqrt{-1}\partial\bar{\partial}$ on both sides, we get

\begin{equation*}
Ric(\omega_{can})=-\omega_{can}+\omega_{WP}^{D}+[N]
\end{equation*}

\noindent here $\omega_{WP}^{D}=\sqrt{-1}\partial\bar{\partial}\log\Theta-\frac{1}{m}\sqrt{-1}\partial\bar{\partial}\log|S_{E}|^{2}$, where $\Theta=\frac{\Omega}{(\Phi^{\dagger})^{*}(\Phi^{\dagger})_{*}\Omega}$.

Now, by taking

\begin{equation*}
\Omega_{can}^{D}=\frac{e^{\Phi^{*}\varphi}}{h_{N}}\left(\sum_{n=0}^{m}\sum_{j_{n}=0}^{d_{n}}|\sigma_{j_{n}}|^{\frac{2}{n}}\right)
\end{equation*}

\noindent and using the definition of pullback of current, we can see that

\begin{equation*}
\sqrt{-1}\partial\bar{\partial}\log\Omega_{can}=\sqrt{-1}\partial\bar{\partial}\log\left(\sum_{n=0}^{m}\sum_{j_{n}=0}^{d_{n}}|\sigma_{j_{n}}|^{\frac{2}{n}}\right)-\sqrt{-1}\partial\bar{\partial}\log h_{N}+\sqrt{-1}\partial\bar{\partial}\Phi^{*}\varphi=\Phi^{*}\omega_{can}^{D}
\end{equation*}
\noindent where $\omega_{can}^{D}=e^{-t}\omega_{WP}^{D}+(1-e^{-t})\omega_{0}+Ric(h_{N})+\sqrt{-1}\partial\bar{\partial}\varphi$

\noindent Hence

\begin{equation*}
\sqrt{-1}\partial\bar{\partial}\log\Omega_{(X,D)/X_{can}^{D}}=\Phi^{*}\omega_{can}^{D}
\end{equation*}

\noindent Note that, we write the relative canonical volume as

\begin{equation*}
\Omega_ {c a n} ^ {D} := \Omega_ {(X, D) / X _ {c a n} ^ {D}}
\end{equation*}

Now, we prove the uniqueness of $\Omega_{can}^{D}$. Assume that there exists two log canonical measures $\Omega_{can}^{D}$ and $\Omega_{can}^{\prime D}$, then $\Omega_{can}^{\prime D} = e^{\varphi'}\Omega_{can}^{D}$, but $\sqrt{-1}\partial \bar{\partial}\log \Omega_{can}^{\prime D} - \sqrt{-1}\partial \bar{\partial}\log \Omega_{can}^{D}$ is a pullback from $X_{can}^{D}$ then on a generic fiber of $\Phi^{\dagger}$, we have $\sqrt{-1}\partial \bar{\partial}\log \Omega_{can}^{\prime D} - \sqrt{-1}\partial \bar{\partial}\log \Omega_{can}^{D} = 0$ hence $\sqrt{-1}\partial_F\bar{\partial}_F(\pi^\dagger)^*\varphi' = 0$ hence $\varphi'$ descends to $X_{can}^{D}$ and satisfies the following Monge-Ampere equation:

\begin{equation*}
(e ^ {- t} \omega_ {W P} ^ {D} + (1 - e ^ {- t}) \omega_ {0} + R i c (h _ {N}) + \sqrt {- 1} \partial \bar {\partial} \varphi) ^ {\kappa} = e ^ {\varphi^ {\prime}} \frac {| S _ {E} | _ {h _ {E}} ^ {2}}{| S _ {D} | _ {h _ {D}} ^ {2}} (\Phi^ {\dagger}) _ {*} \Omega
\end{equation*}

\noindent but this equation has unique solution, so $\Omega_{can}^{D} = \Omega_{can}^{\prime D}$. Now we show that log canonical measure is invariant under birational transformation.

Consider the following diagram

\begin{figure}[H]
\centering
\includegraphics[width=0.45\textwidth,height=0.82\textheight,keepaspectratio]{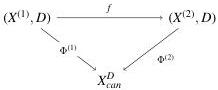}
\end{figure}

\noindent where $f$ is birational and $\Phi^{(1)}$ and $\Phi^{(2)}$ are the log pluricanonical maps. Fix $\Omega^D$ in $(X^{(1)}, D)$. Let $\Omega^{(1)} = e^{\varphi_{(1)}}\Omega^D$ and $\Omega^{(2)} = e^{\varphi_{(2)}}f^*(\Omega^D)$ be the unique canonical measures on $(X^{(1)}, D)$ and $(X^{(2)}, D)$. Then $\varphi_{(1)}$ and $\varphi_{(2)}$ descend to $X_{can}^D$ and satisfy

\begin{equation*}
(e ^ {- t} \omega_ {W P} ^ {D} + (1 - e ^ {- t}) \omega_ {0} + R i c (h _ {N}) + \sqrt {- 1} \partial \bar {\partial} \varphi) ^ {\kappa} = e ^ {\varphi} \frac {| S _ {E} | _ {h _ {E}} ^ {2}}{| S _ {D} | _ {h _ {D}} ^ {2}} (\Phi^ {(1)}) _ {*} \Omega^ {D}
\end{equation*}

\noindent so the uniqueness of the solution of this Monge-Ampere equation implies that $\varphi_{(1)} = \varphi_{(2)}$, and hence $f^{*}\left(\Omega^{(1)}\right) = \Omega^{(2)}$ and so the log canonical measure $\Omega_{can}^{D}$ is invariant under birational transformation $f$. \hfill $\Box$

\vspace{1.5em}
\noindent\textbf{The relation between the Existence of Zariski Decomposition and the Existence of Initial K\"ahler metric along relative K\"ahler Ricci flow:}

\vspace{0.5em}
\noindent Finding an initial Kähler metric $\omega_0$ to run the Kähler Ricci flow is important. Along holomorphic fibration with Calabi-Yau fibres, finding such initial metric is a little bit

\noindent mysterious. In fact, we show that how the existence of initial Kähler metric is related to finite generation of canonical ring along singularities.

Let $\pi : X \to Y$ be an Iitaka fibration of projective varieties $X, Y$, (possibly singular) then is there always the following decomposition

\begin{equation*}
K_Y + \frac{1}{m!} \pi_* \mathcal{O}_X(m! K_{X/Y}) = P + N
\end{equation*}

\noindent where $P$ is semiample and $N$ is effective. The reason is that, If $X$ is smooth projective variety, then as we mentioned before, the canonical ring $R(X, K_X)$ is finitely generated. We may thus assume that $R(X, kK_X)$ is generated in degree 1 for some $k > 0$. Passing to a log resolution of $|kK_X|$ we may assume that $|kK_X| = M + F$ where $F$ is the fixed divisor and $M$ is base point free and so $M$ defines a morphism $f : X \to Y$ which is the Iitaka fibration. Thus $M = f^* \mathcal{O}_Y(1)$ is semiample and $F$ is effective.

In singular case, if $X$ is log terminal. By using Fujino-Mori's higher canonical bundle formula[72], after resolving $X'$, we get a morphism $X' \to Y'$ and a klt pair $K'_{Y} + B'_{Y}$. The Y described above is the log canonical model of $K'_{Y'} + B'_{Y'}$ and so in fact (assuming as above that $Y' \to Y$ is a morphism), then $K'_{Y} + B'_{Y} \sim_{\mathbb{Q}} P + N$ where $P$ is the pull-back of a rational multiple of $\mathcal{O}_Y(1)$ and $N$ is effective (the stable fixed divisor). If $Y' \to Y$ is not a morphism, then P will have a base locus corresponding to the indeterminacy locus of this map.(Thanks of Hacon answer to my Mathoverflow question [57] which is due to E.Viehweg [58])

So the existence of Zariski decomposition is related to the finite generation of canonical ring (when $X$ is smooth or log terminal). Now if such Zariski decomposition exists then, there exists a singular hermitian metric $h$, with semi-positive Ricci curvature $\sqrt{-1}\Theta_h$ on $P$, and it is enough to take the initial metric $\omega_0 = \sqrt{-1}\Theta_h + [N]$ or $\omega_0 = \sqrt{-1}\Theta_h + \sqrt{-1}\delta\partial\bar{\partial}\|S_N\|^{2\beta}$ along relative Kähler Ricci flow

\begin{equation*}
\frac{\partial \omega(t)}{\partial t} = -Ric_{X/Y}(\omega(t)) - \omega(t)
\end{equation*}

\noindent with log terminal singularities.

So when $X, Y$ have at worst log terminal singularities(hence canonical ring is f.g and we have initial Kähler metric to run Kähler Ricci flow with starting metric $\omega_0$) and central fibre is Calabi-Yau variety, and $-K_Y < 0$, then all the fibres are Calabi-Yau varieties and the relative Kähler-Ricci flow converges to $\omega$ which satisfies in

\begin{equation*}
Ric(\omega) = -\omega + f^* \omega_{WP}
\end{equation*}

To get higher estimates of $Ric(\omega) = -\omega + f^* \omega_{WP}$ We need just to think about central fiber to control singularities of general fibers and total space. In fact from Kawamata,

\noindent Let $f: X \to S$ be a flat morphism from a germ of an algebraic variety to a nonsingular germ of an algebraic variety, and suppose that the central fiber $X_0 = f^{-1}(0)$ has at worst canonical singularities (respectively, terminal singularities). Then $X$ itself has at worst canonical singularities (respectively, terminal singularities). In particular, every fiber $X_s = f^{-1}(s)$ has at worst canonical singularities (respectively, terminal singularities).

For the existence of Kähler-Einstein metric when our variety is of general type, we need to the nice deformation of Kähler-Ricci flow and for intermediate Kodaira dimension we need to work on relative version of Kähler Ricci flow. i.e

\begin{equation*}
\frac {\partial \omega}{\partial t} = - R i c _ {X / Y} (\omega) - \omega
\end{equation*}

\noindent take the reference metric as $\tilde{\omega}_t = e^{-t}\omega_0 + (1 - e^{-t})Ric\left(\frac{\omega_{SRF}^n\wedge\pi^*\omega_{can}^m}{\pi^*\omega_{can}^m}\right)$ then the version of Kähler Ricci flow is equivalent to the following relative Monge-Ampere equation

\begin{equation*}
\frac {\partial \phi_ {t}}{\partial t} = \log \frac {(\tilde {\omega} _ {t} + \sqrt {- 1} \partial \bar {\partial} \phi_ {t}) ^ {n} \wedge \pi^ {*} \omega_ {c a n} ^ {m}}{\omega_ {S R F} ^ {n} \wedge \pi^ {*} \omega_ {c a n} ^ {m}} - \phi_ {t}
\end{equation*}

\noindent Take the relative canonical volume form $\Omega_{X / Y} = \frac{\omega_{SRF}^n\wedge\pi^*\omega_{can}^m}{\pi^*\omega_{can}^m}$ and $\omega_{t} = \tilde{\omega}_{t} + \sqrt{-1}\partial \bar{\partial}\phi_{t}$, then

\begin{equation*}
\frac {\partial \omega_ {t}}{\partial t} = \frac {\partial \tilde {\omega} _ {t}}{\partial t} + \sqrt {- 1} \partial \bar {\partial} \frac {\partial \phi_ {t}}{\partial t}
\end{equation*}

\noindent By taking $\omega_{\infty} = -Ric(\Omega_{X / Y}) + \sqrt{-1}\partial \bar{\partial}\phi_{\infty}$ we obtain after using estimates

\begin{equation*}
\log \frac {\omega_ {\infty} ^ {n}}{\Omega_ {X / Y}} - \phi_ {\infty} = 0
\end{equation*}

\noindent By taking $-\sqrt{-1}\partial\bar{\partial}$ of both sides we get

\begin{equation*}
R i c _ {X / Y} (\omega_ {\infty}) = - \omega_ {\infty}
\end{equation*}

\noindent hence by the definition of relative Kähler-metric and higher canonical bundle formula we have the Song-Tian metric[23, 24, 25, 26]

\begin{equation*}
R i c (\omega_ {\infty}) = - \omega_ {\infty} + \pi^ {*} (\omega_ {W P}) + [ N ]
\end{equation*}

\vspace{1em}
\noindent\textbf{Uniqueness result of Relative Kähler-Einstein metric:} Uniqueness of the solutions of relative Kahler Ricci flow along Iitaka fibration or $\pi : X \to X_{can}$ or along log canonical

model $\pi : (X, D) \to X_{can}^{D}$. Let $\phi_0$ and $\psi_0$ be be $\omega$-plurisubharmonic functions such that $\mathfrak{v}(\phi_0, x) = 0$ for all $x \in X$, let $\phi_t$, and $\psi_t$ be the solutions of relative Kähler Ricci flow starting from $\phi_0$ and $\psi_0$, respectively. Then in [48] it has been proven that if $\phi_0 < \psi_0$ then $\phi_t < \psi_t$ for all $t$. In particular, the flow is unique. So from the deep result of Tsuji-Schumacher[40], it has been shown that Weil-Petersson metric has zero Lelong number on moduli space of Calabi-Yau varieties, and by the same method we can show that logarithmic Weil-Petersson metric has zero Lelong number on moduli space of log Calabi-Yau varieties, hence by taking the initial metric to be Weil-Petersson metric or logarithmic Weil-Petersson metric and since Weil-Petersson metric or logarithmic Weil-Petersson metric are Kahler and semi-positive hence we get the uniquenessof the solutions of relative Kähler Ricci flow.

\textbf{Conjecture:} The twisted Kähler-Einstein metric $Ric(\omega) = -\omega + \alpha$ where $\alpha$ is a semi-positive current has unique solution if and only if $\alpha$ has zero Lelong number

For the log-Calabi-Yau fibration $f : (X, D) \to Y$, such that the fibers $(X_{t}, D_{t})$ are log Calabi-Yau varieties, if $(X, \omega)$ be a Kähler variety with Poincaré singularities then the semi-Ricci flat metric has $\omega_{SRF}|_{X_{t}}$ is quasi-isometric with the following model which we call it \textbf{fibrewise Poincaré singularities.}

\begin{equation*}
\frac {\sqrt {- 1}}{\pi} \sum_ {k = 1} ^ {n} \frac {d z _ {k} \wedge d \bar {z} _ {k}}{| z _ {k} | ^ {2} (\log | z _ {k} | ^ {2}) ^ {2}} + \frac {\sqrt {- 1}}{\pi} \frac {1}{(\log | t | ^ {2} - \sum_ {k = 1} ^ {n} \log | z _ {k} | ^ {2}) ^ {2}} \left(\sum_ {k = 1} ^ {n} \frac {d z _ {k}}{z _ {k}} \wedge \sum_ {k = 1} ^ {n} \frac {d \bar {z} _ {k}}{\bar {z} _ {k}}\right)
\end{equation*}

We can define the same \textbf{fibrewise conical singularities} and the semi-Ricci flat metric has $\omega_{SRF}|_{X_t}$ is quasi-isometric with the following model

\begin{equation*}
\frac {\sqrt {- 1}}{\pi} \sum_ {k = 1} ^ {n} \frac {d z _ {k} \wedge d \bar {z} _ {k}}{| z _ {k} | ^ {2}} + \frac {\sqrt {- 1}}{\pi} \frac {1}{(\log | t | ^ {2} - \sum_ {k = 1} ^ {n} \log | z _ {k} | ^ {2}) ^ {2}} \left(\sum_ {k = 1} ^ {n} \frac {d z _ {k}}{z _ {k}} \wedge \sum_ {k = 1} ^ {n} \frac {d \bar {z} _ {k}}{\bar {z} _ {k}}\right)
\end{equation*}

\textbf{Remark:}Note that the log semi-Ricci flat metric $\omega_{SRF}^{D}$ is not continuous in general. But if the central fiber has at worst canonical singularities and the central fiber $(X_{0}, D_{0})$ be itself as Calabi-Yau pair, then by open condition property of Kahler-Einstein metrics, semi-Ricci flat metric is smooth in an open Zariski subset.

\textbf{Remark:}So by applying the previous remark, the relative volume form

\begin{equation*}
\Omega_ {(X, D) / Y} = \frac {(\omega_ {S R F} ^ {D}) ^ {n} \wedge \pi^ {*} \omega_ {c a n} ^ {m}}{\pi^ {*} \omega_ {c a n} ^ {m} | S | ^ {2}}
\end{equation*}

is not smooth in general, where $S \in H^0(X, L_N)$ and $N$ is a divisor which come from canonical bundle formula of Fujino-Mori [72].

Now we try to extend the Relative Ricci flow to the fiberwise conical relative Ricci flow. We define the conical Relative Ricci flow on pair $\pi : (X, D) \to Y$ where $D$ is a simple normal crossing divisor as follows

\begin{equation*}
\frac {\partial \omega}{\partial t} = - R i c _ {(X, D) / Y} (\omega) - \omega + [ N ]
\end{equation*}

where $N$ is a divisor which come from canonical bundle formula of Fujino-Mori[72].

Take the reference metric as $\tilde{\omega}_t = e^{-t}\omega_0 + (1 - e^{-t})Ric\left(\frac{\omega_{SRF}^n\wedge\pi^*\omega_{can}^m}{\pi^*\omega_{can}^m}\right)$ then the conical relative Kähler Ricci flow is equivalent to the following relative Monge-Ampere equation

\begin{equation*}
\frac {\partial \phi_ {t}}{\partial t} = \log \frac {(\tilde {\omega} _ {t} + R i c (h _ {N}) + \sqrt {- 1} \partial \bar {\partial} \phi_ {t}) ^ {n} \wedge \pi^ {*} \omega_ {c a n} ^ {m} | S _ {N} | ^ {2}}{(\omega_ {S R F} ^ {D}) ^ {n} \wedge \pi^ {*} \omega_ {c a n} ^ {m}} - \phi_ {t}
\end{equation*}

With cone angle $2\pi \beta$, $(0 < \beta < 1)$ along the divisor $D$, where $h$ is an Hermitian metric on line bundle corresponding to divisor $N$, i.e., $L_N$. This equation can be solved. Take, $\omega = \omega(t) = \omega_B + (1 - \beta)Ric(h) + \sqrt{-1}\partial \bar{\partial} v$ where $\omega_B = e^{-t}\omega_0 + (1 - e^{-t})Ric(\frac{(\omega_{SRF}^D)^n\wedge\pi^*\omega_{can}^m}{\pi^*\omega_{can}^m}) = e^{-t}\omega_0 + (1 - e^{-t})\omega_{WP}$, by using Poincare-Lelong equation,

\begin{equation*}
\sqrt {- 1} \partial \bar {\partial} \log | s _ {N} | _ {h} ^ {2} = - c _ {1} (L _ {N}, h) + [ N ]
\end{equation*}

we have

\begin{align*}
Ric(\omega) &= \\
&= - \sqrt {- 1} \partial \bar {\partial} \log \omega^ {m} \\
&= - \sqrt {- 1} \partial \bar {\partial} \log \pi_ {*} \Omega_ {(X, D) / Y} - \sqrt {- 1} \partial \bar {\partial} v - (1 - \beta) c _ {1} ([ N ], h) + (1 - \beta) \{N \}
\end{align*}

and

\begin{align*}
\sqrt {- 1} \partial \bar {\partial} \log \pi_ {*} \Omega_ {(X, D) / Y} + \sqrt {- 1} \partial \bar {\partial} v &= \\
&= \sqrt {- 1} \partial \bar {\partial} \log \pi_ {*} \Omega_ {(X, D) / Y} + \omega - \omega_ {B} - R i c (h)
\end{align*}

Hence, by using

\begin{equation*}
\omega_ {W P} ^ {D} = \sqrt {- 1} \partial \bar {\partial} \log (\frac {(\omega_ {S R F} ^ {D}) ^ {n} \wedge \pi^ {*} \omega_ {c a n} ^ {m}}{\pi^ {*} \omega_ {c a n} ^ {m} | S | ^ {2}})
\end{equation*}

we get

\begin{align*}
\sqrt {- 1} \partial \bar {\partial} \log \pi_ {*} \Omega_ {(X, D) / Y} + \sqrt {- 1} \partial \bar {\partial} v &= \\
&= \omega - \omega_ {W P} ^ {D} - (1 - \beta) c _ {1} (N)
\end{align*}

So,

\begin{equation*}
Ric(\omega) = - \omega + \omega_ {W P} ^ {D} + (1 - \beta) [ N ]
\end{equation*}

which is equivalent with

\begin{equation*}
Ric_{(X, D) / Y} (\omega) = - \omega + [ N ]
\end{equation*}

Consider the following followwing relative Kähler Ricci flow

\begin{equation*}
\frac {\partial \omega}{\partial t} = - R i c _ {(X, D) / Y} (\omega) - \omega + [ N ]
\end{equation*}

where $N$ is a divisor which come from canonical bundle formula of Fujino-Mori.[72]

Take the reference metric as $\tilde{\omega}_t = e^{-t}\omega_0 + (1 - e^{-t})Ric\left(\frac{\omega_{SRF}^n\wedge\pi^*\omega_{can}^m}{\pi^*\omega_{can}^m}\right)$ then the conical relative Kähler Ricci flow is equivalent to the following relative Monge-Ampere equation

\begin{equation*}
\frac {\partial \phi_ {t}}{\partial t} = \log \frac {(\tilde {\omega} _ {t} + R i c (h _ {N}) + \sqrt {- 1} \partial \bar {\partial} \phi_ {t}) ^ {n} \wedge \pi^ {*} \omega_ {c a n} ^ {m} | S _ {N} | ^ {2}}{(\omega_ {S R F} ^ {D}) ^ {n} \wedge \pi^ {*} \omega_ {c a n} ^ {m}} - \phi_ {t}
\end{equation*}

Now we prove the $C^0$-estimate for this relative Monge-Ampere equation due to Tian's $C^0$-estimate

By approximation our Monge-Ampere equation, we can write

\begin{equation*}
\frac {\partial \varphi_ {\epsilon}}{\partial t} = \log \frac {(\omega_ {t , \epsilon} + \sqrt {- 1} \partial \bar {\partial} \varphi_ {t}) ^ {m} \wedge \pi^ {*} \omega_ {c a n} ^ {n} (| | S | | ^ {2} + \epsilon^ {2}) ^ {(1 - \beta)}}{(\omega_ {S R F} ^ {D}) ^ {m} \wedge \pi^ {*} \omega_ {c a n} ^ {n}} - \delta (| | S | | ^ {2} + \epsilon^ {2}) ^ {\beta} - \varphi_ {t, \epsilon}
\end{equation*}

So by applying maximal principle we get an upper bound for $\varphi_{t,\epsilon}$ as follows

\begin{equation*}
\frac {\partial}{\partial t} \sup \varphi_ {\epsilon} \leq \sup \log \frac {\omega_ {t , \epsilon} ^ {m} \wedge \pi^ {*} (\omega_ {c a n}) ^ {n} \left(| | S | | ^ {2} + \epsilon^ {2}\right) ^ {(1 - \beta)}}{(\omega_ {S R F} ^ {D}) ^ {m} \wedge \pi^ {*} \omega_ {c a n} ^ {n}} - \delta \left(| | S | | ^ {2} + \epsilon^ {2}\right) ^ {\beta}
\end{equation*}

and by expanding $\omega_{t,\epsilon}^{n}$ we have a constant $C$ independent of $\epsilon$ such that the following expression is bounded if and only if the Song-Tian-Tsuji measure be bonded, so to get $C^{0}$ estimate we need special fiber has mild singularities in the sense of MMP

\begin{equation*}
\frac {\omega_ {t , \epsilon} ^ {m} \wedge \pi^ {*} (\omega_ {c a n}) ^ {n} \left(| | S | | ^ {2} + \epsilon^ {2}\right) ^ {(1 - \beta)}}{(\omega_ {S R F} ^ {D}) ^ {m} \wedge \pi^ {*} \omega_ {c a n} ^ {n}} \approx C
\end{equation*}

and also $\delta \to 0$ so $\delta \left(||S||^2 +\epsilon^2\right)^\beta$ is too small. So we can get a uniform upper bound for $\varphi_{\epsilon}$. By applying the same argument for the lower bound, and using maximal principle again, we get a $C^0$ estimate for $\varphi_{\epsilon}$. Moreover if central fiber $X_0$ has canonical singularities then Song-Tian-Tsuji measure is continuous. So this means that we have $C^0$-estimate for relative Kähler-Ricci flow if and only if the central fiber has at worst canonical singularities.

Note that to get $C^\infty$-estimate we need just check that our reference metric is bounded and Song-Tian-Tsuji measure is $C^\infty$ smooth. So it just remains to see that $\omega_{WP}$ is bounded. But when fibers are not smooth in general, Weil-Petersson metric is not bounded and Yoshikawa in Proposition 5.1 in [39] and Tian[69] showed that under the some additional condition when central fiber $X_0$ is reduced and irreducible and has only canonical singularities we have

\begin{equation*}
0 \leq \omega_ {W P} \leq C \frac {\sqrt {- 1} | s | ^ {2 r} d s \wedge d \bar {s}}{| s | ^ {2} (- \log | s |) ^ {2}}
\end{equation*}

So we can get easily by ancient method! the $C^{\infty}$-solution.

Note that the main difficulty of the solution of $C^\infty$ for the solution of relative Kähler-Einstein metric is that the null direction Vafa-Yau semi Ricci flat metric $\omega_{SRF}$ gives a following Tsuji's foliation [62] along Iitaka fibration $\pi : X \to Y$ and we call it fiberwise Calabi-Yau foliation and can be defined as follows

\begin{equation*}
\mathcal {F} = \{\theta \in T X | \omega_ {S R F} (\theta , \bar {\theta}) = 0 \}
\end{equation*}

and along log Iitaka fibration $\pi:(X,D)\to Y$, we can define the following foliation

\begin{equation*}
\mathcal {F} ^ {\prime} = \{\theta \in T X ^ {\prime} | \omega_ {S R F} ^ {D} (\theta , \bar {\theta}) = 0 \}
\end{equation*}

where $X' = X \setminus D$. In fact the method of Song-Tian (only in fiber direction and they couldn't prove the estimates in horizontal direction which is the main part of computation) works when $\omega_{SRF} > 0$. For the null direction we need to an extension of Monge-Ampere foliation method of Gang Tian in [51]. It will appear in my new paper [52]

A complex analytic space is a topological space such that each point has an open neighborhood homeomorphic to some zero set $V(f_1, \ldots, f_k)$ of finitely many holomorphic functions in $\mathbb{C}^n$, in a way such that the transition maps (restricted to their appropriate domains) are biholomorphic functions.

\textbf{Lemma:} Fiberwise Calabi-Yau foliation is complex analytic space and its leaves are also complex analytic spaces. See [52][62]

\textbf{Lemma:}(See [52][62]) Let $\mathcal{L}$ be a leaf of $f_*\mathcal{F}'$, then $\mathcal{L}$ is a closed complex submanifold and the leaf $\mathcal{L}$ can be seen as fiber on the moduli map

\begin{equation*}
\eta : \mathcal{Y} \to \mathcal{M}_{can}^D
\end{equation*}

where $\mathcal{M}_{can}^D$ is the moduli space of log calabi-Yau fibers with at worst canonical singularites and

\begin{equation*}
\mathcal{Y} = \{y \in Y_{reg} | (X_y, D_y) \text{ is Kawamata log terminal pair}\}
\end{equation*}

Now we apply Theorem 1.12. when $(X, D)$ is a minimal elliptic surface of $\text{kod}(X, D) = 1$. Let first recall Kodaira's canonical bundle formula for minimal elliptic surfaces.

Now we can have the following definition of KE-stability due to Takayama-Tosatti and Wang [67] (when dimension of base is one) and when dimension of base is bigger than one, then Takayama-Tosatti result does not functional and we need to consider $C^\infty$ smoothness and boundedness of Song-Tian-Tsuji measure as follows which seems to be natural.

Now if the dimension of smooth base be bigger than one, then we don't know the semi-stable reduction and instead we can use weak Abramovich-Karu reduction [66] or Kawamata's unipotent reduction theorem[42]. In fact when the dimension of base is one we know from Fujino's recent result[73] that if we allow semi-stable reduction and MMP on the family of Calabi-Yau varieties then the central fiber will be Calabi-Yau variety. But If the dimension of smooth base be bigger than one on the family of Calabi-Yau fibers, then if we apply MMP and weak Abramovich-Karu semi-stable reduction [31] then the special fiber can have simple nature as stable variety with slc singularities at worst. But if the dimension of base be singular then we don't know


about semi-stable reduction which seems is very important for finding canonical metric along Iitaka fibration.

\textbf{Definition 2.17} Let $\pi : X \to B$ be a family of K\"ahler-Einstein varieties, then we introduce the new notion of stability and call it fiberwise KE-stability, if the Weil-Petersson distance $d_{WP}(B,0) < \infty$ (or boundedness of Song-Tian-Tsuji measure when the dimension of base is bigger than one) and also Song-Tian-Tsuji measure be $C^\infty$. Note that when central fiber has at worst canonical singularities then Song-Tian-Tsuji measure is $C^0$-continuous.

Note when fibers are Calabi-Yau varities, Takayama and Tosatti, by using Tian's K\"ahler-potential for Weil-Petersson metric for moduli space of Calabi-Yau varieties showed that Fiberwise KE-Stability is as same as when the central fiber is Calabi-Yau variety with at worst canonical singularities(after base change and birational modification ). Note that when the dimension of base is bigger than one, we can replace the boundedness of Weil-Petersson distance with boundedness of Song-Tian-Tsuji measure.

So along canonical model $\pi : X \to X_{can}$ for mildly singular variety $X$, we have $Ric(\omega) = -\omega + \omega_{WP}$ if and only if our family of fibers be fiberwise KE-stable

\textbf{Theorem 2.18} (\textit{Song-Tian [10]}) Let $f : (X, D) \to \Sigma$ be a minimal elliptic surface such that its multiple fibres are $(X_{s_1} = m_1F_1, D_{s_1}), (X_{s_2} = m_2F_2, D_{s_2}), \ldots, (X_{s_k} = m_kF_k, D_{s_k})$.Then
\begin{equation*}
K_X + D = f^*\left(L \otimes \mathcal{O}_X \left(\sum (1 - m_i) F_i\right)\right)
\end{equation*}
where $L$ is a line bundle on $\Sigma$

Song and Tian [10], by using Kodaira's canonical bundle formula for minimal elliptic surface $f : (X, D) \to \Sigma$ showed that there exists a unique twisted K\"ahler Einstein metric $\omega_\infty$ on regular part of $\Sigma$ such that
\begin{equation*}
Ric(\omega_\infty) = -\omega_\infty + \omega_{WP}^D + \sum \left(1 - \frac{1}{m_i}\right)[D_i]
\end{equation*}
So, if in Theorem 0.8, we take $(1 - \beta)[D] = \sum (1 - \beta_i)[D_i]$ and $\beta_i = \frac{1}{m_i}$ and using Kodaira's canonical bundle formula, then we see that we have Song-Tian's theorem about twisted K\"ahler Einstein metric on minimal elliptic surfaces.


\textbf{Theorem 2.19} \textit{Let $f : (X, D) \to \Sigma$ be a minimal elliptic surface of log Kodaira dimension $\text{kod}(X, D) = 1$ such that its multiple fibres are $(X_{s_1} = m_1F_1, D_{s_1}), (X_{s_2} = m_2F_2, D_{s_2}), \ldots, (X_{s_k} = m_kF_k, D_{s_k})$. Then for any initial conical K\"ahler metric, the conical K\"ahler-Ricci flow has global solution $\omega(t, .)$ and there exists a unique twisted K\"ahler Einstein metric $\omega_\infty$ on regular part of $\Sigma$ such that}
\begin{equation*}
Ric(\omega_\infty) = -\omega_\infty + \omega_{WP}^D + \sum \left(1 - \frac{1}{m_i}\right)[D_i]
\end{equation*}
Moreover, Let $Z$ be a smooth projective complex threefold that admits an abelian surfaces fibration onto a curve $C$, i.e $\pi : Z \to C$. Then
\begin{equation*}
K_Z = \pi^* \left( K_C + \sum_{p \in B} \left(1 - \frac{1}{m_p}\right) p + M_C + R_C \right)
\end{equation*}
where $M_C$, called moduli part correspond to $j^*\mathcal{O}_{\mathbb{P}^1}(1)$ where $j$ is the $j$-function associated to the elliptic fibration and $\sum_{p\in B}\left(1 - \frac{1}{m_p}\right)p + R_C$ is called the discriminant part.

and in this case we have also
\begin{equation*}
Ric(\omega_\infty) = -\omega_\infty + \omega_{WP}^D + \sum \left(1 - \frac{1}{m_i}\right)[D_i]
\end{equation*}

\textbf{Example:} A Kodaira surface is defined to be a smooth compact complex surface X admitting a Kodaira fibration, that is, there exists a connected fibration $\pi: X \to C$ over a smooth compact Riemann surface C such that $\pi$ is everywhere of maximal rank, and the associated Kodaira-Spencer map $\rho_{t}: T_{t}C \to H^{1}(X_{t}, TX_{t})$ at each point $t \in C$ is non-zero. Here, $X_{t} := \pi^{-1}(t)$ denotes the fiber of $\pi$ at t. It is known that each fiber $X_{t}$ is necessarily a smooth compact Riemann surface of genus $g \geq 2$

By identifying a point $x$ in a fiber $X_t$ with the corresponding punctured Riemann surface $X_t \setminus \{x\}$, one obtains a map from $X$ to the moduli space $\mathcal{M}_{g,1}$ of punctured Riemann surfaces of type $(g, 1)$. The map lifts to local immersions to the corresponding Teichmuller space $\mathcal{T}_{g,1}$, and this enables us to induce a metric on $X$ from the Weil-Petersson metric $g_{WP}$ on $\mathcal{T}_{g,1}$.

So on Kodaira fibration surface $f : X \to C$ with C has negative first Chern class then on total space we have
\begin{equation*}
Ric(\omega) = -\omega + f^*\omega_{WP}
\end{equation*}


which Weil-Petersson metric is a metric on moduli space $\mathcal{M}_{g,1}$. We have the same result also for log Kodaira fibration surface

\textbf{Remark:} In the papers of Song-Tian [10],[11], they asked that if for holomorphic submersion $\pi : X \to B$ fibers and base are of general type then we can apply the same method of Song-Tian to get twisted K\"ahler Einstein metric on $X$, but the fact is that we don't need such twisted metric. Let $\pi : X \to B$ be a holomorphic submersion where $B$ and the general fibers of $\pi$ are of general type. Then $X$ is of general type due to E. Viehweg and hence we have unique K\"ahler-Einstein metric on $X$ with negative Ricci curvature. In fact this is a special case of Iitaka conjecture, which was proved by Viehweg, see [33]

Now, we try to obtain the relation between first Chern class of $X$ and cohomology class of $[\pi^{*}\omega_{B}]$. Lets first mention the Grauert's theorem which will be essential for deriving our next theorem.

\textbf{Theorem 2.20} \textit{Let $\pi : (X, D) \to B$ is a holomorphic submersion map with connected fibers and for every $y \in B$, $(X_y, D_y)$ are log Calabi-Yau manifolds , then the relative log pluri-canonical bundle $L = \pi_*((K_{X/B} + D)^n)$ is a Hermitian semi-positive Line bundle on $B$ and logarithmic Weil-Petersson form satisfies in $[\omega_{WP}^D] = \frac{1}{m} c_1(L)$ where $m$ is a smallest positive integer such that $(K_{X_y} + D_y)^{\otimes m}$ be trivial}

\textbf{Theorem 2.21} \textit{Let $\pi : (X, D) \to B$ is a holomorphic submersion map with connected fibers and for every $y \in B$, $\pi^{-1}(y) = (X_y, D_y)$ are log Calabi-Yau manifolds and $c_1(B) < 0$. Then we have}
\begin{equation*}
[\pi^*\omega_B] = -c_1(X, D) + (1 - \beta)\left(\sum_P (b(1 - t_P^D))[\pi^*(P)] + [B^D]\right)
\end{equation*}
\textbf{Proof} By the definition of relative pluri-canonical bundle, we have
\begin{equation*}
\pi_*((K_X + D)^m) = \pi_*((K_{X/B} + D)^m) \otimes K_B^m
\end{equation*}
and by using $K_X^m = \pi^*\pi_*((K_X^m)$ and Grauert's theorem we conclude that $c_1(X, D) = \pi^*c_1(B) - [\pi^*\omega_{WP}^D]$ and we have $\pi^*c_1(\omega_B) = -[\pi^*\omega_B] + [\pi^*\omega_{WP}^D] + (1 - \beta)[D]$ so, by combining these two relations we get the desired result. \hfill $\square$

\textbf{Theorem 2.22} \textit{The maximal time existence $T$ for the solutions of relative K\"ahler Ricci flow is}
\begin{equation*}
T = \sup \left\{t \mid e^{-t}[\omega_0] + (1 - e^{-t})c_1(K_{X/Y} + D) \in \mathcal{K}((X, D)/Y) \right\}
\end{equation*}
where $\mathcal{K}((X,D)/Y)$ denote the relative K\"ahler cone of $f : (X, D) \to Y$


Now we give a motivation that why the geometry of pair $(X, D)$ must be interesting. The first one comes from algebraic geometry, in fact for deforming the cone angle we need to use of geometry of pair $(X, D)$. In the case of minimal general type manifold the canonical bundle of $X$, i.e., $K_X$ is nef and we would like $K_X$ to be ample and it is not possible in general and what we can do is that to add a small multiple of ample bundle $\frac{1}{m} A$, i.e., $K_X + \frac{1}{m} A$ and then we are deal with the pair $(X, \frac{1}{m} H)$ which $H$ is a generic section of it. The second one is the works of Chen-Sun-Donaldson and Tian on existence of K\"ahler Einstein metrics for Fano varieties which they used of geometry of pair $(X, D)$ for their proof.

On pair $(X, D)$, when the first Chern class is definite we can formulate K\"ahler-Einstein metric as
\begin{equation*}
Ric(\omega) = \lambda \omega + [D]
\end{equation*}
On pair $(X, D)$, when log Kodaira dimension is positive and first Chern class has no sign we \textbf{can not} formulate the generalized K\"ahler-Einstein metric on holomorphic fiber space $\pi : (X, D) \to X_{can}^D$ as
\begin{equation*}
Ric(\omega) = \lambda \omega + \alpha + [D]
\end{equation*}
we must explicitly find the additional term along singularities which comes from fibers. In fact Song-Tian showed that along holomorphic fiber space $X \to X_{can}$ the generalized K\"ahler Einstein metric is
\begin{equation*}
Ric(\omega) = -\omega + \omega_{WP}
\end{equation*}
We can expect that on pair $(X, D)$ we have
\begin{equation*}
Ric(\omega) = -\omega + \omega_{WP} + [D]
\end{equation*}
but it is \textbf{not true}

The fact is that the Weil-Petersson metric $\omega_{WP}$ changes to logarithmic Weil-Petersson metric $\omega_{WP}^D$ and current of integration $[D]$ also must be replaced to something else and we will find it explicitly.

Now we explain Tian-Yau program to how to construct model metrics in general, like conical model metric, Poincare model metric, or Saper model metric,....

Tian-Yau program: Let $\mathbb{C}^n = \mathbb{C}^n(z_1, ..., z_n)$ be a complex Euclidian space for some $n > 0$. For a positive number $\epsilon$ with $0 < \epsilon < 1$ consider
\begin{equation*}
\overline{X} = \overline{X}_\epsilon = \left\{ z = (z_1, ..., z_n) \in \mathbb{C}^n \mid |z_i| < \epsilon \right\}
\end{equation*}


Now, let $D_i = \{z_i = 0\}$ be the irreducible divisors and take $D = \sum_i D_i$ where
\begin{equation*}
D = \{ z \in \overline{X} \mid z_1 z_2 \dots z_k = 0 \}
\end{equation*}
and take $X = \overline{X} \setminus D$. In polar coordinate we can write $z_i = r_i e^{i\theta_i}$. Let $g$ be a K\"ahler metric on $D$ such that the associated K\"ahler form $\omega$ is of the following form
\begin{equation*}
\omega = \sqrt{-1} \sum_i \frac{1}{|dz_i|^2} dz_i \wedge d\bar{z}_i
\end{equation*}
Then the volume form $dv$ associated to $\omega$ is written in the form;
\begin{equation*}
dv = (\sqrt{-1})^n \prod_{i=1}^n \frac{1}{|dz_i|^2} \prod_i dz_i \wedge d\bar{z}_i \ , \ v = \frac{1}{|dz_i|^2}
\end{equation*}
Let $L$ be a (trivial) holomorphic line bundle defined on $\overline{X}$, with a generating holomorphic section $S$ on $\overline{X}$. Fix a $C^\infty$ hermitian metric $h$ of $L$ over $X$ and denote by $|S|^2$ the square norm of $S$ with respect to $h$. Assume the functions $|S|^2$ and $|dz_i|^2$ depend only on $r_i$, $1 \leq i \leq k$. Set
\begin{equation*}
d(r_1, \dots, r_k) = |S|^2 \cdot v \cdot \prod_{1 \leq i \leq k} r_i
\end{equation*}
and further make the following three assumptions:

\textbf{A1)} The function $d$ is of the form
\begin{equation*}
d(r_1, \dots, r_k) = r_1^{c_1} \dots r_k^{c_k} (\log 1/r_1)^{b_1} \dots (\log 1/r_k)^{b_k} L(r_1, \dots, r_k)^t
\end{equation*}
where
\begin{equation*}
L = L(r_1, \dots, r_k) = \sum_{i=1}^k \log 1/r_i
\end{equation*}
and $c_i, b_j, t$ are real numbers with $t \geq 0$ such that $q_i = b_i + t \neq -1$ if $c_i$ is an odd integer. We set $a_i = (c_i + 1)/2$ and denote by $[a_i]$ the largest integer which does not exceed $a_i$.

\textbf{A2)} If $1 \leq i \leq k$, then $|dz_i|^2$ is either of the following two forms;
\begin{equation*}
|dz_i|^2(r) = r_i^2 (\log 1/r_i)^2, \text{ or } |dz_i|^2(r) = r_i^2 L^2
\end{equation*}

In fact, A2) implies that the Kähler metric $g$ is (uniformly) complete along $D$.

\textbf{A3)} If $k + 1 \leq i \leq n$, then $|dz_i|^{-2}$ is bounded (above) on $X$.

Now, we give some well-known examples of Tian-Fujiki picture, i.e, conical model metric, Poincare model metric, and Saper model metric.

A Kähler current $\omega$ is called a conical Kähler metric (or Hilbert Modular type) with angle $2\pi\beta$, ($0 < \beta < 1$) along the divisor $D$, if $\omega$ is smooth away from $D$ and asymptotically equivalent along $D$ to the model conic metric

\begin{equation*}
\omega_\beta = \sqrt{-1} \left( \frac{dz_1 \wedge d\bar{z}_1}{|z_1|^{2(1-\beta)}} + \sum_{i=2}^n dz_i \wedge d\bar{z}_i \right)
\end{equation*}

here $(z_1, z_2, \dots, z_n)$ are local holomorphic coordinates and $D = \{z_1 = 0\}$ locally.

After an appropriate -singular- change of coordinates, one can see that this model metric represents an Euclidean cone of total angle $\theta = 2\pi\beta$, whose model on $\mathbb{R}^2$ is the following metric: $d\theta^2 + \beta^2 dr^2$. The volume form $V$ of a conical Kähler metric $\omega_D$ on the pair $(X, D)$ has the form

\begin{equation*}
V = \prod_j |S_j|^{2\beta_j - 2} e^f \omega^n
\end{equation*}

where $f \in C^0$.

This asymptotic behaviour of metrics can be translated to the second order asymptotic behaviour of their potentials

\begin{equation*}
\omega_\beta = -\sqrt{-1} \partial \bar{\partial} \log e^{-u}
\end{equation*}

where $u = \frac{1}{2} \left( \frac{1}{\beta^2} |z_1|^{2\beta} + |z_2|^2 + \dots + |z_n|^2 \right)$.

Moreover, if we let $z = re^{i\theta}$ and $\rho = r^\beta$ then the model metric in $\omega_\beta$ becomes

\begin{equation*}
(d\rho + \sqrt{-1}\beta\rho d\theta) \wedge (d\rho - \sqrt{-1}\beta\rho d\theta) + \sum_{i>1} dz_i \wedge d\bar{z}_i
\end{equation*}

and if we set $\epsilon = e^{\sqrt{-1}\beta\theta} (d\rho + \sqrt{-1}\beta\rho d\theta)$ then the conical Kähler metric $\omega$ on $(X, (1-\beta)D)$ can be expressed as

\begin{equation*}
\omega = \sqrt{-1} \left( f\epsilon \wedge \bar{\epsilon} + f_{\bar{j}}\epsilon \wedge d\bar{z}_j + f_{j} dz_j \wedge \bar{\epsilon} + f_{i\bar{j}} dz_i \wedge d\bar{z}_j \right)
\end{equation*}

By the assumption on the asymptotic behaviour we mean there exists some coordinate chart in which the zero-th order asymptotic of the metric agrees with the model metric. In other words, there is a constant $C$, such that

\begin{equation*}
\frac{1}{C}\omega_{\beta} \leq \omega \leq C\omega_{\beta}
\end{equation*}

In this note because we assume certain singularities for the Kähler manifold $X$ we must design our Kähler Ricci flow such that our flow preserve singularities. Now fix a simple normal crossing divisor $D = \sum_{i}(1 - \beta_{i})D_{i}$, where $\beta_{i} \in (0, 1)$ and simple normal crossing divisor $D$ means that $D_{i}$'s are irreducible smooth divisors and for any $p \in \operatorname{Supp}(D)$ lying in the intersection of exactly $k$ divisors $D_{1}, D_{2}, \dots, D_{k}$, there exists a coordinate chart $(U_{p}, \{z_{i}\})$ containing $p$, such that $D_{i}|_{U_{p}} = \{z_{i} = 0\}$ for $i = 1, \dots, k$.

If $S_{i} \in H^{0}(X, \mathcal{O}_{X}(L_{D_{i}}))$ is the defining sections and $h_{i}$ is hermitian metrics on the line bundle induced by $D_{i}$, then Donaldson showed that for sufficiently small $\epsilon_{i} > 0$, $\omega_{i} = \omega_{0} + \epsilon_{i}\sqrt{-1}\partial\overline{\partial}|S_{i}|_{h_{i}}^{2\beta_{i}}$ gives a conic Kähler metric on $X \setminus \operatorname{Supp}(D_{i})$ with cone angle $2\pi\beta_{i}$ along divisor $D_{i}$ and also if we set $\omega = \sum_{i=1}^{N}\omega_{i}$ then, $\omega$ is a smooth Kähler metric on $X \setminus \operatorname{Supp}(D)$ and

\begin{equation*}
||S||^{2(1-\beta)} = \prod_{i=1}^{k} ||S_{i}||^{2(1-\beta_i)}
\end{equation*}

where $S \in H^{0}(X, \mathcal{O}(L_{D}))$. Moreover, $\omega$ is uniformly equivalent to the standard cone metric

\begin{equation*}
\omega_{p} = \sum_{i=1}^{k} \frac{\sqrt{-1}dz_{i} \wedge d\bar{z}_{i}}{|z_{i}|^{2(1-\beta_{i})}} + \sum_{i=k+1}^{N} \sqrt{-1}dz_{i} \wedge d\bar{z}_{i}
\end{equation*}

From Tian-Fujiki theory, $|dz_{i}|^{2} = r_{i}^{2}$ for $1 \leq i \leq k$ and $|dz_{j}|^{2} = 1$ for $k + 1 \leq j \leq n$ so that A2) and A3) are again satisfied.

From now on for simplicity we write just "divisor $D$" instead "simple normal crossing divisor $D$".

We give an example of varieties which have conical singularities. Consider a smooth geometric orbifold given by $\mathbb{Q}$-divisor

\begin{equation*}
D = \sum_{j \in J} \left(1 - \frac{1}{m_{j}}\right)D_{j}
\end{equation*}

where $m_j \geq 2$ are positive integers and $\operatorname{Supp} D = \bigcap_{j \in J} D_j$ is of normal crossings divisor. Let $\omega$ be any Kähler metric on $X$, let $C > 0$ be a real number and $s_j \in H^0(X, \mathcal{O}_X(D_j))$ be a section defining $D_j$. Consider the following expression

\begin{equation*}
\omega_{D} = C \omega + \sqrt{-1} \sum_{j \in J} \partial \bar{\partial} | s_{j} |^{2 / m_{j}}
\end{equation*}

If $C$ is large enough, the above formula defines a closed positive (1, 1)-current (smooth away from $D$). Moreover

\begin{equation*}
\omega_{D} \geq \omega
\end{equation*}

in the sense of currents. Consider $\mathbb{C}^n$ with the orbifold divisor given by the equation

\begin{equation*}
\prod_{j = 1}^{n} z_{j}^{1 - 1 / m_{j}} = 0
\end{equation*}

(with eventually $m_j = 1$ for some $j$). The sections $s_j$ are simply the coordinates $z_j$ and a simple computation gives

\begin{equation*}
\omega_{D} = \omega_{eucl} + \sqrt{-1} \sum_{j = 1}^{n} \partial \bar{\partial} | z_{j} |^{2 / m_{j}} = \omega_{eucl} + \sqrt{-1} \sum_{j = 1}^{n} \frac{d z_{j} \wedge d \bar{z}_{j}}{m_{j}^{2} | z_{j} |^{2 (1 - 1 / m_{j})}}
\end{equation*}

Here we mention also metrics with non-conic singularities. We say a metric $\omega$ is of Poincare type, if it is quasi-isometric to

\begin{equation*}
\omega_{\beta} = \sqrt{-1} \left(\sum_{i = 1}^{k} \frac{d z_{i} \wedge d \bar{z}_{i}}{| z_{i} |^{2} \log^{2} | z_{i} |^{2}} + \sum_{i = k + 1}^{n} d z_{i} \wedge d \bar{z}_{i}\right)
\end{equation*}

It is always possible to construct a Poincare metric on $M \setminus D$ by patching together local forms with $C^\infty$ partitions of unity. Now, from Tian-Fujiki theory $|dz_i|^2 = r_i^2 (\log 1 / r_i)^2$, $1 \leq i \leq k$ and $|dz_j|^2 = 1$, $k + 1 \leq j \leq n$ so that A2) and A3) above are satisfied; we have

\begin{equation*}
v = \prod_{1 \leq i \leq k} r_{i}^{-2} (\log 1 / r_{i})^{-2}
\end{equation*}

Let $\Omega_P$ be the volume form on $X \setminus D$, then, there exists a locally bounded positive continuous function $c(z)$ on polydisk $\mathbb{D}^n$ such that

\begin{equation*}
\Omega_{P} = c(z) \sqrt{-1} \left(\bigwedge_{i = 1}^{k} \frac{d z_{i} \wedge d \bar{z}_{i}}{| z_{i} |^{2} \log^{2} | z_{i} |^{2}} + \bigwedge_{i = k + 1}^{n} d z_{i} \wedge d \bar{z}_{i}\right)
\end{equation*}

holds on $\mathbb{D}^n \cap (X \setminus D)$

\textbf{Remark:} Note that if $\Omega_P$ be a volume form of Poincare growth on $(X, D)$, with $X$ compact. If $c(z)$ be $C^2$ on $\mathbb{D}^n$, then $-Ric(\Omega_P)$ is of Poincare growth.

We say that $\omega$ is the homogeneous Poincare metric if its fundamental form $\omega_{\beta}$ is described locally in normal coordinates by the quasi-isometry

\begin{equation*}
\omega_{\beta} = \sqrt{-1} \left(\frac{1}{(\log | z_{1} z_{2} \dots z_{k} |^{2})^{2}} \sum_{i = 1}^{k} \frac{d z_{i} \wedge d \bar{z}_{i}}{| z_{i} |^{2}} + \sum_{i = 1}^{n} d z_{i} \wedge d \bar{z}_{i}\right)
\end{equation*}

and we say $\omega$ has Ball Quotient singularities if it is quasi-isometric to

\begin{equation*}
\omega_{\beta} = \sqrt{-1} \frac{d z_{1} \wedge d \bar{z}_{1}}{(| z_{1} | \log (1 / | z_{1} |))^{2}} + \sqrt{-1} \sum_{j = 2}^{n} \frac{d z_{j} \wedge d \bar{z}_{j}}{\log 1 / | z_{1} |}
\end{equation*}

It is called also Saper's distinguished metrics.

\begin{equation*}
\left| d z_{1} \right|^{2} = r_{1}^{2} (\log 1 / r_{1})^{2}, \quad \left| d z_{j} \right|^{2} = \log 1 / r_{1}, \quad k + 1 \leq j \leq n
\end{equation*}

so that A2) and A3) are satisfied; also we have the volume form as

\begin{equation*}
v = r_{1}^{-2} \left(\log 1 / r_{1}\right)^{-(n + 1)}
\end{equation*}

If $\omega$ is the fundamental form of a metric on the compact manifold $X$, and $\omega_{sap}$ be the fundamental forms of Saper's distinguished metrics and $\omega_{P,hom}$ be the fundamental forms of homogeneous Poincare metric, on the noncompact manifold $M \setminus D$, then $\omega_{sap} + \omega$ and $\omega_{P,hom}$ are quasi-isometric.

\textbf{Definition 2.23} A Kähler metric with cone singularities along $D$ with cone angle $2\pi \beta$ is a smooth Kähler metric on $X \setminus D$ which satisfies the following conditions when we write $\omega_{sing} = \sum_{i,j} g_{i\bar{j}} \sqrt{-1} dz_i \wedge d\bar{z}_j$ in terms of the local holomorphic coordinates $(z_1; \dots; z_n)$ on a neighbourhood $U \subset X$ with $D \cap U = \{z_1 = 0\}$

\begin{enumerate}
\item $g_{1\bar{1}} = F|z_1|^{2\beta - 2}$ for some strictly positive smooth bounded function $F$ on $X \setminus D$
\item $g_{1\bar{j}} = g_{i\bar{1}} = O(|z_1|^{2\beta -1})$
\item $g_{i\bar{j}} = O(1)$ for $i,j\neq 1$
\end{enumerate}

Now we shortly explain Donaldson's linear theory which is useful later in the definition of logarithmic Vafa-Yau's semi ricci flat metrics.

\textbf{Definition 2.24} 1) A function $f$ is in $C^{\gamma, \beta}(X, D)$ if $f$ is $C^\gamma$ on $X \setminus D$, and locally near each point in $D$, $f$ is $C^\gamma$ in the coordinate $(\hat{\zeta} = \rho e^{i\theta} = z_1 |z_1|^{\beta - 1}, z_j)$.

2) A (1,0)-form $\alpha$ is in $C^{\gamma,\beta}(X,D)$ if $\alpha$ is $C^\gamma$ on $X\setminus D$ and locally near each point in $D$, we have $\alpha = f_1\epsilon +\sum_{j > 1}f_jdz_j$ with $f_{i}\in C^{\gamma ,\beta}$ for $1\leq i\leq n$, and $f_{1}\to 0$ as $z_{1}\rightarrow 0$ where $\epsilon = e^{\sqrt{-1}\beta \theta}(d\rho +\sqrt{-1}\beta \rho d\theta)$

3) A $(1,1)$-form $\omega$ is in $C^{\gamma,\beta}(X,D)$ if $\omega$ is $C^\gamma$ on $X\setminus D$ and near each point in $D$ we can write $\omega$ as

\begin{equation*}
\omega = \sqrt{-1} \left(f \epsilon \wedge \bar{\epsilon} + f_{\bar{j}} \epsilon \wedge d\bar{z}_{j} + f_{j} d z_{j} \wedge \bar{\epsilon} + f_{i\bar{j}} d z_{i} \wedge d \bar{z}_{j}\right)
\end{equation*}

such that $f, f_j, f_{\bar{j}}, f_{i\bar{j}} \in C^{\gamma, \beta}$, and $f_j, f_{\bar{j}} \to 0$ as $z_1 \to 0$

4) A function $f$ is in $C^{2,\gamma,\beta}(X,D)$ if $f, \partial f, \partial \bar{\partial} f$ are all in $C^{\gamma,\beta}$

Fix a smooth metric $\omega_0$ in $c_{1}(X)$, we define the space of admissible functions to be

\begin{equation*}
\hat{C}(X, D) = C^{2,\gamma}(X) \cup \bigcup_{0 < \beta < 1} \left(\bigcup_{0 < \gamma < \beta^{-1} - 1} C^{2,\gamma,\beta}(X, D)\right)
\end{equation*}

and the space of admissible Kähler potentials to be

\begin{equation*}
\hat{\mathcal{H}}(\omega_{0}) = \{\phi \in \hat{C}(X, D) \mid \omega_{\phi} = \omega_{0} + \sqrt{-1}\partial\bar{\partial}\phi > 0 \}
\end{equation*}

Note that

\begin{equation*}
\mathcal{H}(\omega_{0}) \subset \hat{\mathcal{H}}(\omega_{0}) \subset \mathcal{PSH}(\omega_{0}) \cap L^{\infty}(X)
\end{equation*}

Where $\mathcal{PSH}(\omega_0)\cap L^\infty (X)$ is the space of bounded $\omega_0$-plurisubharmonic functions and

\begin{equation*}
\mathcal{PSH}(\omega_{0}) = \left\{\phi \in L_{loc}^{1}(X) \mid \phi \text{ is u.s.c and } \omega_{0} + \sqrt{-1}\partial\bar{\partial}\phi > 0 \right\}
\end{equation*}

The Ricci curvature of the Kählerian form $\omega_{D}$ on the pair $(X,D)$ can be represented as:


\begin{equation*}
Ric(\omega_D) = 2\pi \sum_j (1-\beta_j)[D_j] + \theta + \sqrt{-1}\partial\bar{\partial}\psi
\end{equation*}

with $\psi \in C^0(X)$ and $\theta$ is closed smooth $(1,1)$-form.

We have also $dd^c$-lemma on $X = \overline{X} \setminus D$. Let $\Omega$ be a smooth closed $(1,1)$-form in the cohomology class $c_1(K_{\overline{X}}^{-1} \otimes L_D^{-1})$. Then for any $\epsilon > 0$ there exists an explicitly given complete Kähler metric $g_\epsilon$ on $M$ such that

\begin{equation*}
Ric(g_\epsilon) - \Omega = \frac{\sqrt{-1}}{2\pi}\partial\bar{\partial}f_\epsilon \text{ on } X
\end{equation*}

where $f_\epsilon$ is a smooth function on $X$ that decays to the order of $O(\|S\|^\epsilon)$. Moreover, the Riemann curvature tensor $R(g_\epsilon)$ of the metric $g_\epsilon$ decays to the order of $O\left((-n\log \|S\|^2)^{-\frac{1}{n}}\right)$

Now we explain the logarithmic Weil-Petersson metric on moduli space of log Calabi-Yau manifolds(if it exists. for special case of rational surfaces it has been proven that such moduli space exists). The logarithmic Weil-Petersson metric has pole singularities [7] and we can introduce it also by elements of logarithmic Kodaira-Spencer tensors which represent elements of $H^1\left(X, \Omega_X^1(\log(D))^\vee\right)$. More precisely, Let $X$ be a complex manifold, and $D \subset X$ a divisor and $\omega$ a holomorphic $p$-form on $X \setminus D$. If $\omega$ and $d\omega$ have a pole of order at most one along $D$, then $\omega$ is said to have a logarithmic pole along $D$. $\omega$ is also known as a logarithmic $p$-form. The logarithmic $p$-forms make up a subsheaf of the meromorphic $p$-forms on $X$ with a pole along $D$, denoted

\begin{equation*}
\Omega_X^p(\log D)
\end{equation*}

and for the simple normal crossing divisor $D = \{z_1 z_2 \dots z_k = 0\}$ we can write the stalk of $\Omega_X^1(\log D)$ at $p$ as follows

\begin{equation*}
\Omega_X^1(\log D)_p = \mathcal{O}_{X,p} \frac{dz_1}{z_1} \oplus \dots \oplus \mathcal{O}_{X,p} \frac{dz_k}{z_k} \oplus \mathcal{O}_{X,p} dz_{k+1} \oplus \dots \oplus \mathcal{O}_{X,p} dz_n
\end{equation*}

Since, fibers are log Calabi-Yau manifolds and by recent result of Jeffres-Mazzeo-Rubinstein [71], we have Ricci flat metric on each fiber $(X_y, D_y)$ and hence we can have log semi-Ricci flat metric and by the same method of previous theorem, the proof of Theorem 0.8 is straightforward.


\textbf{Theorem 2.25} \textit{Let $(M, \omega_0)$ be a compact Kähler manifold with $D \subset M$ a smooth divisor and suppose we have topological constraint condition $c_1(M) = (1 - \beta)[D]$ where $\beta \in (0,1]$ then there exists a conical Kähler Ricci flat metric with angle $2\pi\beta$ along $D$. This metric is unique in its Kähler class. This metric is polyhomogeneous; namely, the Kähler Ricci flat metric $\omega_0 + \sqrt{-1}\partial\bar{\partial}\varphi$ admits a complete asymptotic expansion with smooth coefficients as $r \to 0$ of the form}

\begin{equation*}
\varphi(r, \theta, Z) \sim \sum_{j,k \geq 0} \sum_{l=0}^{N_{j,k}} a_{j,k,l}(\theta, Z) r^{j + k/\beta} (\log r)^l
\end{equation*}

\textit{where $r = |z_1|^\beta / \beta$ and $\theta = \arg z_1$ and with each $a_{j,k,l} \in C^\infty$}

Now we can introduce Logarithmic Yau-Vafa semi Ricci flat metrics. The volume of fibers $(X_y, D_y)$ are homological constant independent of $y$, and we assume that it is equal to 1. Since fibers are log Calabi-Yau varieties, so $c_1(X_y, D_y) = 0$, hence there is a smooth function $F_y$ such that $Ric(\omega_y) = \sqrt{-1}\partial\bar{\partial}F_y$. The function $F_y$ vary smoothly in $y$. By Jeffres-Mazzeo-Rubinstein's theorem[71], there is a unique conical Ricci-flat Kähler metric $\omega_{SRF,y}$ on $X_y \setminus D_y$ cohomologous to $\omega_0$. So there is a smooth function $\rho_y$ on $X_y \setminus D_y$ such that $\omega_0|_{X_y \setminus D_y} + \sqrt{-1}\partial\bar{\partial}\rho_y = \omega_{SRF,y}$ is the unique Ricci-flat Kähler metric on $X_y \setminus D_y$. If we normalize $\rho_y$, then $\rho_y$ varies smoothly in $y$ and defines a smooth function $\rho^D$ on $X \setminus D$ and we let

\begin{equation*}
\omega_{SRF}^D = \omega_0 + \sqrt{-1}\partial\bar{\partial}\rho^D
\end{equation*}

which is called as Log Semi-Ricci Flat metric.

Let $f : X \setminus D \to S$, be a smooth family of quasi-projective Kähler manifolds. Let $x \in X \setminus D$, and $(\sigma, z_2, \dots, z_n, s^1, \dots, s^d)$, be a coordinate centered at $x$, where $(\sigma, z_2, \dots, z_n)$ is a local coordinate of a fixed fiber of $f$ and $(s^1, \dots, s^d)$ is a local coordinate of $S$, such that

\begin{equation*}
f(\sigma, z_2, \dots, z_n, s^1, \dots, s^d) = (s^1, \dots, s^d)
\end{equation*}

Now consider a smooth form $\omega$ on $X \setminus D$, whose restriction to any fiber of $f$, is positive definite. Then $\omega$ can be written as

\begin{align*}
\omega(\sigma, z, s) = & \sqrt{-1} \left( \omega_{i\bar{j}} ds^i \wedge d\bar{s}^j + \omega_{i\bar{\beta}} ds^i \wedge d\bar{z}^\beta + \omega_{\alpha\bar{j}} dz^\alpha \wedge d\bar{s}^j + \omega_{\alpha\bar{\beta}} dz^\alpha \wedge d\bar{z}^\beta + \omega_{\sigma\bar{j}} d\sigma \wedge d\bar{s}^j \right. \\
& \left. + \omega_{i\bar{\sigma}} ds^i \wedge d\bar{\sigma} + \omega_{\sigma\bar{\sigma}} d\sigma \wedge d\bar{\sigma} + \omega_{\sigma\bar{j}} d\sigma \wedge d\bar{z}^j + \omega_{i\bar{\sigma}} dz^\alpha \wedge d\bar{\sigma} \right)
\end{align*}


Since $\omega$ is positive definite on each fibre, hence
\begin{equation*}
\sum_{\alpha, \beta = 2} \omega_{\alpha \bar{\beta}} d z^{\alpha} \wedge d \bar{z}^{\beta} + \omega_{\sigma \bar{\sigma}} d \sigma \wedge d \bar{\sigma} + \sum_{j = 2} \omega_{\sigma \bar{j}} d \sigma \wedge d \bar{z}^{j} + \sum_{i = 2} \omega_{i \bar{\sigma}} d z^{i} \wedge d \bar{\sigma}
\end{equation*}
gives a Kähler metric on each fiber $X_{s} \setminus D_{s}$. So
\begin{equation*}
\det (\omega_{\lambda \bar{\eta}}^{-1} (\sigma, z, s)) = \det \left( \begin{array}{cccc} \omega_{\sigma \bar{\sigma}} & \omega_{\sigma \bar{2}} & \ldots & \omega_{\sigma \bar{n}} \\ \omega_{2 \bar{\sigma}} & \omega_{2 \bar{2}} & \ldots & \omega_{2 \bar{n}} \\ \vdots & \vdots & \ddots & \vdots \\ \omega_{n \bar{\sigma}} & \omega_{n \bar{2}} & \ldots & \omega_{n \bar{n}} \end{array} \right)^{-1}
\end{equation*}
gives a hermitian metric on the relative line bundle $K_{X'/S}$ and its Ricci curvature can be written as $\sqrt{-1}\partial \bar{\partial}\log \det \omega_{\lambda \bar{\eta}}(\sigma ,z,s)$

\subsection*{2.1 \quad Vafa-Yau's semi Ricci-flat metric}

The volume of fibers $\pi^{-1}(y) = X_y$ is a homological constant independent of $y$, and we assume that it is equal to 1. Since fibers are Calabi-Yau manifolds so $c_1(X_y) = 0$, hence there is a smooth function $F_y$ such that $Ric(\omega_y) = \sqrt{-1}\partial \bar{\partial} F_y$ and $\int_{X_y}(e^{F_y} - 1)\omega_y^{n - m} = 0$. The function $F_y$ vary smoothly in $y$. By Yau's theorem there is a unique Ricci-flat Kähler metric $\omega_{SRF,y}$ on $X_y$ cohomologous to $\omega_0$. So there is a smooth function $\rho_y$ on $\pi^{-1}(y) = X_y$ such that $\omega_0|_{X_y} + \sqrt{-1}\partial \bar{\partial}\rho_y = \omega_{SRF,y}$ is the unique Ricci-flat Kähler metric on $X_y$. If we normalize by $\int_{X_y}\rho_y\omega_0^n |_{X_y} = 0$ then $\rho_y$ varies smoothly in $y$ and defines a smooth function $\rho$ on $X$ and we let
\begin{equation*}
\omega_{SRF} = \omega_{0} + \sqrt{-1} \partial \bar{\partial} \rho
\end{equation*}
which is called as Semi-Ricci Flat metric. Such Semi-Flat Calabi-Yau metrics were first constructed by Greene-Shapere-Vafa-Yau on surfaces [24]. More precisely, a closed real $(1,1)$-form $\omega_{SRF}$ on open set $U \subset X \setminus S$, (where $S$ is proper analytic subvariety contains singular points of $X$) will be called semi-Ricci flat if its restriction to each fiber $X_{y} \cap U$ with $y \in f(U)$ be Ricci-flat. Notice that $\omega_{SRF}$ is semi-positive in fiber direction but in horizontal direction it is not known the semi-positivity of such relative form.

Define
\begin{equation*}
\operatorname{Scal}(\omega) = \frac{n Ric(\omega) \wedge \omega^{n-1}}{\omega^n}
\end{equation*}


and take a holomorpic submerssion $\pi : X \to \mathbb{D}$ over some disc $\mathbb{D}$. It is enough to show smoothness of semi Ricci flat metric over a disc $\mathbb{D}$. Define a map $F: \mathbb{D} \times C^{\infty}(\pi^{-1}(b)) \to C^{\infty}(\pi^{-1}(b))$ by
\begin{equation*}
F(b, \rho) = Scal(\omega_0|_b + \sqrt{-1}\partial_b \bar{\partial}_b \rho)
\end{equation*}
By using Lemma 2.9 in [34] we can extend $F$ to a smooth map $\mathbb{D} \times L_{k+4}^2(\pi^{-1}(b)) \to L_k^2(\pi^{-1}(b))$. From the definition $F(b, \rho_b)$ is zero.

Recall that, if $L$ denotes the linearization of $\rho \to Scal(\omega_{\rho})$, then
\begin{equation*}
L(\rho) = \mathcal{D}^* \mathcal{D}(\rho) + \nabla Scal \cdot \nabla \rho
\end{equation*}
where $\mathcal{D}$ is defined by
\begin{equation*}
\mathcal{D} = \bar{\partial} \circ \nabla : C^{\infty} \rightarrow \Omega^{0,1}(T)
\end{equation*}
and $\mathcal{D}^*$ is $L^2$ adjoint of $\mathcal{D}$. The linearisation of $F$ with respect to $\rho$ at $b$ is given by
\begin{equation*}
\mathcal{D}_b^* \mathcal{D}_b : C^{\infty}(\pi^{-1}(b)) \rightarrow C^{\infty}(\pi^{-1}(b))
\end{equation*}
Leading order term of $L$ is $\Delta^2$. So $L$ is elliptic and we can use the standard theory of elliptic partial differential equations. Since the central fiber is Ricci flat metric and hence such Laplace has zero Kernel, hence $\mathcal{D}_b^*\mathcal{D}_b$ is an isomorphism. By the implicit function theorem, the map $b \to \rho_b$ is a smooth map
\begin{equation*}
\mathbb{D} \rightarrow L_k^2(\pi^{-1}(b)) \ ; \ \forall k.
\end{equation*}
By Sobolev embedding, it is a smooth map $\mathbb{D} \to C^r(\pi^{-1}(b))$ for any $r$. Hence $\rho_b$ is smooth in $b$. So semi Ricci flat metric $\omega_{SRF}$ is smooth (by assuming central fiber is smooth CY).

\textbf{Remark $\mathcal{A}$}: Semi Ricci flat metric $\omega_{SRF}$ is smooth(by assuming central fiber is smooth CY ).

Moreover, by the same method and using perturbation(when fibers are singular), log Semi Ricci flat metric $\omega_{SRF}^D$ is continuous when central fiber has canonical singularities and is CY itself .

Song-Tian introduced a new volume by using semi-Ricci flat metric and since volume is semi-positive , hence we must show that semi-Ricci flat metric is semi-positive when $\kappa(X)=dimX-1=n-1$ or $\omega_{SRF}$ is semi-positive when $dimX-\kappa(X)$ is odd number.

\textbf{Remark $\mathcal{B}$}: Semi Ricci flat metric $\omega_{SRF}$ and Log semi Ricci flat metric $\omega_{SRF}^D$ is semi-positive in both smooth(as form) and singular (as current) when Kodaira dimension


and log Kodaira dimension is positive respectively when $dimX - \kappa(X)$ is odd number. In general $\omega_{SRF}^{dimX - \kappa(X)}$ is semi-positive which is enough to show that Song-Tian volume is well defined.

Proof. R. Berman[14], showed that for Iitaka fibration $F : X \to Y$ a canonical sequence of Bergman type measures
\begin{equation*}
\nu_k = \int_{X^{N_k - 1}} \mu^{(N_k)}
\end{equation*}
where
\begin{equation*}
\mu^{(N_k)} = \frac{1}{Z_k} |S^{(k)}(z_1, \dots, z_{N_k})|^{2/k} d z_1 \wedge d \bar{z}_1 \dots d z_{N_k} \wedge d \bar{z}_{N_k}
\end{equation*}
and $N_{Z_k}$ is the normalizing constant ensuring that $\mu^{(N_k)}$ is a probability measure, and
\begin{equation*}
S^{(k)}(x_1, \dots, x_{N_k}) := \det \left( s_i^{(k)}(x_j) \right)_{1 \leq i, j \leq N_k}
\end{equation*}
where $s_i^{(k)}$ is a basis in $H^0(X, kK_X)$

then $v_k$ converges weakly to Song-Tian measure
\begin{equation*}
\mu_X = F^*(\omega_Y)^{\kappa(X)} \wedge \omega_{SRF}^{n - \kappa(X)}
\end{equation*}
We know, if $(a_n)$ be a convergent sequence of positive real numbers then the limit is semi-positive. So $\mu_X$ is semi-positive.

Now $\pi : X \to X_{can} = \operatorname{Proj}(R(X, K_X))$ and $\pi' : (X, D) \to X_{can}^D = \operatorname{Proj} \oplus_{m \geq 0} H^0(X, mK_X + \lfloor mD \rfloor)$ are Iitaka fibrations and since fibers are Calabi-Yau varieties and log Calabi-Yau varieties respectively. Hence $\omega_{SRF}^{n - \kappa(X)}$ and $(\omega_{SRF}^D)^{n - \kappa(X,D)}$ are semi-positive.

\textbf{Remark $\mathcal{C}$}: A result of Demailly implies that every Kähler current can be approximated (in the weak topology) by Kähler currents with analytic singularities(we introduce it later) in the same cohomology class. So the semi-Ricci flat metric as current has non-algebraic singularities due to variation of Hodge structure. So such Kähler currents with analytic singularities for $\omega_{SRF}$ are as same as Berman's sequences $v_k$.

\textbf{Remark $\mathcal{D}$}: We know that, for $X$ be a smooth variety such that the logarithmic Kodaira dimension of $X$ zero. Then the quasi-Albanese map $\alpha : X \to A = H^0(X, \Omega_X^1)^* / H_1(M, \mathbb{Z})$ is dominant and has irreducible general fibers and hence general fibers are connected. Note that when log Kodaira dimension is zero, we have quasi-Albanese map such that general fibers have zero log Kodaira dimension. If we assume

fibers be semi-ample, then fibers are log Calabi-Yau varieties. Hence we have smooth and positive semi-Ricci flat metric in this case.

\noindent \textbf{Conjecture:} When the log-Kodaira dimension is negative, then the general fibers are log Fano varieties and if we assume they admit Kähler-Einstein with positive Ricci curvature then the semi-Kähler Einstein metric(a form on total space which its restriction to each fiber is Kähler-Einstein metric with positive Ricci curvature) is no longer positive.

Note that $K_{X / Y} + D$ is nef and pseudoeffective. By a theorem of Boucksom[37], for pseudoeffective class $\alpha \in H^{1,1}_{\bar{\partial} \partial}(X,\mathbb{R})$, $\alpha$ is nef if and only if the Lelong number vanishes $\nu(\alpha ,x) = 0$ for all $x\in X$.

\noindent \textbf{Remark:} The Lelong number of semi-Ricci flat metric $\omega_{SRF}$ is zero

The Logarithmic Weil-Petersson metric along log canonical model $\pi : (X, D) \to X_{can}$

\begin{equation*}
Ric_{X / Y}^{h_{X / Y}^{(\omega_{SRF}^{D})}} (\omega_{SRF}) = \sqrt{-1} \partial \bar{\partial} \log \left(\frac{(\omega_{SRF}^{D})^{n} \wedge \pi^{*} \omega_{can}^{m}}{\pi^{*} \omega_{can}^{m}}\right) = \omega_{WP}^{D}
\end{equation*}

By the same computation in [38], the logarithmic Weil-Petersson metric on the moduli space of log-Calabi-Yau varieties(if it exists) is Kähler metric. The logarithmic Weil-Petersson metric has pole singularities and we can introduce it also by elements of logarithmic Kodaira-Spencer tensors which represent elements of $H^1 \left( X, \Omega_X^1 (\log(D))^{\vee} \right)$. More precisely, Let $X$ be a complex manifold, and $D \subset X$ a divisor and $\omega$ a holomorphic $p$-form on $X \setminus D$. If $\omega$ and $d\omega$ have a pole of order at most one along $D$, then $\omega$ is said to have a logarithmic pole along $D$. $\omega$ is also known as a logarithmic $p$-form. The logarithmic $p$-forms make up a subsheaf of the meromorphic $p$-forms on $X$ with a pole along $D$, denoted

\begin{equation*}
\Omega_{X}^{p} (\log D)
\end{equation*}

and for the simple normal crossing divisor $D = \{z_1 z_2 \ldots z_k = 0\}$ we can write the stalk of $\Omega_X^1(\log D)$ at $p$ as follows

\begin{equation*}
\Omega_{X}^{1} (\log D)_{p} = \mathcal{O}_{X, p} \frac{d z_{1}}{z_{1}} \oplus \dots \oplus \mathcal{O}_{X, p} \frac{d z_{k}}{z_{k}} \oplus \mathcal{O}_{X, p} d z_{k + 1} \oplus \dots \oplus \mathcal{O}_{X, p} d z_{n}
\end{equation*}

Now take $f: X \setminus D \to S$ where $\omega$ on each fibre $X_s \setminus D_s$ is positive and $S$ is one dimensional. Now consider the form $\omega^{n+1}$. It is well known that

\begin{equation*}
\omega^{n + 1} = c (\omega) . \omega^{n} \wedge \sqrt{-1} d s \wedge d \bar{s}
\end{equation*}

\noindent where $c(\omega)$ is the geodesic curvature. In fact, if $c(\omega)$ be positive then $\omega$ is positive. For proving the semi positivity of log semi Ricci flat metric, it is sufficient to show that $c(\omega_{SRF}^D)$ is semi positive. By applying the same computation of Y-J.Choi [35], and Berman[36], for the log semi Ricci flat metric $\omega_{SRF}^D$, we have the following PDE on $X_y \setminus D_y$

\begin{equation*}
- \Delta_{\omega_{SRF}^{D}} c (\omega_{SRF}^{D}) (\vartheta) = | \bar{\partial} \vartheta_{\omega_{SRF}^{D}} |_{\omega_{SRF}^{D}}^{2} - \Theta_{\vartheta \bar{\vartheta}} (f_{*} (K_{X^{\prime} / S}))
\end{equation*}

\noindent where $\vartheta \in T_yS$ and $\Theta$ is the curvature of $f_{*}(K_{X^{\prime} / S})$.

\noindent where $\bar{\partial}\vartheta_{\omega_{SRF}^D} = -\frac{\partial\omega_{s\bar{\beta}}\omega^{\bar{\beta}\alpha}}{\partial z^{\bar{\beta}}} dz^{\bar{\beta}}\otimes \frac{\partial}{\partial z^{\alpha}}$ is a $T^{(0,1)}(X_y\setminus D_y)$-valued (0,1)-form. By integrating both sides of our previous PDE over $X_{y}\setminus D_{y}$ we have

\begin{equation*}
 | \bar{\partial} \vartheta_{\omega_{SRF}^{D}} |_{\omega_{SRF}^{D}}^{2}  = \Theta_{s \bar{s}}
\end{equation*}

\noindent hence $\inf c(\omega_{SRF}^D) \geq 0$ and so the log semi Ricci flat $\omega_{SRF}^D$ is semi-positive and hence is a Kähler metric.

Singular Weil-Petersson metric has pole singularities. Because we are in deal with singular varieties, so we must consider $L^2$-cohomology and intersection cohomology. Let $\mathcal{M}_g$ be the moduli space of curves of genus $g$ and $\overline{\mathcal{M}}_g$ be the Deligne-Mumford compactification of $\mathcal{M}_g$ then it is known that there is a isomorphism between $L^2$ Cohomology of $\mathcal{M}_g$ with respect to Weil-Petersson metric $\omega_{WP}$ and Intersection cohomology and ordinary cohomology

\begin{equation*}
H_{(2)}^{*} (\mathcal{M}_{g}, \omega_{WP}) \cong H^{*} (\overline{\mathcal{M}}_{g})
\end{equation*}

Note that if we take $X^{reg} = X \setminus D$ where D is a divisor with normal crossings. Endow $X^{reg}$ with a complete Kahler metric which has Poincaré singularities normal to each component of $D$; in local coordinates, if $D = (z_1, \dots, z_k)$, the metric is quasi-isometric to

\begin{equation*}
\sum_{i = 1}^{k} \frac{d z_{i} \wedge d \bar{z}_{i}}{| z_{i} |^{2} (\log | z_{i} |^{2})^{2}} + \sum_{i = k + 1}^{n} d z_{i} \wedge d \bar{z}_{i}
\end{equation*}

The isomorphism $H_{(2)}^{*}(X \setminus D) \cong IH^{*}(X)$ was proved by Zucker

The Weil-Petersson metric on $\mathcal{M}_{g}$ is quasi-isometric with the following model,

\begin{equation*}
\sum_{i = 1}^{k} \frac{d z_{i} \wedge d \bar{z}_{i}}{| z_{i} |^{2} (- \log | z_{i} |)^{3}} + \sum_{i = k + 1}^{n} d z_{i} \wedge d \bar{z}_{i}
\end{equation*}

\noindent But we don't know yet that the Weil-Petersson metric on moduli space of Calabi-Yau manifolds is quasi-isometric with wich model.

\noindent \textbf{Remark A:} Ken-Ichi Yoshikawa [39], found an asymptotic formula for Weil-Petersson metric as follows

\begin{equation*}
\omega_{WP} = \left\{\frac{\ell}{| s |^{2} (\log | s |)^{2}} + O \left(\frac{1}{| s |^{2} (\log | s |)^{3}}\right) \right\} \sqrt{-1} d s \wedge d \bar{s}
\end{equation*}

\noindent \textbf{Remark B:} If $\pi : X \to Y$ be a surjective holomorphic fibre space of projective varieties $X, Y$ (with $\dim Y = 1$) where the central fibre $X_0$ has canonical singularities, then the hermitian metric $h$ (can be seen as the inverse of Song-Tian volume) of $f^*\omega_Y \wedge \omega_{SF}^{m-1}$ has following expression

\begin{equation*}
h = O \left(\left(\log \frac{1}{| \sigma_{D} |^{n}}\right)\right)
\end{equation*}

\noindent where $\sigma_{D}$ is a local defining section of discriminant locus $D$ of $\pi$ and $n$ is positive integer. In fact this is trivial due to Ken-Ichi Yoshikawa's expansion. We need it for our estimates of Monge-Ampere equation to solve the twisted Kähler-Einstein metric, $Ric(\omega) = -\omega + \omega_{WP}$.

\noindent \textbf{Remark C:} Since we are in deal with moduli space of log Calabi-Yau fibres and it is hard conjecture about the existence of such spaces and for special case of rational surfaces it has been done. In the preprint "Moduli of surfaces with an anti-canonical cycle, Mark Gross, Paul Hacking, Sean Keel" there is some discussion of the moduli space of log Calabi-Yau varieties in case $dimX = 2$ in Section 6 of "arxiv.org/abs/1211.6367". Hence the existence of moduli space of log Calabi-Yau varieties is the hard conjecture up to now.

By the same method of Theorem 2 of [40]

\noindent \textbf{Remark D:} The logarithmic Weil-Petersson metric on moduli space of log-Calabi-Yau pairs has zero Lelong number, hence the double dual of direct image of relative line bundle is nef.

We briefly explain about the Weil-Petersson metric on moduli space of polarized Calabi-Yau manifolds. We study the moduli space of Calabi-Yau manifolds via the

\noindent Weil-Petersson metric. We outline the important properties of such metrics here. The Weil-Petersson metric is not complete metric in general but in the case of abelian varieties and $K3$ surfaces, the Weil-Petersson metric turns out to be equal to the Bergman metric of the Hermitian symmetric period domain, hence is in fact complete Kähler Einstein metric. Weil and Ahlfors showed that the Weil-Petersson metric is a Kähler metric and later Tian gave a different proof for it. Ahlfors proved that it has negative holomorphic sectional, scalar, and Ricci curvatures. The quasi-projectivity of coarse moduli spaces of polarized Calabi-Yau manifolds in the category of separated analytic spaces (which also can be constructed in the category of Moishezon spaces) has been proved by Viehweg. By using Bogomolov-Tian-Todorov theorem[25], these moduli spaces are smooth Kähler orbifolds equipped with the Weil-Petersson metrics. Let $X \to M$ be a family of polarized Calabi-Yau manifolds. Lu and Sun showed that the volume of the first Chern class with respect to the Weil-Petersson metric over the moduli space $M$ is a rational number. Gang Tian proved that the Weil-Petersson metric on moduli space of polarized Calabi-Yau manifolds is just pull back of Chern form of the tautological of $\mathbb{CP}^N$ restricted to period domain which is an open set of a quadric in $\mathbb{CP}^N$ and he showed that holomorphic sectional curvature is bounded away from zero. Let $X$ be a compact projective Calabi-Yau manifold and let $f: X \to Y$ be an algebraic fiber space with $Y$ an irreducible normal algebraic variety of lower dimension then Weil-Petersson metric measuring the change of complex structures of the fibers.

Now, consider a polarized Kähler manifolds $\mathcal{X} \rightarrow S$ with Kähler metrics $g(s)$ on $X_{s}$. We can define a possibly degenerate hermitian metric G on S as follows: Take Kodaira-Spencer map

\begin{equation*}
\rho : T_{S, s} \to H^{1} (X, T_{X}) \cong H_{\bar{\partial}}^{0, 1} (T_{X})
\end{equation*}

\noindent into harmonic forms with respect to $g(s)$; so for $v, w \in T_{s}(S)$, we may define

\begin{equation*}
G (v, w) := \int_{\mathcal{X}_{s}} \langle \rho (v), \rho (w) \rangle_{g (s)}
\end{equation*}

\noindent When $\mathcal{X} \to S$ is a polarized Kähler-Einstein family and $\rho$ is injective $G_{WP} := G$ is called the Weil-Petersson metric on S. Tian-Todorov, showed that if we take $\pi : \mathcal{X} \to S$, $\pi^{-1}(0) = X_0 = X$, $\pi^{-1}(s) = X_s$ be the family of X, then S is a non-singular complex analytic space such that

\begin{equation*}
dim_{\mathbb{C}} S = dim_{\mathbb{C}} H^{1} (X_{s}, T X_{s})
\end{equation*}

Note that in general, if $f : X \to S$ be a smooth projective family over a complex manifold S. Then for every positive integer m,

\begin{equation*}
P_{m} (X_{s}) = dimH^{0} (X_{s}, \mathcal{O}_{X_{s}} (m K_{X_{s}}))
\end{equation*}

\noindent is locally constant function on $S$.

It is worth to mention that the fibers $X_{s}$ are diffeomorphic to each other and if fibers $X_{s}$ be biholomorphic then $\pi$ is holomorphic fiber bundle and Weil-Petersson metric is zero in this case in other words the Kodaira-Spencer maps

\begin{equation*}
\rho : T_{S, s} \rightarrow H^{1} (X_{s}, T_{X_{s}}) \cong H_{\bar{\partial}}^{0, 1} (T_{X_{s}})
\end{equation*}

\noindent are zero. In special case, let $dimX_{s}=1$, then the fibers are elliptic curves and $\pi$ is holomorphic fiber bundle and hence the Weil-Petersson metric is zero. In general, the Weil-Petersson metric is semipositive definite on the moduli space of Calabi-Yau varieties. Note that Moduli space of varieties of general type has Weil-Petersson metric also.

\noindent \textbf{Remark:} Let $(E, \|\cdot\|)$ be the direct image bundle $f_{*}(K_{X^{\prime} / S})$, where $X^{\prime} = X \setminus D$, of relative canonical line bundle equipped with the $L^2$ metric $\|\cdot\|$. Then the fibre $E_y$ is $H^0(X_y \setminus D_y, K_{X_y \setminus D_y})$. Since the pair $(X_y, D_y)$ is Calabi-Yau pair, hence $H^0(X_y \setminus D_y, K_{X_y \setminus D_y})$ is a 1-dimensional vector space. This implies that $E$ is a line bundle.

\noindent \textbf{Definition 2.26} (Tian's formula for Weil-Petersson metric) Take holomorphic fiber space $\pi : X \to B$ and assume $\Psi_y$ be any local non-vanishing holomorphic section of Hermitian line bundle $\pi_*(K_{X/B}^l)$, then the Weil-Petersson (1,1)-form on a small ball $N_r(y) \subset B$ can be written as

\begin{equation*}
\omega_{WP} = - \sqrt{-1} \partial_{y} \bar{\partial}_{y} \log \left((\sqrt{-1})^{n^{2}} \int_{X_{y}} (\Psi_{y} \wedge \overline{\Psi_{y}})^{\frac{1}{l}}\right)
\end{equation*}

\noindent Note that $\omega_{WP}$ is globally defined on $B$

Now because we are in deal with Calabi-Yau pair $(X, D)$ which $K_X + D$ is numerically trivial so we must introduce Log Weil-Petersson metrics instead Weil-Petersson metric. Here we introduce such metrics on moduli space of paired Calabi-Yau fibers $(X_y, D_y)$. Let $i: D \hookrightarrow X$ and $f: X \to Y$ be holomorphic mappings of complex manifolds such that $i$ is a closed embedding and $f$ as well as $f \circ i$ are proper and smooth. Then a holomorphic family $(X_y, D_y)$ are the fibers $X_y = f^{-1}(y)$ and $D_y = (f \circ i)^{-1}(y)$. There exists the moduli space of $\mathcal{M}$ of such family because any $(X_y, D_y)$ with trivial canonical bundle is non-uniruled. Now $X \setminus D$ is quasi-projective so we must deal with quasi-coordinate system instead of coordinate system. Let $(X, D)$ be a Calabi-Yau pair and take $X' = X \setminus D$ equipped with quasi-coordinate system. We say that a tensor $A$ on $X'$ which are covariant of type $(p, q)$ is quasi-$C^{k,\lambda}$-tensor, if it is of class $C^{k,\lambda}$ with respect to quasi-coordinates. Now we construct the logarithmic version of Weil-Petersson metric on moduli space of paired Calabi-Yau fibers $f: (X, D) \to Y$.[70]


\subsection*{Formulating logarithmic Weil-Petersson metric}

Now we are ready to state our theorem. We must mention that The result of Tian was on Polarized Calabi-Yau fibers and in this theorem we consider non-polarized fibers. [68][69]

\textbf{Theorem 2.27} \textit{Let $\pi : X \to Y$ be a smooth family of compact K\"{a}hler manifolds whith Calabi-Yau fibers. Then Weil-Petersson metric can be written as}
\begin{equation*}
\omega_{WP} = -\sqrt{-1}\partial_{y}\bar{\partial}_{y}\log \int_{X_{y}} |\Omega_{y}|^{2}
\end{equation*}
\textit{where $\Omega_y$ is a holomorphic $(n,0)$-form on $\pi^{-1}(U)$, where $U$ is a neighborhood of $y$}

\textbf{Proof:} For proof, We need to recall the Yau-Vafa semi Ricci flat metrics. Since fibers are Calabi-Yau varieties, so $c_{1}(X_{y}) = 0$, hence there is a smooth function $F_{y}$ such that $Ric(\omega_y) = \sqrt{-1}\partial \bar{\partial} F_y$ . The function $F_{y}$ vary smoothly in $y$. By Yau's theorem, there is a unique Ricci-flat K\"{a}hler metric $\omega_{SRF,y}$ on $X_{y}$ cohomologous to $\omega_0$ where $\omega_0$ is a K\"{a}hler metric attached to $X$. So there is a smooth function $\rho_y$ on $X_{y}$ such that $\omega_0|_{X_y} + \sqrt{-1}\partial \bar{\partial}\rho_y = \omega_{SRF,y}$ is the unique Ricci-flat K\"{a}hler metric on $X_{y}$. If we normalize $\rho_y$, then $\rho_y$ varies smoothly in $y$ and defines a smooth function $\rho$ on $X$ and we let
\begin{equation*}
\omega_{SRF} = \omega_{0} + \sqrt{-1}\partial \bar{\partial}\rho
\end{equation*}
which is called as semi-Ricci flat metric. Robert Berman and Y.J.Choi independently showed that the semi-Ricci flat metric is semi positive along horizontal direction. Now for semi Ricci flat metric $\omega_{SRF}$, we have
\begin{equation*}
\omega_{SRF}^{n+1} = c(\omega_{SRF}) \cdot \omega_{SRF}^{n} dy \wedge d\bar{y}
\end{equation*}
Here $c(\omega_{SRF})$ is called a geodesic curvature of semi $\omega_{SRF}$. Now from Berman and Choi formula, for $V \in T_yY$, the following PDE holds on $X_y$
\begin{equation*}
-\Delta_{\omega_{SRF}}c(\omega_{SRF})(V) = |\bar{\partial}V_{\omega_{SRF}}|_{\omega_{SRF}}^{2} - \Theta_{V\bar{V}}(\pi_{*}(K_{X/Y}))
\end{equation*}
$\Theta_{V\bar{V}}$ is the Ricci curvature of direct image of relative line bundle( which is a line bundle, since fibers are Calabi Yau manifolds ). Now by integrating on both sides of this PDE, since
\begin{equation*}
\int_{X} \Delta_{\omega_{SRF}}c(\omega_{SRF})(V) = 0
\end{equation*}


{\sloppy and from the definition of Weil-Petersson metric and this PDE we get $\pi^{*}\omega_{WP} = Ric(\pi_{*}(K_{X/Y}))$ and Choi showed that for some holomorphic $(n, 0)$-form $\Omega_{y}$ on $\pi^{-1}(U)$, where $U$ is a neighborhood of $y$ we have\par}
\begin{equation*}
Ric_{\omega_{SRF}}(\pi_{*}(K_{X/Y})) = -\sqrt{-1}\partial\bar{\partial}\log\|\Omega\|_{y}^{2}
\end{equation*}
Hence
\begin{equation*}
\omega_{WP} = -\sqrt{-1}\partial_{y}\bar{\partial}_{y}\log\int_{X_{y}}|\Omega_{y}|^{2}
\end{equation*}
and we obtain the desired result.

By the same method we can introduce the logarithmic Weil-Petersson metric on $\pi ; (X, D) \to Y$ with assuming fibers to be log Calabi-Yau manifolds and snc divisor $D$ has conic singularities, then we have
\begin{equation*}
\omega_{WP}^{D} = -\sqrt{-1}\partial_{y}\bar{\partial}_{y}\log\int_{X_{y}\setminus D_{y}}\frac{\Omega_{y}\wedge\bar{\Omega}_{y}}{\|S_{y}\|^{2}}
\end{equation*}
where $S_{y}\in H^{0}(X_{s},L_{D_{s}})$.

\textbf{Remark:} The fact is that the solutions of relative K\"{a}hler-Einstein metric or Song-Tian metric $Ric(\omega) = -\omega + f^{*}\omega_{WP}$ may not be $C^{\infty}$. In fact we have $C^{\infty}$ of solutions if and only if the Song-Tian measure or Tian's K\"{a}hler potential be $C^{\infty}$. Now we explain that under some following algebraic condition we have $C^{\infty}$-solutions for
\begin{equation*}
Ric(\omega) = -\omega + f^{*}\omega_{WP}
\end{equation*}
along Iitaka fibration. We recall the following Kawamata's theorem. [49]

\textbf{Theorem 2.28} \textit{Let $f: X \to B$ be a surjective morphism of smooth projective varieties with connected fibers. Let $P = \sum_{j} P_{j}$, $Q = \sum_{l} Q_{l}$, be normal crossing divisors on $X$ and $B$, respectively, such that $f^{-1}(Q) \subset P$ and $f$ is smooth over $B \setminus Q$. Let $D = \sum_{j} d_{j} P_{j}$ be a $\mathbb{Q}$-divisor on $X$, where $d_{j}$ may be positive, zero or negative, which satisfies the following conditions A,B,C:}

\vspace{0.5em}
\noindent\textbf{A)} $D = D^{h} + D^{v}$ such that any irreducible component of $D^{h}$ is mapped surjectively onto $B$ by $f$, $f: Supp(D^{h}) \to B$ is relatively normal crossing over $B \setminus Q$, and $f(Supp(D^v)) \subset Q$. An irreducible component of $D^{h}$ (resp. $D^{v}$) is called horizontal (resp. vertical)

\vspace{0.5em}
\noindent\textbf{B)} $d_{j} < 1$ for all $j$


\noindent \textbf{C)} The natural homomorphism $\mathcal{O}_B\to f_*\mathcal{O}_X(\lceil -D\rceil)$ is surjective at the generic point of $B$.\\
\textbf{D)} $K_{X} + D\sim_{\mathbb{Q}}f^{*}(K_{B} + L)$ for some $\mathbb{Q}$-divisor $L$ on $B$

\vspace{12pt}
\noindent Let
\begin{align*}
f^*Q_l &= \sum_j w_{lj} P_j \\
\bar{d}_j &= \frac{d_j + w_{lj} - 1}{w_{lj}}, \text{ \textit{if} } f(P_j) = Q_l \\
\delta_l &= \max \{ \bar{d}_j ; f(P_j) = Q_l \}. \\
\Delta &= \sum_l \delta_l Q_l. \\
M &= L - \Delta.
\end{align*}

\noindent Then $M$ is nef.

\vspace{12pt}
\noindent The following theorem is straightforward from Kawamata's theorem

\vspace{12pt}
\noindent \textbf{Theorem 2.29} \textit{Let $d_j < 1$ for all $j$ be as above in Theorem 2.28, and fibers be log Calabi-Yau pairs, then}
\begin{equation*}
\int_{X_s \setminus D_s} (-1)^{n^2/2} \frac{\Omega_s \wedge \overline{\Omega_s}}{|S_s|^2}
\end{equation*}
\textit{is continuous on a nonempty Zariski open subset of $B$.}

\vspace{12pt}
\noindent Since the inverse of volume gives a singular hermitian line bundle, we have the following theorem from Theorem 2.28

\vspace{12pt}
\noindent \textbf{Theorem 2.30} \textit{Let $K_{X} + D \sim_{\mathbb{Q}} f^{*}(K_{B} + L)$ for some $\mathbb{Q}$-divisor $L$ on $B$ and}
\begin{align*}
f^*Q_l &= \sum_j w_{lj} P_j \\
\bar{d}_j &= \frac{d_j + w_{lj} - 1}{w_{lj}}, \text{ \textit{if} } f(P_j) = Q_l \\
\delta_l &= \max \{ \bar{d}_j ; f(P_j) = Q_l \}.
\end{align*}


\begin{align*}
\Delta &= \sum_l \delta_l Q_l. \\
M &= L - \Delta.
\end{align*}

\noindent Then
\begin{equation*}
\left(\int_{X_s \setminus D_s} (-1)^{n^2/2} \frac{\Omega_s \wedge \overline{\Omega_s}}{|S_s|^2}\right)^{-1}
\end{equation*}
is a continuous hermitian metric on the $\mathbb{Q}$-line bundle $K_B + \Delta$ when fibers are log Calabi-Yau pairs.

\vspace{12pt}
\noindent \textbf{Remark:} Note that Yoshikawa [50],and befor him Tian [69] showed that when the base of Calabi-Yau fibration $f: X \to B$ is a disc and central fibre $X_0$ is reduced and irreducible and pair $(X, X_0)$ has only canonical singularities then Tian's Kähler potential can be extended to a continuous Hermitian metric lying in the following class
\begin{equation*}
\mathcal{B}(B) = C^\infty(S) \oplus \bigoplus_{r \in \mathbb{Q} \cap (0, 1]} \bigoplus_{k=0}^n |s|^{2r} (\log |s|)^k C^\infty(B)
\end{equation*}

\vspace{12pt}
\noindent \textbf{Remark:} Note that hermitian metric of Yau-Vafa semi Ricci flat metric $\omega_{SRF}$ is in the class of $\mathcal{B}(B)$ as soon as the central fibre $X_0$ is reduced and irreducible and pair $(X, X_0)$ has only canonical singularities

\vspace{12pt}
\noindent Now in next theorem we will find the relation between logarithmic Weil-Petersson metric and fiberwise Ricci flat metric which can be considered as the logarithmic version of Song-Tian formula [10, 11].

\vspace{12pt}
\noindent \textbf{Theorem 2.31} \textit{Let $\pi : (X, D) \to Y$ be a holomorphic family of log Calabi-Yau pairs $(X_s, D_s)$ for the Kähler varieties $X, Y$. Then we have the following relation between logarithmic Weil-Petersson metric and fiberwise Ricci flat metric.}
\begin{equation*}
\sqrt{-1} \partial \bar{\partial} \log \left(\frac{f^* \omega_Y^m \wedge (\omega_{SRF}^D)^{n-m}}{|S|^2}\right) = -f^* \operatorname{Ric}(\omega_Y) + f^* \omega_{WP}^D
\end{equation*}
where $S \in H^0(X, \mathcal{O}(L_D))$.


\noindent \textbf{Proof:} Take $X' = X \setminus D$. Choose a local nonvanishing holomorphic section $\Psi_y$ of $\pi_*(K_{X'/Y}^{\otimes l})$ with $y \in U \subset X'$. We define a smooth positive function on $\pi(U)$ by
\begin{equation*}
u(y) = \frac{(\sqrt{-1})^{(n-m)^2} (\Psi_y \wedge \overline{\Psi_y})^{\frac{1}{l}}}{(\omega_{SRF}^D)^{n-m}|_{X_y \setminus D_y}}
\end{equation*}
But the Numerator and Denominator of $u$ are Ricci flat volume forms on $X_y \setminus D_y$, so $u$ is a constant function. Hence by integrating $u(y)(\omega_{SRF}^D)^{n-m}|_{X_y \setminus D_y}$ over $X_y \setminus D_y$ we see that
\begin{equation*}
u(y) = \frac{(\sqrt{-1})^{(n-m)^2} \int_{X_y \setminus D_y} \frac{(\Psi_y \wedge \overline{\Psi_y})^{\frac{1}{l}}}{|S_y|^2}}{\int_{X_y \setminus D_y} \frac{(\omega_{SRF}^D)^{n-m}}{|S_y|^2}}
\end{equation*}
where $S_y \in H^0(X', \mathcal{O}(L_{D_y}))$.

\vspace{12pt}
\noindent But $y \mapsto \int_{X_y \setminus D_y} \frac{(\omega_{SRF}^D)^{n-m}}{|S_y|^2}$ is constant over $Y$. Hence the Logarithmic Weil-Petersson can be written as
\begin{equation*}
-\sqrt{-1} \partial \bar{\partial} \log u = \omega_{WP}^D \quad (*)
\end{equation*}
Now, to finish the proof we can write $\Psi_y = F(\sigma, y, z)(d\sigma \wedge dz_2 \wedge \ldots \wedge dz_{n-m})$ where $F$ is holomorphic and non-zero. Hence by substituting $\Psi_y$ in $u$ and rewriting $\sqrt{-1}\partial \bar{\partial}\log \left(\frac{f^*\omega_Y^m\wedge(\omega_{SRF}^D)^{n-m}}{|S|^2}\right)$ and using (*) we get the desired result.

\vspace{12pt}
\noindent Now we show that Song-Tian-Tsuji measure is bounded if and only if central fiber of Iitaka fibration has canonical singularities at worst. Not that Song-Tian just showed that such measure is in $L^{1,\epsilon}$ and they couldn't prove the boundedness of such measure.

\vspace{12pt}
\noindent \textbf{Theorem 2.32} \textit{Song-Tian-Tsuji measure}
\begin{equation*}
\Omega_{X/Y} = \left(\frac{\omega_{SRF}^n \wedge \pi^* \omega_{can}^m}{\pi^* \omega_{can}^m}\right)
\end{equation*}
\textit{on Iitaka fobration is bounded if and only if the central fiber $X_0$ has log terminal singularities.}

\vspace{12pt}
\noindent \textbf{Proof:} R. Berman in [14], showed that for Iitaka fibration $F : X \to Y$ a canonical sequence of Bergman type measures

\setcounter{page}{60}

\begin{equation*}
\nu_{k} = \int_{X^{N_{k}-1}} \mu^{(N_{k})}
\end{equation*}

where

\begin{equation*}
\mu^{(N_{k})} = \frac{1}{Z_{k}} | S^{(k)} (z_{1},\dots, z_{N_{k}}) |^{2 / k} d z_{1} \wedge d \bar{z}_{1}\dots d z_{N_{k}} \wedge d \bar{z}_{N_{k}}
\end{equation*}

and \(N_{Z_k}\) is the normalizing constant ensuring that \(\mu^{(N_k)}\) is a probability measure, and

\begin{equation*}
S^{(k)} (x_{1}, \dots , x_{N_{k}}) := \det \left(s_{i}^{(k)} (x_{j})\right)_{1 \leq i, j \leq N_{k}}
\end{equation*}

where \(s_i^{(k)}\) is a basis in \(H^0 (X,kK_X)\)

then \( v_{k} \) converges weakly to Song-Tian measure

\begin{equation*}
\mu_{X} = F^{*} (\omega_{Y})^{\kappa (X)} \wedge \omega_{SRF}^{n - \kappa (X)}
\end{equation*}

But from Proposition 1.17 [54], ([53] also), we know that If \( K_{X_0} \) is a \( \mathbb{Q} \)-Cartier divisor then \( X_0 \) has log terminal singularities if and only if \( \mathcal{L}_{X_0}^{2,m} = \mathcal{O}_{X_0}(m(K_{X_0})) \) for all \( m \geq 1 \) where \( L^{2,m} \) is the sheaf of locally \( 2/m \)-integrable \( m \)-fold \( d \)-norms. So Song-Tian measure is bounded if and only if central fiber has log terminal singularities (in log Iitaka fibration case from Proposition 1.19 [54], we need to assume log canonical singularities for central fiber to obtain boundedness of Song-Tian measure).

Mumford in his celebrated paper [55] introduced the notion of good metric (which we call such metrics as Mumford metric) to extend the Chern-Weil theory to quasi-projective manifolds. First we recall Mumford metrics. For fixing our notation. Let \(\tilde{E}\) be a holomorphic vector bundle of rank \(l\) over \(X\) and \(E = \tilde{E}|_X\) and \(h\) an Hermitian metric on \(E\) which may be singular near divisor \(D\). Now cover a neighborhood of \(D \subset X\) by finitely many polydiscs \(\{U_{\alpha} = (\Delta^n, (z_1, \dots, z_n))\}\) such that \(V_{\alpha} = U_{\alpha} \setminus D = (\Delta^*)^k \times \Delta^{n-k}\). Namely, \(U_{\alpha} \cap D = \{z_1\dots z_k = 0\}\). Let \(U = \cup_{\alpha} U_{\alpha}\) and \(V = \cup V_{\alpha}\) then on each \(V_{\alpha}\) we have the local Poincaré metric.

\textbf{Definition 2.33} Let \(\omega\) be a smooth local \(p\)-form defined on \(V_{\alpha}\) we say \(\omega\) has Poincaré growth if there is a constant \(C_{\alpha} > 0\) depending on \(\omega\) such that

\begin{equation*}
\left| \omega (t_{1}, t_{2},\dots, t_{p}) \right|^{2} \leq C_{\alpha} \prod_{i = 1}^{p} \left| \left| t_{i} \right| \right|^{2}
\end{equation*}

for any point \( z \in V_{\alpha} \) and \( t_1, \dots, t_p \in T_zX \) where \( ||.|| \) is taken on Poincaré metric. We say \( \omega \) is Mumford metric if both \( \omega \) and \( d\omega \) have Poincaré growth.


Now we recall the notion Hermitian Mumford metric \(h\) on the vector bundle \(E\).

\textbf{Definition 2.34} An Hermitian metric \(h\) on \(E\) is good if for all \( z \in V \) , assuming \( z \in V_{\alpha} \) , and for all basis \( (e_{1}, \dots, e_{n}) \) of \(E\) over \( U_{\alpha} \) , if we let \( h_{i\bar{j}} = h(e_{i}, e_{\bar{j}}) \) , then

I) \( |h_{i\bar{j}}|, (\det h)^{-1} \leq C(\sum_{i=1}^{k} \log |z_i|)^{2n} \) for some \( C > 0 \)

II) The local 1-forms \((\partial h.h^{-1})_{\alpha \lambda}\) are Mumford on \(V_{\alpha}\), namely the local connection and curvature forms of \(h\) have Poincaré growth.

If \( h \) is a good metric on \( E \), the Chern forms \( c_{i}(E, h) \) are Mumford forms.

\textbf{Theorem 2.35} \textit{Fiberwise Calabi-Yau metric on Iitaka fibration and hermitian metric corresponding to Song-Tian-Tsuji measure (inverse of Song-Tian-Tsuji measure \( h = \Omega_{X / Y}^{-1} \) gives a singular hermitian metric) is a good metric in the sense of Mumford, i.e., it is Mumford metric when central fiber has log terminal singularities at worst.}

\textbf{Proof:} From the Berman's formula in our proof of Theorem 2.33, and boundedness of Song-Tian-Tsuji measure and using Theorem 30, Lemma 36 in [56], we get the desired result.

\textbf{Remark:} Fiberwise log Calabi-Yau metric on log Iitaka fibration and hermitian metric corresponding to Song-Tian-Tsuji measure (inverse of Song-Tian-Tsuji measure \( h = \Omega_{(X,D)/Y}^{-1} \) gives a singular hermitian metric) is a good metric in the sense of Mumford, i.e., it is Mumford metric when central fiber has log canonical singularities at worst.

The following Theorem, gives a result to the goodness of canonical metric in the sense of Mumford, on moduli spaces of Calabi-Yau varieties with mild singularities in the sense of Minimal Model Prgram.

\textbf{Theorem 2.36} \textit{Since if a singular hermitian metric be a good metric in the sense of Mumford (Mumford metric), its first Chern current is also good metric in the sense of Mumford, hence from Theorem 2.32, the singular Weil-Petersson metric on moduli space of Calabi-Yau varieties is a good metric when central fiber has log terminal singularities at worst. Moreover if central fiber has log canonical singularities at worst, then the singular logarithmic Weil-Petersson metric on the moduli spaces of log Calabi-Yau varieties is a good metric in the sense of Mumford}


\textbf{Conjecture: : } If a singular Hermitian metric on a Kähler manifold be a good metric in the sense of Mumford then it has vanishing Lelong number.

\subsection*{Invariance of Plurigenera via Song-Tian program}

Now we talk about semi-positivity of logarithmic-Weil-Petersson metric via Invariance of plurigenera.

Let \(\pi : (X, D) \to Y\) be a smooth holomorphic fibre space whose fibres have pseudoeffective canonical bundles. Suppose that

\begin{equation*}
\frac{\partial \omega(t)}{\partial t} = -Ric_{X'/Y}(\omega(t)) - \omega(t) + [N]
\end{equation*}

be a Kähler ricci flow that starts with semi-positive Kähler form \([\omega(t)] = e^{-t}\omega_0 + (1 - e^{-t})\omega_{WP}^D\) and \(X' = X \setminus D\). Note that the logarithmic Weil-Petersson metric \(\omega_{WP}^D\) is semi-positive.

Then the invariance of plurigenera holds true if and only if the solutions \(\omega(t) = e^{-t}\omega_0 + (1 - e^{-t})\omega_{WP}^D\) are semi-positive(see the Analytical approach of Tsuji, Siu, Song-Tian). In fact an answer to this question leads to invariance of plurigenera in Kähler setting which still is open. Thanks to Song-Tian program. If our family of fibers be fiberwise KE-stable, then invariance of plurigenera holds true from \(L^2\)-extension theorem and also due to this fact that if the central fiber be pseudo-effective, then all the general fibers are pseudo-effective.

\textbf{Theorem 2.37} (\textit{\(L^2\)-extension theorem}) \textit{Let \(X\) be a Stein manifold of dimension \(n\), \(\psi\) a plurisubharmonic function on \(X\) and \(s\) a holomorphic function on \(X\) such that \(ds \neq 0\) on every branch of \(s^{-1}(0)\). We put \(Y = s^{-1}(0)\) and \(Y_0 = \{X \in Y; ds(x) \neq 0\}\) Let \(g\) be a holomorphic \((n-1)\)-form on \(Y_0\) with}

\begin{equation*}
c_{n-1} \int_{Y_0} e^{-\psi} g \wedge \bar{g} < \infty
\end{equation*}

\textit{where \(c_k = (-1)^{k(k-1)/2} (\sqrt{-1})^k\) Then there exists a holomorphic \(n\)-form \(G\) on \(X\) such that \(G(x) = g(x) \wedge ds(x)\) on \(Y_0\) and}

\begin{equation*}
c_n \int_X e^{-\psi} (1 + |s|^2)^{-2} G \wedge \bar{G} < 1620\pi c_{n-1} \int_{Y_0} e^{-\psi} g \wedge \bar{g}
\end{equation*}

\textbf{Theorem 2.38} \textit{(Siu) Assume \(\pi : X \to B\) is smooth, and every \(X_t\) is of general type. Then the plurigenera \(P_m(X_t) = \dim H^0(X_t, mK_{X_t})\) is independent of \(t \in B\) for any \(m\).}


After Siu, an ``algebraic proof'' is given, and applied to the deformation theory of certain type of singularities which appear in MMP by Kawamata.

\textbf{Definition 2.39} Let \( B \) be a normal variety such that \( K_B \) is \( \mathbb{Q} \)-Cartier, and \( f: X \to B \) a resolution of singularities. Then,

\begin{equation*}
K_{X} = f^{*} (K_{B}) + \sum_{i} a_{i} E_{i}
\end{equation*}

where \(a_{i} \in \mathbb{Q}\) and the \(E_{i}\) are the irreducible exceptional divisors. Then the singularities of \(B\) are terminal, canonical, log terminal or log canonical if \(a_{i} > 0, \geq 0, > -1\) or \(\geq -1\), respectively.

\textbf{Theorem 2.40} \textit{(Kawamata) If \( X_0 \) has at most canonical singularities, then \( X_t \) has canonical singularities at most for all \( t \in B \). Moreover, if all \( X_t \) are of general type and have canonical singularities at most, then \( P_m(X_t) = \dim H^0(X_t, mK_{X_t}) \) is independent of \( t \in B \) for all \( m \)}

\textbf{Remark:} If along holomorphic fiber space \((X, D) \to B\) (with some stability condition on \(B\)) the fibers are of general type then to get

\begin{equation*}
Ric(\omega) = \lambda \omega + \omega_{WP} + \text{additional term which come from higher canonical bundle formula}
\end{equation*}

, (here Weil-Petersson metric is a metric on moduli space of fibers of general type) when fibers are singular and of general type then from Theorem 1.35 we must impose that the central fiber \((X_0, D_0)\) must have canonical singularities and be of general type to obtain such result.

\textbf{Theorem 2.41} \textit{(Nakayama) If \( X_0 \) has at most terminal singularities, then \( X_t \) has terminal singularities at most for all \( t \in B \). Moreover, If \(\pi : X \to B\) is smooth and the ``abundance conjecture'' holds true for general \(X_t\), then \(P_m(X_t) = \dim H^0(X_t, mK_{X_t})\) is independent of \(t \in B\) for all \(m\).}

Takayama, showed the following important theorem

\textbf{Theorem 2.42} \textit{Let all fibers \( X_{t} = \pi^{-1}(t) \) have canonical singularities at most, then \( P_{m}(X_{t}) = \dim H^{0}(X_{t}, mK_{X_{t}}) \) is independent of \( t \in B \) for all \( t \)}

\textbf{Theorem 2.43} \textit{Let \(\pi : X \to Y\) be a proper smooth holomorphic fiber space of projective varieties such that all fibers \(X_{y}\) are of general type, then \(\omega_{WP}\) is semi-positive}


\textbf{Proof.} Let \(\pi : X \to Y\) be a smooth holomorphic fibre space whose fibres are of general type. Suppose that

\begin{equation*}
\frac{\partial \omega(t)}{\partial t} = -Ric_{X/Y}(\omega(t)) - \omega(t)
\end{equation*}

be a Kähler ricci flow that starts with semi-positive Kähler form \(\omega_0\) (take it Weil-Petersson metric).

Then since Siu's therems holds true for invariance of plurigenera, so the pseudo-effectiveness of \(K_{X_0}\) gives the pseudo-effectiveness of \(K_{X_t}\). The solutions of \(\omega(t)\) are semi-positive. But by cohomological characterization we know that \([\omega(t)] = e^{-t}\omega_{WP} + (1 - e^{-t})[\omega_0]\) and since \(\omega_0\) and \(\omega(t)\) are semi-positive, hence \(\omega_{WP}\) is semi-positive. \(\square\)

We consider the semi-positivity of singular Weil-Petersson metric \(\omega_{WP}\) in the sense of current.

\textbf{Theorem 2.44} \textit{Let \(\pi : X \to Y\) be a proper holomorphic fiber space such that all fibers \(X_{y}\) are of general type and have at worse canonical singularities, then the Weil-Petersson metric \(\omega_{WP}\) is semi-positive}

\textbf{Proof.} Suppose that

\begin{equation*}
\frac{\partial \omega(t)}{\partial t} = -Ric_{X/Y}(\omega(t)) - \omega(t)
\end{equation*}

be a Kähler Ricci flow. Then since Kawamata's therems say's that ``If all fibers \(X_t\) are of general type and have canonical singularities at most, then \(P_m(X_t) = \dim H^0(X_t, mK_{X_t})\) is independent of \(t \in B\) for all \(m\)'' hence invariance of plurigenera hold's true, and the solutions of \(\omega(t)\) are semi-positive by invariance of plurigenera. But by cohomological characterization we know that \([\omega(t)] = e^{-t}\omega_{WP} + (1 - e^{-t})[\omega_0]\) and since \(\omega_0\) and \(\omega(t)\) are semi-positive, hence \(\omega_{WP}\) is semi-positive. \(\square\)

From Nakayama's theorem, if \(X_0\) has at most terminal singularities, then \(X_t\) has terminal singularities at most for all \(t \in B\). Moreover, If \(\pi : X \to B\) is smooth and the ``abundance conjecture'' holds true for general \(X_t\), then \(P_m(X_t) = \dim H^0(X_t, mK_{X_t})\) is independent of \(t \in B\) for all \(m\). So when fibers are of general type then the solutions of Kähler Ricci flow \(\omega(t)\) is semi-positive and hence by the same method of the proof of previous Theorem, the Weil-Petersson metric \(\omega_{WP}\) is semi-positive.


\textbf{Remark II} Now assume that fibers are Fano Kähler-Einstein metrics and such family is fiberwise KE-stable, then
\begin{equation*}
\frac{\partial \omega(t)}{\partial t} = -Ric_{X/Y}(\omega(t)) + \omega(t)
\end{equation*}
be a Kähler Ricci flow that starts with semi-positive Kähler form $\omega_0$ (take it Weil-Petersson metric). Then semi-Kähler Einstein metric in fiber direction is semi-positive and in horizontal direction is not positive conjecturally. The fact is that the solutions of Kähler Ricci flow is no longer semi-positive.

So we give the following conjecture about Tian's K-stability via positivity theory.

Let $\pi : X \to \mathbb{D}$ be a proper holomorphic fibre space which the Fano fibers have unique Kähler-Einstein metric with positive Ricci curvature. Let the family of fibers is fiberwise KE-stable, then, there exists a unique Semi-Kähler Einstein metric $\omega_{SKE}$ on total space $X$ (as relative Kähler metric) such that its restriction to each fibre $X_s$ is Fano Kähler Einstein metric.

\textbf{Conjecture:} The Fano variety $X$ is K-poly stable if and only if for the proper holomorphic fibre space $\pi : X \to \mathbb{D}$ which the Fano fibers have unique Kähler-Einstein metric with positive Ricci curvature, then the fiberwise Kähler Einstein metric (Semi-Kähler Einstein metric) $\omega_{SKE}$ be smooth and semi-positive. Note that if fibers are K-poly stable then by Schumacher's result we have
\begin{equation*}
-\Delta_{\omega_i} c(\omega_{SKE}) - c(\omega_{SKE}) = |A|^2_{\omega_i}
\end{equation*}
{\sloppy where $A$ represents the Kodaira-Spencer class of the deformation and since $\omega_{SKE}^{n+1} = c(\omega_{SKE})\omega_{SKE}^n ds \wedge d\bar{s}$ so $c(\omega_{SKE})$ and $\omega_{SKE}$ have the same sign. By the minimum principle $\inf \omega_{SKE} < 0$. But our conjecture says that the fibrewise Fano Kähler-Einstein metric $\omega_{SKE}$ is smooth and semi-positive if and only if $X$ be K-poly stable.\par}

\section*{3 \quad Analytical MMP program via Gromov Invariant of Ruan-Tian}

In this section we study the long-time behavior of the conical Kähler-Ricci flow on pair $(X, D)$. We consider the classification of the singularity type of the conical Kähler Ricci flow on pair $(X, D)$. We classify the singularity type of the solutions of the conical Kähler Ricci flow into two classes called IIb, III which extends the results of Tosatti-Zhang [41]. We explain how the classification of the singularity type of the solutions


of the conical Kähler Ricci flow on pair $(X, D)$ is related to existence of a rational curve and hence we explain how Song-Tian program is related to the Gromov-Invariants of Ruan-Tian [75]

R. Hamilton [76] (see also [77]) classified the solutions of the Ricci flow as follows

We say:

The long-time solution of the \textbf{unnormalized conical Kähler-Ricci flow}
\begin{equation*}
\frac{\partial}{\partial t} \omega = -Ric(\omega) + [D] \quad , \omega(0) = \omega_X^D
\end{equation*}
is said to be of type IIb if
\begin{equation*}
\sup_{X \setminus D \times [0, \infty)} t|Rm(\omega(t))|_{\omega(t)} = +\infty
\end{equation*}
and of type III if
\begin{equation*}
\sup_{X \setminus D \times [0, \infty)} t|Rm(\omega(t))|_{\omega(t)} < +\infty
\end{equation*}

Moreover, we say:

The long-time solution of the \textbf{normalized conical Kähler-Ricci flow}
\begin{equation*}
\frac{\partial}{\partial t} \omega = -Ric(\omega) - \omega + [D] \quad , \omega(0) = \omega_X^D
\end{equation*}
is said to be of type IIb if
\begin{equation*}
\sup_{X \setminus D \times [0, \infty)} |Rm(\omega(t))|_{\omega(t)} = +\infty
\end{equation*}
and of type III if
\begin{equation*}
\sup_{X \setminus D \times [0, \infty)} |Rm(\omega(t))|_{\omega(t)} < +\infty
\end{equation*}

Note that if we replace $t$ to $T - t$, then we can classify short time solution of unnormalized conical Kähler Ricci flow and we can classify them into two classes I, IIa when $\sup_{X \setminus D \times [0, \infty)} (T - t) |Rm(\omega(t))|_{\omega(t)} < +\infty$ and $\sup_{X \setminus D \times [0, \infty)} (T - t) |Rm(\omega(t))|_{\omega(t)} = +\infty$.

\textbf{Definition 3.1} A log rational curve on a log pair $(X, D)$ is a rational curve $f: \mathbb{P}^1 \to X$ which meets $D$ at most once.


Existence of rational curves are very important in the classification of the solution of Kähler Ricci flows. Siu-Yau showed that when the bisectional curvature of a compact Kähler variety is positive then there exists a rational curve. Birkar, et al, showed that when $X$ is a Moishezon manifold which is not projective. Then $X$ contains a rational curve. Moreover we have the same result when $X$ is a normal Moishezon variety with analytic $\mathbb{Q}$-factorial singularities.

Now we explain how the existence of rational curve is related to Gromov-Witten Invariants.

Let $X$ be a smooth projective variety over $\mathbb{C}$. Fix a point $x \in X$, a homology class $A \in H_2(X(\mathbb{C}), \mathbb{Z})$ and very ample divisors in general position $H_i \subset X$, $i = 1, \dots, k$.

Let $y_0, \dots, y_k \in \mathbb{C}\mathbb{P}^1$ be general points. Take maps $f : \mathbb{C}\mathbb{P}^1 \to X$ such that $f_*[\mathbb{C}\mathbb{P}^1] = A$, $f(y_0) = x$ and $f(y_i) \in H_i$

Ruan-Tian, introduced an invariant
\begin{equation*}
\tilde{F}_{A,X}(x, H_1, \dots, H_k; y_0, \dots, y_k) := \text{number of such maps}
\end{equation*}

Gromov-Witten theory of pseudo-holomorphic curves shows that one can make a similar definition where $X$ is replaced by a symplectic manifold $(M, \omega)$. Then the corresponding invariant is denoted by
\begin{equation*}
\tilde{\Phi}_{A,X}(x, H_1, \dots, H_k; y_0, \dots, y_k) := \text{number of such maps}
\end{equation*}

Ruan-Tian, showed that if $\tilde{\Phi}_{A,X}(x, H_1, \dots, H_k; y_0, \dots, y_k) := \text{number of such maps}$ then there is a rational map $f : \mathbb{C}\mathbb{P}^1 \to X$ such that $f_*[\mathbb{C}\mathbb{P}^1] = A$, $f(y_0) = x$ and $f(y_i) \in H_i$

Moreover, $\tilde{\Phi}_{A,X}(x, H_1, \dots, H_k; y_0, \dots, y_k) = \tilde{F}_{A,X}(x, H_1, \dots, H_k; y_0, \dots, y_k)$ if the following conditions are satisfied.
\begin{itemize}
\item If $C_1, \dots, C_m \subset X$ are rational curves such that $\sum [C_i] = A$ and $x \in C_1$, then $m = 1$
\item If $g : \mathbb{C}\mathbb{P}^1 \to X$ is any map such that $g_*[\mathbb{C}\mathbb{P}^1] = A$, $g(y_0) = x$ then $H^1(\mathbb{C}\mathbb{P}^1, g^*TX) = 0$
\end{itemize}

Existence of rational curves are very important for the classification of the solution of Kähler Ricci flow. We think that existence of log rational curves is related to non-vanishing of logarithmic Gromov-Witten invariant of Gross, and Ruan-Tian. We will formulate it in future.


\textbf{Definition 3.2} Let $L \to X$ be a holomorphic line bundle over a projective manifold $X$. $L$ is said to be semi-ample if the linear system $|kL|$ is base point free for some $k \in \mathbb{Z}^+$. $L$ is said to be big if the Iitaka dimension of $L$ is equal to the dimension of $X$. $L$ is called numerically effective (nef) if $L.C \geq 0$ for any irreducible curve $C \subset X$

Now we prove the following theorem which is the generalization of the result of Tosatti-Zhang [41] on pair $(X, D)$.

We need to the following lemma due to McDuff-Salamon[?]

\textbf{Lemma 3.3} \textit{Let $r > 0$ and $a \geq 0$, If $w: B_r \to \mathbb{R}$ is a $C^2$ function that satisfies the inequalities}
\begin{equation*}
\Delta w \geq -aw^2, \quad w \geq 0, \quad \int_{B_r} w < \frac{\pi}{8a}
\end{equation*}
\textit{then}
\begin{equation*}
w(0) \leq \frac{8}{\pi r^2} \int_{B_r} w
\end{equation*}

\textbf{Theorem 3.4} \textit{Let $X$ be a compact Kähler manifold with simple normal crossing divisor $D$ with conic singularities. Let $K_X + (1 - \beta)D$ is nef and contains a possibly singular log rational curve $C \subset X$ such that $\int_C c_1(K_X + D) = 0$ which is equivalent to $\int_C c_1(X) = \int_C [D]$. Then any solution of the unnormalized conical Kähler Ricci flow on $(X, D)$ must be of type IIb.}

\textbf{Proof} Let $K_X + (1 - \beta)D$ is nef, Shen[74], showed that if $\omega(t)$ be the solution of the normalized conical Kähler Ricci flow
\begin{equation*}
\left\{ \begin{array}{l} \frac{\partial \omega(t)}{\partial t} = -Ric(\omega(t)) - \omega(t) + 2\pi(1 - \beta)[D] \\ \omega_{\beta}(., 0) = \omega_0 + \sqrt{-1}\delta \partial \bar{\partial} \|S\|^{2\beta} \end{array} \right.
\end{equation*}
then $\omega(t)$ exists for all the time. By the same method of Tosatti[41], by using Chern-Gauss-Bonnet theorem and Lemma 4.3, if $(X,D)$ contain a log rational curve and $\int_C c_1(K_X + D) = 0$, we have
\begin{equation*}
\sup_{X \setminus D} Bisec_{\omega(t)} \geq \frac{\pi}{8 \int_C \omega} > 0
\end{equation*}

Take $X' = X \setminus D$, if there are points $x_k \in X'$, times $t_k \to \infty$, 2-planes $\pi_k \subset T_{x_k}X'$ and a constant $\kappa > 0$ such that
\begin{equation*}
Sec_{\omega(t_k)}(\pi_k) \geq \kappa
\end{equation*}


for all $k$, then in particular $\sup_{X'} |Rm(\omega(t_k))|_{\omega(t_k)} \geq \kappa$ and so
\begin{equation*}
\sup_{X'} |Rm(\omega(t))|_{\omega(t)} = +\infty
\end{equation*}

Now, by using the definition of bisectional curvature, and since bisectional curvature is the sum of two sectional curvature, we get the result and proof is complete. $\square$

For classifying the solutions of the conical Kähler Ricci flow when logarithmic Kodaira dimension is zero, we need to the Kawamata's theorem on the quasi-Albanese maps for varieties of the logarithmic Kodaira dimension zero which is based on Deligne's theory of mixed Hodge structures.

\textbf{Theorem 3.5 (Kawamata[42])} \textit{Let $X$ be a smooth variety such that the logarithmic Kodaira dimension $\kappa(X)$ of $X$ is zero. Then the quasi-Albanese map $\alpha: X \to A = Alb(X)$ is dominant and has irreducible general fibers.}

In fact, by taking a suitable base $\{\omega_j\}_{j=1}^{p}$ of $H^0(X, \Omega_X^1(\log D))$ we have the following Iitaka's quasi-Albanese map
\begin{equation*}
\alpha : x \in X \setminus D \to \left(\int_{x_0}^x \omega_1, \dots, \int_{x_0}^x \omega_p\right) \in \mathbb{C}^p / \Gamma
\end{equation*}
where
\begin{equation*}
\Gamma = \left\{\left(\int_{\lambda} \omega_1, \dots, \int_{\lambda} \omega_p\right), \lambda \in H^1(X \setminus D, \mathbb{Z}) \right\}
\end{equation*}
The Abelian variety $Alb(X) = \mathbb{C}^p / \Gamma$ is called quasi-Albanese variety which is semi-Torus

We need to the following lemma to state the next theorem.

\textbf{Lemma 3.6} \textit{Let $K_X + D$ be semi-ample and logarithmic Kodaira dimension $\kappa(X, D) = 0$, then $(K_X + D)^m$ is trivial for some $m \geq 1$}

The non-logarithmic version of the following theorem is due to Zhang-Tosatti [?]

\textbf{Theorem 3.7} \textit{Let $X$ be a compact Kähler $n$-manifold with snc divisor $D$ with conic singularities such that $K_X + D$ semi-ample, and logarithmic Kodaira dimension $\kappa(X, D) = 0$ and consider a solution of conical Kähler Ricci flow}
\begin{equation*}
\left\{ \begin{array}{l} \frac{\partial \omega(t)}{\partial t} = -Ric(\omega(t)) - \omega(t) + 2\pi(1 - \beta)[D] \\ \omega_{\beta}(., 0) = \omega_0 + \sqrt{-1}\delta \partial \bar{\partial} \|S\|^{2\beta} \end{array} \right.
\end{equation*}

\begin{itemize}
\item \textit{A) If $X$ is not a finite quotient of a quasi-Albanese variety, then the solution is of type IIb}
\item \textit{B) If $X$ is a finite quotient of a quasi-Albanese variety, then the solution is of type III}
\end{itemize}

\noindent \textbf{Proof} For the proof we apply the Zhang-Tosatti's method [41]. From the Lemma 4.6, $K_{X} + D$ is torsion. From Jeffres-Mazzeo-Rubinstein's Theorem 1.7 [71], we have the Ricci flat metric $\omega_{\infty}$ in the class $[\omega]$ which satisfies in $Ric(\omega_{\infty}) = [D]$. Let $X$ is not a finite quotient of a quasi-Albanese variety, then $\omega_{\infty}$ is not flat. Let $x \in X' := X \setminus D$ be a point and $\pi \subset T_xX'$ a 2-plane with $Sec_{\omega_{\infty}}(\pi) \geq \kappa > 0$ for some constant $\kappa$. Since $\omega_{\infty}$ is not flat and $Ric(\omega_{\infty}) = 0$, there must be a positive sectional curvature at $x$. Now, since by the result of Shen[74], $\omega(t)$ converges to $\omega_{\infty}$, then also $Sec_{\omega(t)}(\pi) > 0$. We know, if there are points $x_k \in X'$, times $t_k \to \infty$, 2-planes $\pi_k \subset T_{x_k}X'$ and a constant $\kappa > 0$ such that
\begin{equation*}
Sec_{\omega(t_k)}(\pi_k) \geq \kappa
\end{equation*}
for all $k$, then in particular $\sup_{X'} |Rm(\omega(t_k))|_{\omega(t_k)} \geq \kappa$ and so
\begin{equation*}
\sup_{X'} |Rm(\omega(t))|_{\omega(t)} = +\infty
\end{equation*}
When $X$ is the finite quotient of the quasi-Albanese variety, then it is known that the solution of the conical Kähler Ricci flow converges to $\omega_{\infty}$ exponentially fast and we have $|Rm(\omega(t))|_{\omega(t)} \leq Ce^{-\eta t}$ and the flow is of type III. \hfill $\square$

\vspace{1em}
\noindent \textbf{Definition 3.8} Moishezon manifold $M$ is a compact complex manifold such that the field of meromorphic functions on each component $M$ has transcendence degree equal the complex dimension of the component:
\begin{equation*}
\dim_{\mathbf{C}} M = a(M) = \operatorname{tr.deg.}_{\mathbf{C}} \mathbf{C}(M)
\end{equation*}
Complex algebraic varieties have this property,

\vspace{1em}
\noindent \textbf{Lemma 3.9} \textit{Let $X$ be a normal compact complex space, such that $K_{X}$ is big line bundle, then $X$ is Moishezon. Moreover, if $X$ be Kähler then $X$ is projective}

\vspace{1em}
\noindent Now, we need to recall the following theorem of Y.Kawamata

\vspace{1em}
\noindent \textbf{Theorem 3.10} \textit{(Kawamata [42]) Let $f: X \to Y$ be a projective surjective morphism, $D$ an effective $\mathbb{Q}$-divisor on $X$, and $E$ an irreducible component of}
\begin{equation*}
\{x \in X; \dim_{x} f^{-1} f(x) > \dim X - \dim Y \}
\end{equation*}
\textit{Suppose that the pair $(X, D)$ has only log-terminal singularities and that $-(K_X + D)$ is $f$-nef. Then $E$ is covered by a family of rational curves}

\noindent \textbf{Theorem 3.11} \textit{Let $X$ be a compact Kähler $n$-manifold with snc $\mathbb{Q}$-effective divisor $D$ with conic singularities such that $K_{X} + D$ semi-ample, pair $(X, D)$ has only log-terminal singularities, and $(X, D)$ is of log general type, i.e logarithmic Kodaira dimension $\kappa(X, D) = \dim X$ and consider a solution of conical Kähler Ricci flow}
\begin{equation*}
\left\{ \begin{array}{l} \frac{\partial \omega(t)}{\partial t} = -Ric(\omega(t)) - \omega(t) + 2\pi(1-\beta)[D] \\ \omega_{\beta}(., 0) = \omega_{0} + \sqrt{-1}\delta\partial\bar{\partial}\|S\|^{2\beta} \end{array} \right.
\end{equation*}
\begin{itemize}
\item \textit{A) If $K_{X} + D$ is ample, then the soultion is of type III}
\item \textit{B) If $K_{X} + D$ is not ample, then the solution is of type IIb}
\end{itemize}

\noindent \textbf{Proof} For our proof, we apply the method of Zhang-Tosatti [41]. Let $K_{X} + D$ is ample, then by the result of Tian-Yau, Cheng-Yau, the conical Kähler Ricci flow converges to conical Kahler Einstein metric $\omega_{\infty}$, which satisfies in $Ric(\omega_{\infty}) = -\omega_{\infty} + [D]$. We can construct an initial metric which has bounded bisectional curvature and since sectional curvature remain bounded along conical Kähler Ricci flow, hence the flow has solution of type III.

Now, let $K_{X} + D$ is not an ample divisor. From Lemma 3.9. $X$ is Moishezon and Projective. Hence the linear system $|m(K_X + D)|$ gives a rational map $f: (X, D) \to \mathbb{CP}^N$. It is clear that $f$ is not an isomorphism with its image $f(X)$, since we assumed $K_{X} + D$ is not ample. Now, take a fibre $F = (X_s, D_s)$ of $f$. From the Theorem 4.10, $F$ is uniruled and since for a log rational curve $C \subset F = (X_s, D_s)$, $f(C) = \{a\}$, so $\int_C c_1(K_X + D) = 0$. By the Theorem 4.4, the solution of the conical Kähler Ricci flow must be of type IIb. \hfill $\square$

\vspace{1em}
\noindent \textbf{Remark:} It is worth to mention that when $0 < \kappa(X) < n$, we can extend the Tosatti's inequality as follows
\begin{equation*}
\sup_{X \setminus D} Bisec_{\omega(t)} \geq \frac{\pi}{8\left(\int_{C} \omega + \int_{C} \omega_{WP}\right)} > 0
\end{equation*}
where $\omega(t)$ is the solution of relative conical Kähler Ricci flow.

One of the main important problem in Kähler geometry is to prove the same following theorem for Kähler varieties. We give one of applications of the classification of the solutions of Kähler Ricci flow

\vspace{1em}
\noindent \textbf{Theorem:} \textit{Let $X$ be a complex projective manifold. Then $X$ is non-uniruled if and only if its canonical divisor $K_{X}$ is pseudo-effective.}

\noindent One of idea for the solution in Kähler setting is to use Kähler Ricci flow theory. In fact we can translate is as follows. If $X$ has Kodaira dimension zero, then $K_X$ is pseudo effective iff Kähler Ricci flow has long time solution and when $X$ is finite quotient of Torus, then there is no rational curve in $X$ and we get non-uniruleness. When $X$ is not finite quotient of Torus then we can not have uniruleness, since the solution is of type IIb. When we have intermidiate Kodaira dimension, we can translate it in the relative Kähler Ricci flow theory along Iitaka fibration. i.e. $K_{X/Y}$ is pseudo effective iff relative Kähler Ricci flow has long-time solution. Note that

The long-time solution of the \textbf{relative Kähler-Ricci flow} along Iitaka fibration when $0 < \kappa(X) < n$
\begin{equation*}
\frac{\partial}{\partial t}\omega = -Ric_{X/Y}(\omega) - \omega
\end{equation*}
is said to be of type IIb if
\begin{equation*}
\sup |Rm_{X/Y}(\omega(t))|_{\omega(t)} = +\infty
\end{equation*}
and of type III if
\begin{equation*}
\sup |Rm_{X/Y}(\omega(t))|_{\omega(t)} < +\infty
\end{equation*}
Note that such solution exists when family of fibers is fiberwise KE-stable and central fiber has canonical singularities at worst.

\subsection*{Stable family and Relative Kähler-Einstein metric}

For compactification of the moduli spaces of polarized varieties Alexeev, and Kollar-Shepherd-Barron, [63] started a program by using new notion of moduli space of "stable family". They needed to use the new class of singularities, called semi-log canonical singularities. See [60], [64][65].

Let $X$ be an equidimensional algebraic variety that satisfies Serre's $S_2$ condition and is normal crossing in codimension one. Let $\Delta$ be an effective $\mathbb{R}$-divisor whose support does not contain any irreducible components of the conductor of $X$. The pair $(X, \Delta)$ is called a semi log canonical pair (an slc pair, for short) if
\begin{enumerate}
\item[(1)] $K_X + \Delta$ is $\mathbb{R}$-Cartier;
\item[(2)] $(X^\nu, \Theta)$ is log canonical, where $\nu : X^\nu \to X$ is the normalization and $K_{X^\nu} + \Theta = \nu^*(K_X + \Delta)$
\end{enumerate}
Note that, the conductor $C_X$ of $X$ is the subscheme defined by, $\operatorname{cond}_X := \mathcal{H}om_{\mathcal{O}_X}(\nu_*\mathcal{O}_{X^\nu}, \mathcal{O}_X)$.

\noindent A morphism $f: X \to B$ is called a weakly stable family if it satisfies the following conditions:
\begin{enumerate}
\item $f$ is flat and projective
\item $\omega_{X / B}$ is a relatively ample $\mathbb{Q}$-line bundle
\item $X_{b}$ has semi log canonical singularities for all $b\in B$
\end{enumerate}
A weakly stable family $f: X \to B$ is called a stable family if it satisfies Kollar's condition, that is, for any $m \in \mathbb{N}$
\begin{equation*}
\omega_{X / B}^{[ m ]} \vert_{X_{b}} \cong \omega_{X_{b}}^{[ m ]}.
\end{equation*}
Note that, if the central fiber be Gorenstein and stable variety, then all general fibers are stable varieties, i.e., stability is an open condition

\vspace{1em}
\noindent \textbf{Conjecture:} Weil-Petersson metric (or logarithmic Weil-Petersson metric) on stable family is semi-positive as current and such family has finite distance from zero i.e $d_{WP}(B,0) < \infty$ when central fiber is stable variety also.

Moreover we predict the following conjecture holds true.

\vspace{1em}
\noindent \textbf{Conjecture:} Let $f: X \to B$ is a stable family of polarized Calabi-Yau varieties, and let $B$ is a smooth disc. then if the central fiber be stable variety as polarized Calabi-Yau variety, then we have following canonical metric on total space.
\begin{equation*}
Ric(\omega) = -\omega + f^* (\omega_{WP}) + [ N ]
\end{equation*}
Moreover, if we have such canonical metric then our family of fibers is stable.

We predict that if the base be singular with mild singularites of general type(for example $B = X_{can}$) then we have such canonical metric on the stable family and central fiber need to be stable after base change and relative canonical model

\vspace{1em}
\noindent \textbf{Conjecture:} The twisted Kähler-Einstein metric $Ric(\omega) = -\omega + \alpha$ where $\alpha$ is a semi-positive current has unique solution if and only if $\alpha$ has zero Lelong number

Now the following formula is cohomological characterization of Relative Kähler-Ricci flow due to Tian

\vspace{1em}
\noindent \textbf{Theorem 3.12} \textit{The maximal time existence $T$ for the solutions of relative Kähler Ricci flow is}
\begin{equation*}
T = \sup \left\{t \mid e^{-t} \left[ \omega_{0} \right] + (1 - e^{-t}) c_{1} \left(K_{X / Y} + D\right) \in \mathcal{K} ((X, D) / Y) \right\}
\end{equation*}
\textit{where $\mathcal{K}\left((X,D) / Y\right)$ denote the relative Kähler cone of $f:(X,D)\to Y$}

\section*{4 \quad Canonical metrics on Arithmetic varieties}

In this section we give an Arithmetic version of Song-Tian program on arithmetic varieties and call it \textbf{Arithmetic Song-Tian program}. We apply the Arithmetic Minimal Model program which was introduced by Yuji Odaka[59], for finding twisted Kähler Einstein metric(canonical metrics) on Arithmetic varieties which do not have definite Arithmetic first Chern class.

Philosophically, Song-Tian believe that

\begin{center}
PDE surgery $\Longleftrightarrow$ Geometric surgery $\Longleftrightarrow$ Algebraic surgery $\Longleftrightarrow$ Arithmetic surgery
\end{center}

For our purpos, we use the language of Arakelov geometry in Arithmetic geometry and Arakelov-Gillet-Soule's intersection theory. In final by applying Boucksom-Jonsson's works on Monge-Ampere equation, we give non-archimedean version of Song-Tian program.

Arakelov theory, also known as arithmetic intersection theory, is used to study number theoretic problems from a geometrical point of view. Arakelov defined an intersection theory on arithmetic surfaces over the ring of integers of a number field. He showed that geometry over number fields in addition with dierential geometry on some corresponding complex manifolds behaves like geometry over a compact variety

In 1987, P. Deligne generalized in [45] the arithmetic intersection theory of Arakelov. Deligne opened the way to a higher dimensional generalization. In 1991, H. Gillet and C. Soulé in [46] extended the arithmetic intersection theory to higher dimensions by translating the theory of Green's functions to the more manageable notion of Green's currents.

Broadly speaking, Arakelov geometry can be seen as

\begin{center}
\textbf{Arakelov geometry = Grothendieck algebraic geometry of scheme+Hermitian complex geometry}
\end{center}

Now, we review the fundamental notions of Arakelov geometry, as developed in Arakelov's paper [43] and Falting's paper[44].

\vspace{1em}
\noindent \textbf{Definition 4.1} An arithmetic variety $\pi : \mathcal{X} \to \operatorname{Spec}\mathcal{O}_K$ is a reduced, regular scheme $\mathcal{X}$, which is flat and projective over $\operatorname{Spec}\mathcal{O}_K$, where as usual $\mathcal{O}_K$ denotes the ring of integers of a number field $K$. Moreover, we assume that the generic fibre $X$ is geometrically connected. Let $d \in \mathbb{N}$ be the relative dimension of $\mathcal{X}$, hence $\dim(\mathcal{X}) = d + 1$. If $d = 1$, we call $\mathcal{X}$ as arithmetic surface.

\noindent\textbf{Definition 4.2} Let $\mathcal{X}_{\infty} = X(\mathbb{C})$ denote the set of complex points of the generic fibre $X$.

\begin{equation*}
X(\mathbb{C}) = \coprod_{\sigma : K \hookrightarrow \mathbb{C}} X_{\sigma}(\mathbb{C})
\end{equation*}

\noindent On $\mathcal{X}_{\infty}$ let $A^{(p,q)}(\mathcal{X}_{\infty})$ be the space of smooth $(p,q)$-forms endowed with the Schwartz topology. This means that a sequence $(\eta_n)$ in $A^{(p,q)}(\mathcal{X}_{\infty})$ converges to $\eta$ in $A^{(p,q)}(\mathcal{X}_{\infty})$ if and only if there exists a compact set $K$ such that for any $n$ we have $\mathrm{supp}(\eta_n) \subset K$ and any derivation of $\eta_n$ converges uniformly to the corresponding derivation of $\eta$. The space of $(p,q)$-currents $D^{(p,q)}(\mathcal{X}_{\infty})$ is the continuous dual space of $A^{(d - p,d - q)}(\mathcal{X}_{\infty})$.

\vspace{1em}
\noindent An example of an arithmetic surface $\mathcal{X}$ is visualized in the following picture

\begin{figure}[H]
  \centering
  \includegraphics[width=0.55\textwidth,height=0.82\textheight,keepaspectratio]{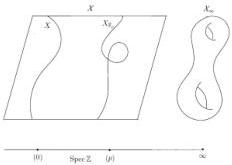}
\end{figure}

\noindent\textbf{Definition 4.3} For an integer $p \geq 0$ let $Z^p(\mathcal{X})$ be the group of cycles $Z$ in $\mathcal{X}$ of codimension $p$. Any cycle $Z = \sum n_i Z_i$, where $n_i \in \mathbb{Z}$, is a formal sum of irreducible cycles $Z_i$, i.e. irreducible closed subschemes of $\mathcal{X}$. An irreducible cycle $Z \in Z^p(\mathcal{X})$ defines a real current $\delta_{Z(\mathbb{C})} \in D^{(p,p)}(\mathcal{X})$ by integration along the smooth part of $Z(\mathbb{C})$. More explicitly, the current $Z(\mathbb{C})$ is defined by

\begin{equation*}
\delta_{Z(\mathbb{C})}(\omega) = \int_{Z(\mathbb{C})_{\mathrm{reg}}} \vartheta^*(\omega)
\end{equation*}

\noindent where $\vartheta : Z(\mathbb{C})_{\mathrm{reg}} \to Z(\mathbb{C})$ is a desingularization of $Z(\mathbb{C})$ along the set of singular points of $Z(\mathbb{C})$. Such desingularization exists due to theorem of Hironaka.

\vspace{1em}
\noindent Now we define Green current.

\noindent\textbf{Definition 4.4} Suppose $Z \in Z^p(\mathcal{X})$. A Green's current for $Z$ is a current $g_Z \in D^{(p - 1,p - 1)}(\mathcal{X})$ such that

\begin{equation*}
dd^c g_Z + \delta_{Z(\mathbb{C})} = [\omega_{g_Z}]
\end{equation*}

\noindent for a smooth form $\omega_{g_Z} \in A^{(p,p)}(\mathcal{X})$.

\vspace{1em}
\noindent Note that for any cycle $Z$ there exists a Green's current for $Z$. Consider the arithmetic variety $X = \mathbb{P}_{\mathbb{Z}}^d = \operatorname{Proj}\mathbb{Z}[x_0, \dots, x_d]$ and the cycle

\begin{equation*}
Z = \left\{ \langle a^{(0)}, x \rangle = \dots = \langle a^{(p - 1)}, x \rangle = 0 \right\}
\end{equation*}

\noindent of codimension $p$, where for $0 \leq i \leq p - 1$ and $x = (x_0, \dots, x_d)$ we set $\langle a^{(i)}, x \rangle = a_0^{(i)}x_0 + \dots + a_d^{(i)}x_d \in \mathbb{Z}[x_0, \dots, x_d]$. The Levine form $g_Z$ associated to $Z$ is defined by

\begin{equation*}
g_Z(x) = -\log \left(\frac{\sum_{i = 0}^{p - 1} | \langle a^{(i)}, x \rangle |^2}{\sum_{i = 0}^d |x_i|^2}\right) \cdot \left(\sum_{j = 0}^{p - 1} \left(dd^c \log \left(\sum_{i = 0}^{p - 1} | \langle a^{(i)}, x \rangle |^2\right)\right)^j \wedge \omega_{FS}^{p - 1 - j}\right)
\end{equation*}

\noindent where $\omega_{FS} = dd^c\left(|x_0|^2 + \dots + |x_d|^2\right)$ is the Fubini-Study form on $\mathbb{P}_{\mathbb{C}}^d$. The Levine form $g_Z$ associated to $Z$ satisfies $dd^c g_Z + \delta_{Z(\mathbb{C})} = [\omega_{FS}]$.

\vspace{1em}
\noindent\textbf{Definition 4.5} Let $\tilde{Z}^p(\mathcal{X})$ be the group of pairs $(Z, g_Z)$, where $Z \in Z^p(\mathcal{X})$ and $g_Z$ is a Green's current for $Z$. Let $\tilde{R}^p(\mathcal{X})$ be the subgroup of $\tilde{Z}^p(\mathcal{X})$ generated by elements of the form

\vspace{0.5em}
\noindent i) $(0, \partial u + \bar{\partial} v)$, where $u \in D^{(p-2,p-1)}(\mathcal{X}_{\infty})$ and $v \in D^{(p-1,p-2)}(\mathcal{X}_{\infty})$

\vspace{0.5em}
\noindent 2) $\left(\mathrm{div}(f), [-\log |f|^2]\right)$ then the $p$-codimensional arithmetic Chow group of $\mathcal{X}$ is defined by

\begin{equation*}
\tilde{CH}^p(\mathcal{X}) := \tilde{Z}^p(\mathcal{X}) / \tilde{R}^p(\mathcal{X})
\end{equation*}

\noindent From [67], there exists an associative, commutative, bilinear pairing,

\begin{equation*}
\tilde{CH}^p(\mathcal{X}) \times \tilde{CH}^q(\mathcal{X}) \to \tilde{CH}^{p + q}(\mathcal{X})_{\mathbb{Q}} := \tilde{CH}^{p + q}(\mathcal{X}) \otimes_{\mathbb{Z}} \mathbb{Q}
\end{equation*}

\begin{equation*}
([Y, g_Y], [Z, g_Z]) \to [Y, g_Y] \cdot [Z, g_Z],
\end{equation*}

\noindent which makes $\bigoplus_{p \geq 0} \tilde{CH}^p(\mathcal{X})_{\mathbb{Q}}$ into a graded ring.

\noindent\textbf{Definition 4.6} A line bundle $\bar{\mathcal{L}} = (\mathcal{L}, h)$ is a $C^\infty$-hermitian line bundle over $\mathcal{X}$ if $\mathcal{L}$ is a line bundle over $\mathcal{X}$ and $h$ is a $C^\infty$-metric over $\mathcal{L}_{\mathbb{C}}$ that is invariant under the complex conjugation, i.e., for $x \in \mathcal{X}(\mathbb{C})$, let $F_\infty : \mathcal{L}_x \to \mathcal{L}_x$ be the isomorphism induced by complex conjugation $F_\infty$. Then $h$ is said to be invariant under the complex conjugation if $h_{\bar{x}}(F_\infty s, F_\infty t) = \overline{h_x(s, t)}$ holds for any $x \in \mathcal{X}(\mathbb{C})$ and $s, t \in \mathcal{L}_x$. Let $s$ and $s_{\mathbb{C}}$ are non-zero rational sections of $\mathcal{L}$ and $\mathcal{L}(\mathbb{C})$ respectively, then we have

\begin{equation*}
\hat{c}_1(\bar{\mathcal{L}}) = (\operatorname{div}(s), [-\log h(s_{\mathbb{C}}, s_{\mathbb{C}})]) \in \tilde{CH}^1(\mathcal{X})
\end{equation*}

\noindent and we call it the arithmetic first Chern class of $\bar{\mathcal{L}}$.

\vspace{1em}
\noindent When $k = \mathbb{C}$, for every separated scheme of finite type $X$ over $\mathbb{C}$ there is an associated complex analytic space $X^{an}$, whose underlying topological space equals the set of complex points $X(\mathbb{C})$. The scheme $X$ is proper over $\mathbb{C}$ if, and only if, $X^{an}$ is compact. Also, $X$ is a non-singular variety over $\mathbb{C}$ if, and only if, $X^{an}$ is a complex analytic manifold.

\vspace{1em}
\noindent Now we give Yau's theorem in Arithmetic setting.

\vspace{1em}
\noindent\textbf{Theorem 4.7} Let $\pi : \mathcal{X} \to \operatorname{Spec}\mathcal{O}_K$ be an arithmetic variety over $\operatorname{Spec}\mathcal{O}_K$, where as usual $\mathcal{O}_K$ denotes the ring of integers of a number field $K$. Suppose that the line bundle $K_{\mathcal{X}}$ is ample line bundle, then, for every $\sigma : K \hookrightarrow \mathbb{C}$, $K_{\mathcal{X}_{\sigma,\mathbb{C}}}$ is ample and there exists a unique Kähler-Einstein metric on $X_{\sigma}(\mathbb{C})$ with constant negative Ricci curvature $-1$.

\vspace{1em}
\noindent For simplicity, throughout this section, we assume that $\mathcal{X}$ is normal. We put $K_{\mathcal{X}^{sm}/C} := \wedge_{\mathcal{O}_{\mathcal{X}^{sm}}}^{n} \Omega_{\mathcal{X}^{sm}/\mathcal{O}_K}$ where $\mathcal{X}^{sm} \subset \mathcal{X}$ denotes the open dense subset of $\mathcal{X}$ where $\pi$ is smooth. Then we further assume, for simplicity, the "$\mathbb{Q}$-Gorenstein condition" i.e. with some $m \in \mathbb{Z}_{>0}$, $(K_{\mathcal{X}^{sm}/C})^{\otimes m}$ extends to an invertible sheaf on whole $\mathcal{X}$. From now on, instead of $K_{\mathcal{X}^{sm}/C}$ we just write $K_{\mathcal{X}}$ for simplicity.

\vspace{1em}
\noindent Now, the \textbf{arithmetic Kodaira dimension} of an arithmetic variety $\mathcal{X}$ can be defined as

\begin{equation*}
\hat{\kappa}(\mathcal{X}) = \limsup_{m \to \infty} \frac{\log \dim H^0(\mathcal{X}, \bar{K}_{\mathcal{X}}^{\otimes m})}{\log m}
\end{equation*}

\noindent where $\dim H^0(\mathcal{X}, \bar{\mathcal{L}}^{\otimes m}) = \dim \{s \in H^0(\mathcal{X}, \mathcal{L}^{\otimes m}) \mid \forall \sigma : K \to \mathbb{C}, \|s\|_{\sigma, \sup} \leq 1\}$ where $\|s\|_{\sigma, \sup} = \sup_{x \in \mathcal{X}(\mathbb{C})} \|s(x)\|$.

\vspace{1em}
\noindent If $K_{\mathcal{X}_0} \geq 0$ then $\mathcal{X}_0$ is a Minimal Model by definition. Now, let the canonical line bundle $K_{\mathcal{X}_0}$ is not be semi-positive, then $\mathcal{X}_0$ can be replaced by sequence of varieties

\noindent $\mathcal{X}_1, \dots, \mathcal{X}_m$ with finitely many birational transformations, i.e., $\mathcal{X}_i$ isomorphic to $\mathcal{X}_0$ outside a codimension 1 subvariety such that

\begin{equation*}
K_{\mathcal{X}_m} \geq 0
\end{equation*}

\noindent and we denote by $\mathcal{X}_{min} = \mathcal{X}_m$ the minimal model of $\mathcal{X}_0$. Hopefully using minimal model program $\mathcal{X}_0$ of non-negative arithmetic Kodaira dimension can be deformed to its minimal model $\mathcal{X}_{min}$ by finitely many birational transformations and we can therefore classify arithmetic projective varieties by classifying their arithmetic minimal models with semi-positive canonical bundle.

\vspace{1em}
\noindent By abundance conjecture, if the minimal model exists, then the canonical line bundle $\mathcal{X}_{min}$ induces a unique holomorphic map

\begin{equation*}
\pi : \mathcal{X}_{min} \to \mathcal{X}_{can}
\end{equation*}

\noindent where $\mathcal{X}_{can}$ is the unique canonical model of $\mathcal{X}_{min}$. The canonical model completely determined by arithmetic variety $\mathcal{X}$ as follows,

\begin{equation*}
\mathcal{X}_{can} = \operatorname{Proj} \left( \bigoplus_{m \geq 0} H^0(\mathcal{X}, K_{\mathcal{X}_{min}}^m) \right)
\end{equation*}

\noindent So combining MMP and Abundance conjecture, we directly get

\begin{equation*}
\mathcal{X} \to \mathcal{X}_{min} \to \mathcal{X}_{can}
\end{equation*}

\noindent which $\mathcal{X}$ here might not be birationally equivalent to $\mathcal{X}_{can}$ and dimension of $\mathcal{X}_{can}$ is smaller than of dimension $\mathcal{X}$. Note also that, Arithmetic Minimal model is not necessary unique but the arithmetic canonical model $\mathcal{X}_{can}$ is unique.

\vspace{1em}
\noindent and on

\begin{equation*}
\pi : \mathcal{X} \to \mathcal{X}_{can}
\end{equation*}

\noindent when arithmetic Kodaira dimension is positive, we have

\begin{equation*}
\operatorname{Ric}(\omega_{can}) = -\omega_{can} + \pi^* \omega_{WP}
\end{equation*}

\noindent where $\omega_{can}$ is a canonical metric on $\mathcal{X}_{can}(\mathbb{C})$

\vspace{1em}
\noindent By \textbf{Arithmetic Song-Tian program} it turns out that the normalized Kähler Ricci flow

\begin{equation*}
\frac{\partial \omega_t}{\partial t} = -\operatorname{Ric}(\omega_t) + a \omega_t
\end{equation*}

\noindent (in the sense that $\omega_t$ are Kähler currents of $\mathcal{X}_t(\mathbb{C})$) doing exactly same thing to replace $\mathcal{X}$ by its arithmetic minimal model by using finitely many arithmetic surgeries and then deform arithmetic minimal model to arithmetic canonical model such that the limiting of arithmetic canonical model is coupled with generalized Kähler Einstein metric twisted with Weil-Petersson metric gives canonical metric for $\mathcal{X}(\mathbb{C})$.

\vspace{1em}
\noindent In fact when the arithmetic Kodaira dimension is semi-positive then the canonical metric $\omega_{can}$ on arithmetic variety is a metric which is attached to arithmetic canonical model $X_{can}$.

\vspace{1em}
\noindent Arithmetic Song-Tian program via Arithmetic Minimal Model program is that if $X_0$ be a projective arithmetic variety with a smooth Kähler metric $\omega_0$, we apply the Kähler Ricci flow with initial data $(\mathcal{X}_0, \omega_0)$, then there exists $0 < T_1 < T_2 \dots < T_{m+1} \leq \infty$ for some $m \in \mathbb{N}$, such that

\begin{equation*}
(\mathcal{X}_0, \omega_0) \xrightarrow{t \to T_1} (\mathcal{X}_1, g_1) \xrightarrow{t \to T_2} \dots \xrightarrow{t \to T_m} (\mathcal{X}_m, \omega_m)
\end{equation*}

\noindent after finitely many surgeries in Gromov-Hausdorff topology and either $\dim_{\mathbb{C}} \mathcal{X}_m < n$, or $\mathcal{X}_m = \mathcal{X}_{min}$ if not collapsing. In the case $\dim_{\mathbb{C}} \mathcal{X}_m < n$, then $\mathcal{X}_m$ admits Fano fibration which is a morphism of varieties whose general fiber is a Fano variety of positive dimension in other words has ample anti-canonical bundle.

\subsection*{References}

\begin{enumerate}[label={[\arabic*]}, leftmargin=*]
\item Tien-Cuong Dinh and Nessim Sibony, Introduction to the theory of currents, Preprint, 2005
\item Jean-Pierre Demailly, Complex Analytic and Differential Geometry, preprint, 2012
\item Tien-Cuong Dinh, Nessim Sibony, Pull-back of currents by holomorphic maps, manuscripta mathematica, 2007, Volume 123, Issue 3, pp 357-371
\item Claire Voisin, Recent Progresses in Kähler and Complex Algebraic Geometry, 4ECM Stockholm 2004,
\item Jian Song, Riemannian geometry of Kähler-Einstein currents II: an analytic proof of Kawamata's base point free theorem, arXiv:1409.8374
\end{enumerate}

\begin{list}{}{\leftmargin=2.5em\labelwidth=2em\labelsep=0.5em\itemsep=0.8ex\parsep=0pt}
\item[[6]] X. Chen, S. Donaldson, S. Sun, K\"{a}hler-Einstein metrics on Fano manifolds. I, Approximation of metrics with cone singularities, Journal of the American Mathematical Society 28 (1), 183-197

\item[[7]] X. Chen, S. Donaldson, S. Sun, K\"{a}hler-Einstein metrics on Fano manifolds. II: Limits with cone angle less than $2\pi$, Journal of the American Mathematical Society 28 (1), 199-234

\item[[8]] X. Chen, S. Donaldson, S. Sun, K\"{a}hler-Einstein metrics on Fano manifolds. III, Limits as cone angle approaches $2\pi$ and completion of the main proof, Journal of the American Mathematical Society 28 (1), 235-278

\item[[9]] G. Tian. K-stability and K\"{a}hler-Einstein metrics. arXiv:1211.4669

\item[[10]] Jian Song; Gang Tian, The K\"{a}hler-Ricci flow on surfaces of positive Kodaira dimension, Inventiones mathematicae. 170 (2007), no. 3, 609-653.

\item[[11]] Jian Song; Gang Tian, Canonical measures and K\"{a}hler-Ricci flow, J. Amer. Math. Soc. 25 (2012), no. 2, 303-353,

\item[[12]] S. Boucksom and H. Tsuji, Semipositivity of relative canonical bundles via K\"{a}hler-Ricci flows, preprint.

\item[[13]] Hajime Tsuji, Canonical measures and the dynamical systems of Bergman kernels, arXiv:0805.1829

\item[[14]] Robert J. Berman, K\"{a}hler-Einstein metrics, canonical random point processes and birational geometry, arXiv:1307.3634v1

\item[[15]] S. T. Yau, On the Ricci curvature of a compact K\"{a}hler manifold and the complex Monge-Ampere equation. I. Comm. Pure Appl. Math.31(1978), no. 3, 339-411

\item[[16]] H.D.Cao: Deformation of K\"{a}hler metrics to K\"{a}hler-Einstein metrics on compact K\"{a}hler manifolds. Invent. Math.81(1985), no. 2, 359-372.

\item[[17]] T. Aubin, Nonlinear Analysis on Manifolds, Monge-Ampere Equations, ISBN 0-387-90704-1

\item[[18]] G. Perelman, The entropy formula for the Ricci flow and its geometric applications. arXiv:math.DG/0211159.

\item[[19]] G. Perelman, Ricci flow with surgery on three-manifolds. arXiv: math.DG/ 0303109

\item[[20]] Tian, G. and Zhu, X., Uniqueness of K\"{a}hler-Ricci solitons, Acta Math. 184 (2000), no. 2, 271-305

\item[[21]] G. Tian. On K\"{a}hler-Einstein metrics on certain K\"{a}hler manifolds with $C_1(M) > 0$. Invent. Math., 89(2):225\{246, 1987.

\item[[22]] Futaki, A., An obstruction to the existence of Einstein-K\"{a}hler metrics. Invent. Math., 73,437-443 (1983)

\item[[23]] Jian Song; Gang Tian, The K\"{a}hler-Ricci flow on surfaces of positive Kodaira dimension, Inventiones mathematicae. 170 (2007), no. 3, 609-653.
\end{list}

\begin{list}{}{\leftmargin=2.5em\labelwidth=2em\labelsep=0.5em\itemsep=0.8ex\parsep=0pt}
\item[[24]] Jian Song; Gang Tian, Canonical measures and K\"{a}hler-Ricci flow, J. Amer. Math. Soc. 25 (2012), no. 2, 303-353,

\item[[25]] Tian, Gang; Yau, Shing-Tung, Complete K\"{a}hler manifolds with zero Ricci curvature, I, Amer. Math. Soc.(1990) 3 (3): 579-609

\item[[26]] Jian Song, Gang Tian, The K\"{a}hler-Ricci flow through singularities, arXiv:0909.4898

\item[[27]] Caucher Birkar, Paolo Cascini, Christopher D. Hacon, James McKernan, Existence of minimal models for varieties of log general type, J. Am. Math. Soc. 23(2), 405-468 (2010)

\item[[28]] Y-T. Siu, A General Non-Vanishing Theorem and an Analytic Proof of the Finite Generation of the Canonical Ring. arXiv:math.AG/0610740

\item[[29]] A.M.Nadel, Multiplier ideal sheaves and existence of K\"{a}hler-Einstein metrics of positive scalar curvature, Proc Natl Acad Sci U S A. 1989, Oct; 86(19):7299-300.

\item[[30]] A. M. Nadel, Multiplier ideal Sheaves and K\"{a}hler-Einstein metrics of positive scalar curvature, Ann. of Math.132(1990), 549-596

\item[[31]] H. Skoda. Sous-ensembles analytiques d'ordre fini ou infini dans $\mathbb{C}^n$. Bulletin de la Societe Mathematique de France, 100:353-408

\item[[32]] P. Eyssidieux, V. Guedj, A. Zeriahi: Singular K\"{a}hler-Einstein metrics, J. Amer. Math. Soc. 22 (2009), no. 3, 607\{639.

\item[[33]] E. Viehweg, Weak positivity and the additivity of the Kodaira dimension for certain fibre spaces. Algebraic varieties and analytic varieties (Tokyo, 1981), 329-353, Adv. Stud. Pure Math., 1, North-Holland, Amsterdam, 1983

\item[[34]] J. Fine. Constant scalar curvature K\"{a}hler metrics on fibred complex surfaces. J. Differential Geom., 68(3):397-432, 2004

\item[[35]] Young-Jun Choi, Semi-positivity of fiberwise Ricci-flat metrics on Calabi-Yau fibrations, arXiv:1508.00323

\item[[36]] Robert J. Berman, Relative K\"{a}hler-Ricci flows and their quantization, Analysis and PDE, Vol. 6 (2013), No. 1, 131\{180

\item[[37]] S.Boucksom, Higher dimensional Zariski decompositions, Preprint, arXiv:math/0204336

\item[[38]] Georg Schumacher, Moduli of framed manifolds, Invent. math. 134, 229-249 (1998)

\item[[39]] Ken-Ichi Yoshikawa, On the boundary behavior of the curvature of $L^{2}$-metrics, \url{http://arxiv.org/abs/1007.2836}

\item[[40]] Georg Schumacher and Hajime Tsuji, Quasi-projectivity of moduli spaces of polarized varieties, Annals of Mathematics, 159(2004), 597-639

\item[[41]] Valentino Tosatti, Yuguang Zhang, Infinite time singularities of the K\"{a}hler-Ricci flow, Geom. Topol. 19 (2015), no.5, 2925-2948.

\item[[42]] Y. Kawamata, Characterization of abelian varieties, Compositio Math. 43(1981), no. 2, 253\{276.
\end{list}

\begin{list}{}{\leftmargin=2.5em\labelwidth=2em\labelsep=0.5em\itemsep=0.8ex\parsep=0pt}
\item[[43]] S. Y. Arakelov, An intersection theory for divisors on an arithmetic surface, Math. USSR Izvestija 8 (1974), 1167-1180.

\item[[44]] G. Faltings, Calculus on arithmetic surfaces, Ann. of Math. 119 (1984), 387-424.

\item[[45]] P. Deligne: Le d\'{e}terminant de la cohomologie. Current trends in arithmetical algebraic geometry (Arcata, Calif., 1985). Contemp. Math. 67. Amer. Math. Soc., Providence, RI, 1987, 93177. MR 902592 (89b:32038)

\item[[46]] C. Soul\'{e}: Lectures on Arakelov geometry. Cambridge Studies in Advanced Mathematics 33. Cambridge University Press, Cambridge, 1992, With the collaboration of D. Abramovich, J.-F. Burnol and J. Kramer. MR 1208731 (94e:14031)

\item[[47]] Georg Schumacher, Moduli of framed manifolds, Invent. math. 134, 229-249 (1998)

\item[[48]] Eleonora Di Nezza, Uniqueness and short time regularity of the weak K\"{a}hler-Ricci flow, arXiv:1411.7958

\item[[49]] Yujiro Kawamata, Subadjunction of log canonical divisors, II, Amer. J. Math. 120 (1998), 893\{899

\item[[50]] K.-I. Yoshikawa, Singularities and analytic torsion, arXiv:1007:2835

\item[[51]] Chen, X. X.; Tian, G. Geometry of K\"{a}hler metrics and foliations by holomorphic discs. Publ. Math. Inst. Hautes \'{E}tudes Sci. No. 107 (2008),

\item[[52]] Hassan Jolany, Complex Monge-Ampere foliation via Song-Tian-Yau-Vafa foliation. In preparation

\item[[53]] Morales, Marcel. "Resolution of quasi-homogeneous singularities and plurigenera." Compositio Mathematica 64.3 (1987): 311-327.

\item[[54]] Hubert Flenner, Mikhail Zaidenberg, Log\{canonical forms and log canonical singularities, Math. Nachr. 254\{255, 107-125 (2003)

\item[[55]] D. Mumford. Hirzebruch's proportionality theorem in the noncompact case. Invent.Math., 42:239\{272, 1977

\item[[56]] Andrey Todorov, Weil-Petersson Volumes of the Moduli Spaces of CY Manifolds, arXiv:hep-th/0408033

\item[[57]] Hacon (http://mathoverflow.net/users/19369/hacon), A Decomposition for Iitaka fibration, URL (version: 2016-08-02): http://mathoverflow.net/q/245682

\item[[58]] Eckart Viehweg, De-Qi Zhang, Effective Iitaka fibrations, arXiv:0707.4287

\item[[59]] Y.Odaka, Canonical K\"{a}hler metrics and Arithmetics -- Generalising Faltings heights, arXiv:1508.07716

\item[[60]] Yuji Odaka, The GIT stability of polarized varieties via discrepancy, Pages 645-661, Volume 177 (2013), Issue 2, Annals of Mathematics

\item[[61]] Yum Tong Siu, Analyticity of sets associated to Lelong numbers and the extension of closed positive currents. Invent. Math., 27:53-156, 1974

\item[[62]] Tsuji, H., Global generations of adjoint bundles. Nagoya Math. J. Vol. 142 (1996),5-16
\end{list}

\begin{list}{}{\leftmargin=2.5em\labelwidth=2em\labelsep=0.5em\itemsep=0.8ex\parsep=0pt}
\item[[63]] J. Kollar and N. I. Shepherd-Barron : Threefolds and deformations of surface singularities, \textit{Invent. Math.} 91 (1988), no. 2, 299-338.

\item[[64]] K\"{a}hler-Einstein Metrics on Stable Varieties and log Canonical Pairs, \textit{Geometric and Functional Analysis} December 2014, Volume 24, Issue 6, pp 1683\{1730

\item[[65]] G. Tian, Degeneration of K\"{a}hler-Einstein Manifolds, I, \textit{Proc. Symp. in Pure Math.} 54 (1993), 595-609.

\item[[66]] D. AbramovichK. Karu, Weak semistable reduction in characteristic 0, \textit{Inventiones mathematicae}, February 2000, Volume 139, Issue 2, pp 241-273

\item[[67]] C.-L. Wang, Quasi-Hodge metrics and canonical singularities. \textit{Mathematical Research Letter}, 10:57-70, 2003

\item[[68]] Candelas, Philip; H\"{u}bsch, Tristan; Schimmrigk, Rolf, Relation between the Weil-Petersson and Zamolodchikov metrics. \textit{Nuclear Phys.} B 329 (1990), no. 3, 583-590.

\item[[69]] Tian, Gang, Smoothness of the universal deformation space of compact Calabi-Yau manifolds and its Petersson-Weil metric. \textit{Mathematical aspects of string theory} (San Diego, Calif., 1986), 629\{646, \textit{Adv. Ser. Math. Phys.}, 1, World Sci. Publishing, Singapore, 1987

\item[[70]] Georg Schumacher, Moduli of framed manifolds, \textit{Invent. math.} 134, 229-249(1998)

\item[[71]] Thalia Jeffres, Rafe Mazzeo, Yanir A. Rubinstein, K\"{a}hler\{Einstein metrics with edge singularities, \textit{Annals of Mathematics}, 95-176, Volume 183 (2016), Issue 1

\item[[72]] Osamu Fujino and Shigefumi Mori, A Canonical Bundle Formula, \textit{J. Differential Geom. Volume 56, Number 1} (2000), 167-188.

\item[[73]] Osamu Fujino, Semi-stable minimal model program for varieties with trivial canonical divisor, \textit{Proc. Japan Acad. Ser. A Math. Sci. Volume 87, Number 3} (2011), 25-30

\item[[74]] Liangming Shen, PhD thesis, Princeton University, 2015

\item[[75]] Y. Ruan, G. Tian, A mathematical theory of quantum cohomology, \textit{J. Differential Geom.} 42 (2) (1995) 259-367.

\item[[76]] Richard S. Hamilton, Eternal solutions to the Ricci flow, \textit{J. Differential Geom. Volume 38, Number 1} (1993), 1-11.

\item[[77]] R. Bamler, S. Brendle, ``A comparison principle for solutions to the Ricci flow'', \textit{Math. Res. Letters}, 22(4) (2015):983-988
\end{list}

\vspace{2cm}
\begin{flushleft}
\textit{Laboratoire de Math\'{e}matiques Paul Painlev\'{e},\\
Universit\'{e} des Sciences et Technologies de Lille,\\
Lille, France}

\vspace{0.5cm}
\texttt{jolany.hassan.math@gmail.com}
\end{flushleft}

\end{document}